# Developing a 21st Century Global Library for Mathematics Research

Committee on Planning a Global Library of the Mathematical Sciences

Board on Mathematical Sciences and Their Applications

Division on Engineering and Physical Sciences

NATIONAL RESEARCH COUNCIL
*OF THE NATIONAL ACADEMIES*

**THE NATIONAL ACADEMIES PRESS     500 Fifth Street, NW     Washington, DC 20001**

NOTICE: The project that is the subject of this report was approved by the Governing Board of the National Research Council, whose members are drawn from the councils of the National Academy of Sciences, the National Academy of Engineering, and the Institute of Medicine. The members of the committee responsible for the report were chosen for their special competences and with regard for appropriate balance.

This project was supported by the Alfred P. Sloan Foundation under grant number 2011-10-28. Any opinions, findings, conclusions, or recommendations expressed in this publication are those of the author(s) and do not necessarily reflect the views of the organization that provided support for the project.

International Standard Book Number 13:   978-0-309-29848-3
International Standard Book Number 10:   0-309-29848-2

Additional copies of this report are available from the National Academies Press, 500 Fifth Street, NW, Keck 360, Washington, DC 20001; (800) 624-6242 or (202) 334-3313; http://www.nap.edu.

Suggested citation: National Research Council. 2014. *Developing a 21st Century Global Library for Mathematics Research*. Washington, D.C.: The National Academies Press.

# THE NATIONAL ACADEMIES
*Advisers to the Nation on Science, Engineering, and Medicine*

# Acknowledgments

This report has been reviewed in draft form by individuals chosen for their diverse perspectives and technical expertise, in accordance with procedures approved by the National Research Council's Report Review Committee. The purpose of this independent review is to provide candid and critical comments that will assist the institution in making its published report as sound as possible and to ensure that the report meets institutional standards for objectivity, evidence, and responsiveness to the study charge. The review comments and draft manuscript remain confidential to protect the integrity of the deliberative process. The committee wishes to thank the following individuals for their review of this report:

Sara Billey, University of Washington
Thierry Bouche, Cellule MathDoc and Institut Fourier, Université de Grenoble
François G. Dorais, MathOverflow and Dartmouth College
Robion Kirby, University of California, Berkeley
Donald McClure, American Mathematical Society
Jason Rute, Pennsylvania State University
Terence Tao, University of California, Los Angeles
Eva Tardos, Cornell University
Heinz Weinheimer, Springer

Although the reviewers listed above have provided many constructive comments and suggestions, they were not asked to endorse the conclusions or recommendations nor did they see the final draft of the report before its release. The review of this report was overseen by C. David Levermore, University of Maryland, College Park. Appointed by the National Research Council, he was responsible for making certain that an independent examination of this report was carried out in accordance with institutional procedures and that all review comments were carefully considered. Responsibility for the final content of this report rests entirely with the authoring committee and the institution.

The committee also acknowledges the valuable contribution of the following individuals, who provided input at the meetings on which this report is based or by other means:

# Contents

APPENDIXES

# Summary

Like most areas of scholarship, mathematics is a cumulative discipline: new research is reliant on well-organized and well-curated literature. Because of the precise definitions and structures within mathematics, today's information technologies and machine learning tools provide an opportunity to further organize and enhance discoverability of the mathematics literature in new ways, with the potential to significantly facilitate mathematics research and learning. Opportunities exist to enhance discoverability directly via new technologies and also by using technology to capture important interactions between mathematicians and the literature for later sharing and reuse.

In most scientific disciplines, including mathematics, Web-based access to digital resources representing the disciplinary literature is now mature and quite effective. Through a mixture of open and proprietary tools, mathematicians are able to search the enormous and very rapidly growing literature using attributes such as subjects, titles, authors, dates, and keywords; they can follow chains of citations among works backward and forward in time. While much information is contained in individual items in the mathematical literature, a greater amount of information is represented by the way they are linked. This is not just via references but through the interrelation of concepts, insights, and techniques as they are developed, refined, and spread from one mathematical discipline to another. For example, if mathematicians were able to search the literature for instances where a specific equation was used or solved, it would allow them to consider alternative approaches toward solving their own research questions. This search capability could be facilitated through the use of a database of machine-generated and human-cultivated information about the mathematical literature and allow for a variety of other capabilities to be built.

This report discusses how information about what the mathematical literature contains can be formalized and made easier to express, encode, and explore. Many of the tools necessary to make this information system a reality will require much more than indexing and will instead depend on community input paired with machine learning, where mathematicians' expertise can fill the gaps of automatization. The Committee on Planning a Global Library of the Mathematical Sciences proposes the establishment of an organization; the development of a set of platforms, tools, and services; the deployment of an ongoing applied research program to complement the development work; and the mobilization and coordination of the mathematical community to take the first steps toward these capabilities.

Mathematics today has the opportunity to expand and redefine the way in which mathematical knowledge is represented and used, the character of the mathematical literature and how it evolves, and the way that mathematicians interact with this collection of knowledge. This new relationship with the literature and the mathematical knowledge corpus goes beyond new forms of access and analytical tools; it must also include the tools and services to accommodate the creation, sharing, and curation of new kinds of knowledge structures.

To be clear, what the committee proposes builds on the extensive work done by many dedicated individuals under the rubric of the World Digital Mathematical Library,[1] as well as many other community initiatives.[2] Comparing desired capabilities going forward with what has been achieved by these efforts to date, the committee concludes that there is little value in new large-scale retrospective digitization efforts or further aggregations of mathematical science publications (both traditional journal articles and newer preprint, blog, video, and similar resources) beyond the federation of distributed repositories already achieved through existing search services. Nor is another bibliographically based secondary indexing service needed at this time. Necessary incremental improvements will likely continue to occur in these areas, but they do not require an initiative on the scale of what is being called for in this report.

The real opportunity is in offering mathematicians new and more direct ways to discover and interact with mathematical objects and mathematical knowledge through the Web. The committee's consensus is that by some combination of machine learning methods and community-based editorial effort, a significantly greater portion of the information and knowledge in the global mathematical corpus could be made available to researchers as linked open data[3] through a central organizational entity—referred to in this report as the Digital Mathematics Library (DML).

---

[1] The World Digital Mathematics Library rubric has been used by a variety of organizations for many distinct projects. A history of many of these efforts and the current state-of-the-art can be found on the wiki page from the International Mathematics Union's Digital Mathematics Workshop in June 2012, http://ada00.math.uni-bielefeld.de/mediawiki-1.18.1/index.php/.

[2] Examples include the *Encyclopedia of Integer Sequences*, the NIST Digital Library of Mathematical Functions, and the *Guide to Available Mathematical Software*.

The DML would aggregate and make available collections of ontologies, links, and other information created and maintained by human contributors, curators, and specialized machine agents, with significant editorial input from the mathematical community. The DML would enable functionalities and services over the aggregated mathematical information that go well beyond simply making publications available, to include capabilities for annotating, searching, browsing, navigating, linking, computing, and visualizing both copyrighted and openly licensed content. While the DML would store modest amounts of new knowledge structures and indices, it would not generally replicate mathematical literature stored elsewhere. Instead, it would strive to represent the mathematical *knowledge* presented within a publication and illustrate how it is connected with other resources.

While the committee believes that the DML could begin development soon, it notes that this work would need to be complemented by an ongoing research program to fill in gaps, improve quality and performance, increase the robustness of available technologies, and increase the automation of processes that still rely heavily on human intervention.

The DML would facilitate discovery of and interaction with mathematical information from diverse sources with varying levels of copyright. The committee envisions the DML as a growing corpus of public-domain and openly licensed mathematical information, Web services, and software agents, which would coexist with present mathematical publishing and indexing services for the foreseeable future.

A key early issue for the DML organization is how to establish constructive and effective partnerships with existing publishers, Web services, and other resources, both those specific to mathematics and those serving the much broader scholarly community. Some of these partnerships might be challenging because of copyright concerns. However, establishing fruitful partnerships is essential to the success of the DML. While the DML would sometimes provide services and functional features that overlap with existing services and tools provided by both commercial and not-for-profit entities, the committee suggests partnering with current service providers whenever possible rather than replicating capabilities of existing resources.

---

[3] Broadly defined, linked open data are structured data that are published in such a way that makes it easy to interlink them with other data, therefore making it possible to connect them with information from multiple sources. These connected data can provide a user with a more meaningful query of a subject by consolidating relevant information from a variety of places—e.g., in different research papers—and pulling out specific components that the user might be particularly interested in.

For example in MathOverflow,[4] a question-and-answer website for research mathematicians, research articles and papers are often referenced in answers given. While the DML would not want to replicate the interface and social networking features of MathOverflow, it would be wholly appropriate for the DML to instigate and participate in a multi-party collaboration with MathOverflow and publishers of research mathematics to automatically capture citations entered in MathOverflow answers and republish them as linked open data annotations. In this scenario, the DML could help broker standard practices for interoperability and help maintain the software agents and annotation repositories that would allow publishers to make mathematicians coming to their websites aware of MathOverflow discussions potentially relevant to the papers they are viewing. The converse could also be supported. Posts on MathOverflow could be automatically annotated when errata or other commentary is added to the publisher's website for an article mentioned in the MathOverflow post. This illustrates the potential for chains of annotations as a new mode of scholarly discourse (Sukovic, 2008). To visualize how an annotation chain might come about, begin by assuming that a post in MathOverflow referencing a particular article is automatically added as an annotation to this article on the publisher's website. A subsequent reply to this annotation made by a reader of the publisher website is then automatically added to the thread on MathOverflow. A new reply subsequently added to the thread on MathOverflow is then automatically added as a further annotation on the publisher's website, and so on. This would allow users of two disparate services—i.e., one scholar using MathOverflow and the other using only the publisher's website—to nonetheless carry on a substantive discourse about published mathematics research in spite of the fact that each is using a different utility to access the publication being discussed.

Similarly, MathSciNet and Zentralblatt Math (zbMath) already classify research papers according to the Mathematics Subject Classification (MSC)[5] schedule. The DML would not want to replicate this indexing. However, it might be beneficial for the DML to provide complementary indexing on other dimensions—e.g., by the occurrence in articles of well-known special functions (hierarchies of which are maintained by the National Institute for Standards and Technology (NIST)[6] and by Wolfram

---

[4] MathOverflow, http://mathoverflow.net/, accessed January 16, 2014.

[5] American Mathematical Society, 2010 Mathematics Subject Classification, http://www.ams.org/mathscinet/msc/msc2010.html, accessed January 16, 2014.

[6] NIST, Digital Library of Mathematical Functions, Version 1.0.6, release date May 6, 2013, http://dlmf.nist.gov/.

Research[7]). Used in concert, one could then envision a collaboratively built interface that allows refinement of an initial MSC search via attributes such as which special functions are used in the articles that appear in the results from the MSC search.

Such partnerships and collaborations are essential. It is vital that users see a well-integrated interface that incorporates both the DML services and commercial services for those affiliated with institutions that have access to the commercial services. The committee envisions the resources, services, and tools offered by the DML as coexisting with, and often enhancing, the offerings from existing players in the mathematical information landscape. The committee hopes that relevant organizations will contribute to the work of the DML in various ways, such as by providing financial support, allowing appropriate access to their content and services, or by participating in the collaborative development, with the shared goal of enhancing the value of the mathematics literature. Building these partnerships would likely require significant negotiations and collaborations, and the DML organization would have to allocate much time and effort to their planning and execution.

The biggest challenge, however, will be in establishing the technical, organizational, and community-coordinating capabilities to deliver on the construction of the resources, services, and tools described earlier in this summary and then planning and implementing the development and deployment of the necessary systems. Some of the technologies required to build the requisite tools and services do not exist today or are not sufficiently mature. The committee sees the DML as having a minimal direct research role; rather, the committee believes that the establishment of the DML needs to be complemented by a long-term (5 to 10 years) commitment to a focused and applied research program that would encompass both needed technology, tools, and services and (to a lesser extent) independent research to understand how the DML is being used and how well it is working. Ideally, the commitment to fund this program could come in parallel with the commitment for the initial funding for the DML itself (whether from one or multiple sources). These research programs need to be well connected to the work of the DML. This could be achieved either by ensuring that the DML is deeply involved in the development of the calls for proposals and the subsequent proposal evaluation or by actually placing the DML in the role of a re-granting organization (although the committee sees some potential bureaucratic complications with the latter option).

---

[7] Wolfram Research, Inc., The Wolfram Functions Site, http://functions.wolfram.com/, accessed January 16, 2014.

## ORGANIZATION AND RESOURCES NEEDED

The committee's vision of an incremental development of the DML starts with the creation of a small nonprofit organization, referred to here as the DML organization. The DML organization will need a small and dedicated paid staff, including a well-respected mathematician in a senior role, to ensure its development and growth. Other staffing needs may become necessary as the needs and status of the DML evolve, although much of the software development and operations could be contracted out. Ideally, the DML would be attached to and draw support from some host institution (a university, a research laboratory, or other organization) in order to facilitate sharing of services and to reduce overhead. The DML organization could be governed ultimately by the mathematical sciences community through organizations such as the International Mathematical Union and, thence, through their member organizations.

The first and foremost challenge that the DML will face is finding a set of primary funding sources that could support its initial development and early operations (a period of between 5 and 10 years). It is the committee's hope that the DML would become a self-sustaining entity once some of its key capabilities are established and a potential sustainable business model is chosen from among options.[8]

For the first few years, perhaps the best approach would be to split operational governance from high-level, longer-term policy governance, because these two tasks will be quite distinct. Both in the short and the longer term, appropriate connections are needed between funding and revenue sources and governance, and these connections may well need to shift over time. Particularly in the early days, a light and agile governance mechanism is crucial. Upon launching the DML effort, there would likely be a coalition of partners with a commitment to the DML concept.

## CONCLUSION

Like other scientific disciplines, mathematics is now completing a complex multi-decade transition from print to a digital system that closely emulates print for authors and readers. The mathematics community is thus at an inflection point where it has the opportunity to think about how its collective knowledge base is going to be constructed, used, structured, managed, curated, and contributed to in the digital world and how that knowledge base will be related to the existing literature corpus, to authoring practices in the future, and to the social and community practices of doing and learning mathematics. Colleagues in other disciplines—astronomy, molecular biology, genomics, chemistry—are in many cases well advanced in formulating their own disciplinary-specific answers that take into account disciplinary practices (such as the mix of experimental, observational, theoretical, and computational approaches) and the conceptual models that underlie disciplinary thinking.

---

[8] There are many lessons on sustainability to draw upon, including experiences with digital libraries (such as arXiv) and open or community source software as well as work on research data curation.

Mathematics is unusual in many ways; it maintains a healthy and constructive relationship with its past, as documented in the literature of the field going back hundreds of years, and some of its literature has a long "shelf life." The committee believes that investments in refreshing and restructuring the corpus of mathematical literature and abstracting it into a knowledge base for future centuries is a valid and sound investment in the future of mathematical scholarship. The DML proposed in this report provides a platform and a context to achieve this and also offers a critical point of focus for the mathematical community in a genuinely digital environment to engage in discussions about the creation, curation, and management of mathematical knowledge.

# 1

# Introduction

## OVERVIEW

Mathematics is facing a pivotal junction where it can either continue to utilize digital mathematics literature in ways similar to traditional printed literature, or it can take advantage of new and developing technology to enable new ways of advancing knowledge. This report details how information contained in individual items within the literature could be readily extracted and linked to create a comprehensive digital mathematics information resource that is more than the sum of its contributing publications. That resource can serve as a platform and focal point for further development of the mathematical knowledge base.

This new system, referred to throughout the report as the Digital Mathematics Library (DML), could support a wide variety of new functionalities and services over aggregated mathematical information, including dramatically improved capabilities for searching, browsing, navigating, linking, computing, visualizing, and analyzing the literature.

## STUDY DEFINITION AND SCOPE AND THE COMMITTEE'S APPROACH

The Alfred P. Sloan Foundation commissioned this study and charged the committee to:

- Evaluate the potential value of a virtual global library of mathematical science publications;

- Assuming that a stable context for sharing copyrighted information has been achieved, assess the remaining issues to be addressed in setting up such a library;
- Identify a range of desired capabilities of such a library; and
- Characterize resource needs.

While a traditional library is perhaps the oldest formal information resource available, the manifestation of libraries has evolved dramatically over the past few decades. In many cases within mathematics, as for other fields of scholarship, buildings housing paper publications have given way to online collections of downloadable documents. While this increased access is not perfect—not all material is readily available to all researchers, and search tools vary from site to site—widespread digitization has made it easier for many to access the mathematical literature. Overall, a much greater proportion of the mathematical literature is available to more people than at any time before. The research libraries, scholarly societies, and other players that curate and steward this material continue to grapple with issues, such as long-term preservation of digital materials, but it is fair to say there exists a fairly comprehensive, distributed "digital library" for mathematics offering a much improved but not fundamentally different version of what existed in the time of printed books and journals.

The committee has thus taken the term *library* in its charge to mean a system that accumulates and shares knowledge, rather than the more traditional library that houses documents, either digital or physical. The committee's focus has been on *functionality* that can meet the needs of mathematicians facing a rapidly expanding and diversifying knowledge base. The committee has largely ignored traditional issues of assembling and stewardship of those collections, which are being handled well, for the most part, by the existing distributed digital library.

The committee envisions its target digital library users to be working research mathematicians and advanced graduate students beginning their research careers throughout the world (hence the word *global*). The library discussed does not specifically target students below the advanced graduate student level or researchers outside of mathematics, although both sets would likely constitute some of the library's user base. Having a clear understanding of the target user base directly impacts the types of content the library targets and the types of services it provides. The committee also believes that the disciplinary scope of the mathematics that this library could provide is best left undefined for now. Mathematics and the mathematical sciences have diffuse boundaries, and this committee takes no stance on where appropriate content lies. However, this is an issue that will have to be addressed by either a future management organization or the community of users.

The committee believes that there is much room for innovation and progress in the mainstream mathematical information services. To determine which potential areas for innovation are of the most interest to the mathematics community, the committee held three meetings where it heard from outside presenters on issues relevant to mathematics (November 27-28, 2012; February 19-20, 2013; and May 30-31, 2013—agendas for these meetings can be found in Appendix A) and two public data-gathering sessions (at the University of Minnesota on May 6, 2013, and at Northwestern University on May 30, 2013), posted questions on two mathematics discussion forums (MathOverflow[1] and Math 2.0[2]), and wrote a guest entry on Professor Terry Tao's mathematics blog.[3] The committee also referred to the information shared at the World Digital Mathematics Library workshop held by the International Mathematical Union (IMU) on June 1-3, 2012.[4]

The committee made an assessment of what computers can do today, what computers can help mathematicians to do, and how rapidly these capabilities are likely to grow, if provided with some ongoing focused research funding. The committee's consensus is that by some combination of machine learning methods and community-based editorial effort, a significant portion of the information and knowledge in the global mathematical corpus could be made available to researchers as linked open data. Broadly defined, linked open data are structured data that are published in such a way that makes it easy to interlink them with other data, thereby making it possible to connect them with information from multiple sources. This connected data can provide a user with a more meaningful query of a subject by consolidating relevant information from a variety of places (e.g., in different research papers) and pulling out specific components that the user might be particularly interested in. The committee envisions that much of the existing mathematical information can be provided as linked open data through a central organizational entity—referred to in this report as the DML. It should be noted that linked open data are not the only way that this can be accomplished, but they are essentially today's standard for ontologies and other important representations. The committee believes that the DML should make use of current best practices rather than trying to develop some other alternative, whenever possible.

---

[1] I. Daubechies, "Math Annotate Platform?," MathOverflow (question and answer site), February 18, 2013, http://mathoverflow.net/questions/122125/math-annotate-platform.

[2] I. Daubechies, "Math Annotate Platform?," Math2.0 (discussion forum), February 18, 2013, http://publishing.mathforge.org/discussion/163/.

[3] I. Daubechies, "Planning for the World Digital Mathematical Library," *What's New* (blog by Terence Tao), daily archive for May 8, 2013, http://terrytao.wordpress.com/2013/05/08/.

[4] Many of the materials presented at the International Mathematics Union's DML workshop can be found at http://ada00.math.uni-bielefeld.de/mediawiki-1.18.1/index.php/, updated April 23, 2013.

## STRUCTURE OF THE REPORT

This report consists of five main chapters and several appendices. The rest of this chapter discusses previous digital mathematics library efforts, the universe of mathematical information, relevant conceptual tools, and current mathematical resources. Chapter 2 discusses what is missing from the mathematical information landscape and what gaps the DML would fill, and elaborates on the desired DML capabilities from a user's perspective. This includes a discussion of what types of features would make the mathematical literature and current resource capability more meaningful to a mathematical researcher. Chapter 3 discusses some of the broad issues that the DML would face during development, including developing partnerships, managing large data sets, navigating open access, and planning for system and data maintenance. Chapter 4 provides a strategic plan for the development of the DML, including a discussion of fundamental principles, the constitution of a governing organization, steps toward initial development, and resources that would be needed. Chapter 5 discusses some details of entity collections and technical considerations for the DML that will be needed to make the features and capabilities discussed in Chapter 2 a reality.

In preparing this report, the committee reviewed many existing digital resources for mathematics, as well as relevant initiatives in some other sciences. A brief discussion of these tools is given in Appendix C.

## PREVIOUS DIGITAL MATHEMATICS LIBRARY EFFORTS

The idea of a comprehensive digital mathematics library has been around for decades, and there have been several incarnations of the idea with different foci. The first step in this vision was retrospective digitization of the older parts of the literature that did not already exist in digital form, and this has largely been achieved (though the quality, and hence utility, of these converted materials varies widely, ranging from simple page scans to carefully proofread markups).

The Cornell University Digital Mathematics Library Planning Project was funded by the National Science Foundation from 2003 to 2004 as a step "toward the establishment of a comprehensive, international, distributed collection of digital information and published knowledge in mathematics."[5] Its vision statement reads as follows:

> In light of mathematicians' reliance on their discipline's rich published heritage and the key role of mathematics in enabling other scientific disci-

---

plines, the Digital Mathematics Library strives to make the entirety of past mathematics scholarship available online, at reasonable cost, in the form of an authoritative and enduring digital collection, developed and curated by a network of institutions.

A follow-up report from the International Mathematical Union (IMU, 2006) shared this vision of a distributed collection of past mathematical scholarship that served the needs of all science, and it encouraged mathematicians and publishers of mathematics to join together in implementing this vision. However, it was clear within a few years that this vision was not going to become a reality soon. As David Ruddy of Project Euclid wrote (Ruddy, 2009):

> The grand vision of a Digital Mathematics Library, coordinated by a group of institutions that establish policies and practices regarding digitization, management, access, and preservation, has not come to pass. The project encountered two related problems: it was overly ambitious, and the approach to realizing it confused local and community responsibilities. While the vision called for a network of distributed, interoperable repositories, the committee approached and planned the project with the goal of building a single, unified library.

At the time of this study, there has been some progress in this vision of a single, unified library in the form of the European Digital Mathematics Library (EuDML) project.[6] The EuDML project, funded from 2010-2013 by the European Commission, created a network of 12 European repositories acquiring selected mathematical content for preservation and access and made progress in establishing a single distributed library with a collection of about 225,000 unique items, spanning 2.6 million pages. The EuDML succeeded in creating a unified metadata framework[7]—which includes items about a document such as the title, authors, abstract, comments, report number, category, journal reference, direct object identifier, Mathematics Subject Classification (MSC), and Association for Computing Machinery (ACM) computing classification—that is shared by these repositories and providing a single point of access to publications in these repositories, albeit with limited rights to search the full text from some sources. Impressive as the EuDML is, when compared to the full size and scope of the universe of published mathematics (described in the next section), and given the essential requirement to integrate with copyrighted materials and the clear desirability and cost-effectiveness of leveraging existing repositories and services, the EuDML experience only emphasizes the difficulties inherent in aiming for a single, centrally managed and truly comprehensive collection of digitized mathematics as the cornerstone for a comprehensive DML. With the advent of recent advances in technology and the advantage of experience gained on EuDML and other projects, the study committee concluded that a more effective approach going forward would be to partner with existing content providers and focus instead on the innovations and elements of shared infrastructure and knowledge management that are not being adequately addressed by other entities (i.e., rather than on central harvesting and aggregation of primary content). The committee believes that this vision is consistent with the original vision of the EuDML, although it was not realized by that project.

---

[6] T. Bouche, Université de Grenoble, "From EuDML to WDML: Next Steps," Presentation to the committee on November 27, 2012.

[7] European Digital Mathematics Library, "Appendix, EuDML Metadata Schema (Final)/ Tagging Best Practices," in EuDML Metadata Schema Specification (v2.0-final), https://project. eudml.org/sites/default/files/d36-appendix_uncropped.pdf, accessed January 16, 2014.

Another example of an online resource that helps users connect with knowledge is the National Science Digital Library (NSDL).[8] NSDL is an online educational resource for teaching and learning, with current emphasis on the sciences, technology, engineering, and mathematics. NSDL does not hold content directly—instead, it provides structured metadata about Web-based educational resources held on other sites by providers who contribute this metadata to NSDL for organized search and open access to educational resources via NSDL.org and its services.

A discussion of many other efforts and current digital resources can be found in Appendix C.

The Alfred P. Sloan Foundation supported a World Digital Mathematics Library workshop in June 2012,[9] which was planned by the IMU's Committee on Electronic Information and Communication. This workshop provided a wealth of information to the committee on the current state of the art and research efforts aimed at making the World Digital Mathematics Library a reality.

Much of the straightforward work of assembling digital mathematics libraries has been done (e.g., digitizing material, aggregating it into small to medium-sized collections). The difficulties that the EuDML faced in creating a single large aggregation of mathematics literature and the difficulty of other World Digital Mathematics Library efforts in gaining community support indicates that these challenges are unlikely to be overcome soon. The committee notes that there has been sizable ongoing investment from publishers (both commercial and noncommercial) to retrospectively digi-tize historical runs of their copyrighted journals and also, in many cases, even earlier historical materials that are now out of copyright, in order to capture comprehensive representations of their journals. However, broad services such as Google Scholar now provide much of the functionality that many of these specialized efforts had hoped to achieve in building comprehensive and coherent collections of the mathematical literature. Such services achieve this functionality by searching *across* a range of repositories, rather than trying to collect all of the material in one (or a very few) repositories. In the committee's view, efforts to build centralized comprehensive resources are reaching a point of diminishing returns.

---

[8] National Science Digital Library, http://nsdl.org/, accessed January 16, 2014.

[9] International Mathematics Union, "The Future World Heritage Digital Mathematics Library: Plans and Prospects," updated April 23, 2013, http://ada00.math.uni-bielefeld.de/mediawiki-1.18.1/index.php/Main_Page.

> **Finding: The construction of mathematical libraries through centralized aggregation of resources has reached a point of diminishing returns, particularly given that much of this construction has been coupled with retrospective digitization efforts.**

While there is still a substantial amount of historical (mostly out of copyright) mathematical literature that would benefit from retrospective digitization, or higher quality digitization than has currently been done, the committee does not believe that there is justification for a major new program and investment in this area. In particular, although there is value in modest, sustained investment in existing efforts, these will make only incremental contributions. While the fundamental importance of the heritage literature remains, its size, as a fraction of the overall mathematics literature, is diminishing steadily. No amount of additional retrospective digitization will result in a fundamental change in the way that the mathematical literature can be used in new ways or evolved to meet new research needs. Moreover, while the historical (e.g., out of copyright) segments of the mathematical literature are valuable, any genuinely meaningful large-scale change in accessing the mathematical literature and knowledge base *must* encompass not only heritage but also current literature. Thus, the committee believes that a very different set of investments (as described in this report) is where the transformative opportunities await.

The next section provides some more detailed information on the existing landscape of mathematical literature and how much has been digitized.

## THE UNIVERSE OF PUBLISHED MATHEMATICAL INFORMATION

Mathematics shares more with the arts than the sciences, in that its primary data are human creations, perhaps representations of ideas in a platonic realm, rather than data derived by observation or measurement of the physical universe. Mathematical information is primarily mined from its own literature or derived by computation. This section describes the state of mathematical publishing and the world of mathematical objects that exist within the publications.

## Digital Mathematical Publications

Most of the mathematics literature of the 20th century is now available digitally. Through the Jahrbuch Electronic Research Archive for Mathematics[10] project and the independent efforts of publishers and others, much of the most important mathematical research of the last half of the 19th century also has been digitized. Appendix C provides an overview of the many sources for digitized mathematical source material, including repositories and many other types of sources, whether freely accessible or behind paywalls (and thus only accessible to subscribers). A large part of the mathematics literature in electronic form consists of papers written in the past 20 years. This portion of the literature is searchable and navigable by any user of a library with access to the main subscription services controlled by libraries and publishers.

In addition, a considerable body of the heritage literature in mathematics has been digitized over the past 15 years. The most comprehensive listing of the retro-digitized mathematics literature is Ulf Rehmann's list of Retrodigitized Mathematics Journals and Monographs,[11] which is a list of titles of serials and books that have been digitized without metadata.[12] Much of this metadata has found its way into indexes maintained by Google, MathSciNet, and Zentralblatt (zbMATH).[13]

The digital corpus of mathematics literature is extensive. The MathSciNet[14] database includes approximately 2.9 million publications from 1940 to the present, with direct links to 1.7 million of them. MathSciNet currently indexes more than 2,000 journal/serial titles and contains about 100,000 books (post 1960). Of the items currently available on MathSciNet, 2.6 million of them are from the 1970s or later, and 1.7 million are from 1990 onward. The American Mathematical Society has kept track of new journal titles in the field since 1997, and there has been an average growth of about 40 new journal titles per year in mathematics.

---

[10] The Jahrbuch Project, Electronic Research Archive for Mathematics, last modified October 31, 2006, http://www.emis.de/projects/JFM/.

[11] DML: Digital Mathematics Library, http://www.mathematik.uni-bielefeld.de/~rehmann/DML/dml_links.html, accessed January 16, 2014.

[12] Metadata are broadly defined as data about data. In the case of a typical mathematics journal digital publication, metadata may include information such as author, journal name and volume, date of publication, time of file creation, size of file.

[13] zbMATH, http://zbmath.org/, accessed January 16, 2014.

[14] American Mathematical Society, MathSciNet, http://www.ams.org/mathscinet/, accessed January 16, 2014.

zbMATH (1931-present) contains more than 3 million publications and currently indexes approximately 3,500 journals. The annual production of mathematics papers is more difficult to quantify. There has been a steady increase in the number of math papers added to arXiv[15] over the past 5 years (shown in Table 1-1), although it is not clear from these data if this shows an increase in mathematics publications or an increase in mathematicians' willingness to post their papers. Annual entries on MathSciNet and the number of mathematics papers listed in Web of Science[16] have both remained relatively constant around 90,000 and 20,000, respectively (see Tables 1-2 and 1-3).

Components of the digitized corpus of mathematics are increasingly included in a variety of stable, well-curated repositories, although access to much of this corpus remains limited by copyright or other intellectual rights restrictions. For example, in terms of retrospectively digitized works cataloged under the subject heading (or subheading) of "mathematics," the *HathiTrust Digital Library*[17] includes approximately 40,000 bibliographically distinct resources.[18] Of these, only 6,800 were digitized from public-domain works; the rest were digitized from copyrighted originals. These numbers are a mix of monograph titles and serial titles (a serial title in *HathiTrust* typically encompasses a complete run of a journal, edited series, or conference publication series). Each serial run could be expected to include tens or even hundreds of issues, with each issue containing at least several articles or papers. In terms of pages, using the *HathiTrust* repository-wide ratio of pages per bibliographic resource to estimate, this translates to a rough estimate of 25.5 million pages of retrospectively digitized mathematics in *HathiTrust* with approximately 17 percent (6,800 out of 40,000) digitized from public-domain sources.

The basic trends seem clear: more and more of the corpus of mathematical literature will be in digital form, including some with high-quality markup, specifically those items that are "born" digital or retro-digitized to be in a machine readable format and that use typesetting such as LaTeX or MathML (as opposed to page images of publications). As mentioned before, the fraction of the overall corpus that is pre-1970 is rapidly diminishing due to the relative explosion in the annual rates of publication in recent decades (however, this should in no way be seen as diminishing the fundamental importance of heritage literature).

---

[15] arXiv, http://arxiv.org/, accessed January 16, 2014.
[16] Thomson Reuters, "Web of Science Core Collection," http://thomsonreuters.com/web-of-science/, accessed January 16, 2014.
[17] *HathiTrust Digital Library*, http://www.hathitrust.org/, accessed January 16, 2014.
[18] Current as of September 2013.

TABLE 1-1 Number of Mathematics Papers Added to arXiv Annually Between 2008 and 2012

| Year | Mathematics Papers Added to arXiv |
| --- | --- |
| 2008 | 14,373 |
| 2009 | 16,319 |
| 2010 | 18,765 |
| 2011 | 21,287 |
| 2012 | 24,176 |

SOURCE: arXiv, http://arxiv.org/, accessed January 16, 2014.

TABLE 1-2 Number of Articles in Research Journals in MathSciNet Annually Between 2006 and 2012

| Publication Year | Entries in MathSciNet |
| --- | --- |
| 2006 | 76,187 |
| 2007 | 81,638 |
| 2008 | 86,533 |
| 2009 | 87,279 |
| 2010 | 87,162 |
| 2011 | 89,638 |
| 2012 | 92,191 |

NOTE: A steady growth of about 3 percent per year is seen.
SOURCE: American Mathematical Society, MathSciNet, http://www.ams.org/mathscinet/, accessed January 16, 2014.

TABLE 1-3 Mathematics Papers Listed in Web of Science Annually Between 2008 and 2012

| Year | Mathematics Papers Listed in Web of Science |
| --- | --- |
| 2008 | 20,908 |
| 2009 | 22,390 |
| 2010 | 22,079 |
| 2011 | 22,716 |
| 2012 | 23,760 |

SOURCE: Thomson Reuters, "Web of Science Core Collection," http://thomsonreuters.com/web-of-science/, accessed January 16, 2014.

## Objects in the Mathematical Literature

Information found in the mathematical literature is diverse but largely falls into two main categories:

1. Bibliographic information, such as
   a. Documents (e.g., articles, books, proceedings, talks, diagrams, homepages, blogs, videos);
   b. People (e.g., authors, editors, referees, reviewers);
   c. Events (e.g., discoveries, publications, conferences, talks, births, deaths, degrees, awards);
   d. Organizations (e.g., universities, publishers, journals, libraries, service providers);
   e. Subjects (e.g., major branches of mathematics—algebra, geometry, analysis, topology, probability, statistics—as well as their intersections and interactions and their various sub-branches, down to even finer topics and including ubiquitous mathematical terms like "number," "set")
2. Mathematical concepts (e.g., axioms, definitions, theorems, proofs, formulas, equations, numbers, sets, functions) and objects (e.g., groups, rings).

Collecting and aggregating mathematical bibliographic information has been the path many digital libraries and digital resources have taken in the past (Chapter 2 and Appendix C discuss many of these efforts to date). While there are many challenges in collecting this information, the even more difficult work lies in collecting mathematical concepts, which lack the standardization that most bibliographic information has acquired. However, an ability to explore these mathematical objects within the literature offers the potential to uncover currently under-explored connections in mathematics.

The recent National Research Council report *The Mathematical Sciences in 2025* (NRC, 2013) discusses the importance of mathematical structures, which are part of the larger mathematical concepts described above:

> A mathematical structure is a mental construct that satisfies a collection of explicit formal rules on which mathematical reasoning can be carried out. . . . What is remarkable is how many interesting mathematical structures there are, how diverse are their characteristics, and how many of them turn out to be important in understanding the real world, often in unanticipated ways. Indeed, one of the reasons for the limitless possibilities of the mathematical sciences is the vast realm of possibilities for mathematical structures. . . . A striking feature of mathematical structures is their hierarchical nature—it is possible to use existing mathematical structures as a foundation on which to build new mathematical structures
. . . . Mathematical structures provide a unifying thread weaving through
and uniting the mathematical sciences. (pp. 29-30)

Given the size, diversity, and inherent nature of mathematics information in categories 1 and 2 above, it is clearly not sufficient to simply provide undifferentiated access to the universe of mathematics monographs, journal articles, and conference papers. Instead, the online research literature of mathematics must be organized into a well-structured network of resources linked together based on a variety of attributes—bibliographic and topical, of course, but also linked in a highly granular fashion on commonalities of mathematical structures and the shared use of mathematical objects, reasoning, and methodologies. The committee believes that the greatest potential for the DML lies in providing mathematicians access to a well-structured network of information and building services that both enhance and utilize this data. In the context of today's Web environment, a well-structured network implies adherence to the Semantic Web[19] and linked open data principles and to community-endorsed standards and best practices. While the foundation for such a well-structured network of digital research mathematics exists in established repositories and component digital libraries, the underlying thesauri and ontologies of mathematical objects do not yet exist (or have not yet been given permanence and formal identity), and the agreements on best practices for interoperability and the implementation of linked open data principles in the context of research mathematics repositories have not yet been reached.

## CONCEPTUAL TOOLS

General conceptual tools that are used to structure, organize, represent, and share knowledge include the closely related ideas of ontologies, taxonomies, and vocabularies. There is considerable debate about the precise definitions and differences among these tools, although ontologies (most commonly viewed as a tool for defining some classes of objects—the attributes that these objects may have and the way in which these objects may be related to each other) are usually seen as the most general formulation (Gruber, 2009). Taxonomies are specific, usually hierarchical, collections of terms that can be used to describe or classify objects in some contexts—examples of these include subject headings or the naming schemes used in biological systematics. "Controlled" vocabularies are collections of values that can be used to populate specific instances of object attributes within an ontology; in a certain sense, they are equivalent to taxonomies in that they can be used to classify. However, controlled vocabularies are often "flat," without other internal structure among the possible values, whereas taxonomies commonly include very rich internal hierarchical structure. Ontologies, vocabularies, and taxonomies work together. As a simple example, a part of an ontology might define a specific class of objects called documents; each of these has attributes that include subjects and languages. One might have a list of possible language values (a controlled vocabulary) associated with the ontology and also a tree structure of subject headings (a taxonomy, though it could also viewed as a simple vocabulary).

---

[19] W3C, "Semantic Web," http://www.w3.org/standards/semanticweb/.

For instance, within the mathematical sciences, the widely accepted Bibliographic Ontology[20] provides a fairly adequate accounting of the many common relations between objects in categories 1a through 1e listed above. The BibTeX[21] schema that describes the structure of BibTeX records defines a similar ontology. The Citation Typing Ontology (CiTO)[22] is an ontology for description of the citation relation between documents. The Mathematics Subject Classification (MSC2010)[23] provides a very well thought out, largely hierarchical taxonomy for the classification of mathematical documents by subject, and thence for the subjects themselves. OpenMath,[24] discussed further in Chapter 5, offers a potential standard for representing the semantics of mathematical objects that is very relevant to the DML's goals.

The application of such ontologies to a mathematical objects data set can create graphical structures of information that can provide new insights. For instance, citations generate a citation graph, and collaborations generate a collaboration graph. Such graphical structures are commonly embedded in the structure of hyperlinked webpages, thereby connecting literature that was not obviously related otherwise.

Development of new ontologies is a complex process requiring a high level of community effort for consensus, even for limited sets of relations. The committee expects that when communities start to curate various digital collections of records of mathematical entities, there will be some "bottom up" development of at least minimal ontologies for these entities, as has already occurred with MSC2010 and OpenMath. The structure of these ontologies will be reflected in the necessary schemas[25] for description of the objects they involve, and the graphical relations induced by these ontologies will be of potentially great interest in the process of extracting information and knowledge from mathematical publications.

---

[20] The Bibliographic Ontology, "Bibliographic Ontology Specification," dated November 4, 2009, http://bibliontology.com/specification.

[21] BibTeX, http://www.bibtex.org/, accessed January 16, 2014.

[22] CiTO, the Citation Typing Ontology, dated March 7, 2013, http://purl.org/spar/cito/.

[23] Encoded by the Mathematics Subject Classification (MSC2010), American Mathematical Society, http://www.ams.org/mathscinet/msc/msc2010.html, accessed January 16, 2014.

[24] OpenMath Society, OpenMath, http://www.openmath.org/, accessed January 16, 2014.

[25] A schema is broadly defined as a representation of a plan or theory in the form of an outline or model.

## CURRENT MATHEMATICAL RESOURCES

The management of formal representations of mathematical concepts is known as mathematics knowledge management (Carette and Farmer, 2009). In this report, this issue is viewed more broadly as the management of mathematical information and concepts, both formal and informal, including the bibliographic information and mathematical concepts categories of objects introduced in the previous section, only the latter of which can be usefully regarded as part of mathematics itself.

### Bibliographic Resources in Mathematics

Several general bibliographic resources exist, and some of these are described in Appendix C. Among them, mathematicians typically use Google[26] and Google Scholar[27] most often, although CrossRef[28] is "under the hood" whenever a user navigates from one publisher's site to another by a reference link. While many mathematicians heavily utilize these general information services because of their power and ubiquity, some mathematicians prefer the discipline-specific abstracting and indexing services provided by MathSciNet[29] and zbMath.[30] This discipline-specific service preference is partly for historical reasons and partly because the focus and quality of metadata provided by these services in mathematics makes it easier to find publications of interest. Both services offer bibliographic entries in BibTeX,[31] which is machine-readable and reusable, for preparation of reference lists for LaTeX[32] documents, and, with more technical effort, for publication of online bibliographies in HTML[33] or JSON.[34] Using search engines with access to well-curated bibliographic metadata and full-text indexing is how most mathematicians find mathematical primary sources today.

---

[26] Google, https://www.google.com/, accessed January 16, 2014.

[27] Google Scholar, http://scholar.google.com/, accessed January 16, 2014.

[28] CrossRef, http://www.crossref.org/, accessed January 16, 2014.

[29] American Mathematical Society, MathSciNet, http://www.ams.org/mathscinet/, accessed January 16, 2014.

[30] zbMATH, http://www.zentralblatt-math.org/zmath/, accessed January 16, 2014.

[31] BibTeX, http://www.bibtex.org/, accessed January 16, 2014.

[32] LaTeX—A document preparation system, last revised January 10, 2010, http://www.latex-project.org/.

[33] "HTML," *Wikipedia*, http://en.wikipedia.org/wiki/HTML, accessed January 16, 2014.

[34] "Introducing JSON," http://www.json.org/, accessed January 16, 2014.

Services such as MathSciNet, zbMATH, and Google Scholar provide complementary and somewhat overlapping services. One distinct difference is that MathSciNet is organized chronologically and referentially, while Google Scholar is based on "importance" as qualified by page ranks or some variant thereof. Both are important and are used in literature searches. MathSciNet is great for tasks such as listing all articles by an author and listing all articles in a specific mathematical field, and it has high-quality metadata that are needed for many purposes. Its search capabilities are limited because it only searches over metadata. Google Scholar is often better for searches because it searches over full text, including reference lists, and has better ranking or returns for most purposes. One issue that some mathematicians have with Google Scholar is that it is not possible to limit searches to math or subfields of math. MathSciNet, zbMATH, and Google Scholar combined do a good job providing conventional discovery over the corpus of traditionally published mathematical literature, but no services currently provide a finer-grain search capability that allows a user to search for mathematical objects or ideas that cannot be easily defined by text search, such as an equation or the evolution of a specific notation. Ideally, a mathematician should have the best of both capabilities through a single interface, but this is challenging because neither MathSciNet nor Google Scholar currently allow their data to be merged with the other's.

Mathematicians also make extensive use of arXiv as a platform for sharing preprints and keeping up with current research developments. Mathematicians strongly support arXiv in part because the full text is largely indexed and exposed to the Web through search engines. However, arXiv items are not indexed through services such as MathSciNet or zbMATH, which would help connect these items to the rest of the literature. Search tools associated with distinct subsets of the literature, such as arXiv, publisher-based repositories, library catalogs, and academic institutional repositories provide overlapping access to the mathematical literature. Unfortunately, the present configuration of these discipline-specific tools does not provide a single information source where mathematicians can find and access information from diverse sources, and the more general information sources often lack the mathematical metadata and details that make mathematics literature easy to search and browse.

Combining data from multiple information resources (e.g., Google, MathSciNet, zbMATH) is complicated. Partnering organizations would have to allow their data to be collected, reused, or recombined on a large scale, which many services are hesitant to do. Even seemingly open resources (such as arXiv) may have legal restrictions on outside data aggregation, depending on what is done with the data. This collaboration would have to be negotiated between potential partners with the goal of creating a unified view of the mathematics literature. Some approaches toward developing partnerships and relevant examples are discussed in Chapter 3.

Given the central importance of bibliographic data searches and the repeated use of bibliographic information by researchers in preparation of research articles, it is essential for the DML to provide adequate bibliographic support tools with access to the best available bibliographic data in mathematics and related fields. Ideally, it should support advanced bibliographic data processing to detect and identify the structure of networks of papers, authors, topics, and the like. The foundations of such bibliographic data processing are provided by the larger existing bibliographic services in mathematics and beyond, especially MathSciNet, zbMATH, and Google Scholar, which are the most commonly used by mathematicians. At present, none of these services provides an application programming interface (API) for programmatic access, and none of them allow their data to be downloaded in bulk, except with severe restrictions on what can be done with it. To provide the greatest benefit to users of a DML, that would have to change. Both EuDML and Microsoft Academic Search provide steps in a positive direction with more or less open bibliographic data stores with an API for access, which allows tools and services to be built over the corpus.

To seriously engage the mathematics world with a digital library system, extensive coverage of mathematical information is essential. The committee considered whether the DML could initially focus on out-of-copyright material, but it concluded that there would not be community support or interest in this approach because it is too limited. On the other hand, much progress has been made in digitizing heritage content, and it is essential that this be integrated with the rest of the math literature base.

### Specialized Mathematical Information Resources

General bibliographic services provide limited support for navigating and searching mathematical literature below the top five bibliographic classes (documents, people, events, organizations, subjects) discussed above. Beyond these five universal classes, information storage and retrieval for math-specific entities is fragmented and typically does not have links or references to the main indexing services.[35]

Research mathematics literature includes a diverse range of special objects—e.g., theorems, lemmas, functions, sequences—that are not represented adequately, or sometimes at all, in full-text indexing and article-level subject classification systems. Currently, these objects are computationally expensive and difficult to recognize through machine-based methods alone. Ontologies of objects—such as reference volumes that enumerate classes of functions, sequences, and other objects—have been developed and curated by mathematicians for centuries. These resources include mathematical handbooks, some of the most famous being the following:

- Abramowitz and Stegun (1972) and the subsequent Digital Library of Mathematical Functions,[36]
- The Bateman Manuscript,[37]
- Gradshteyn and Ryzhik (2007),
- Borodin and Salminen (2002), and
- The Princeton Companion to Mathematics (Gowers et al., 2008).

---

[35] MathSciNet and zbMATH share the MSC2010 subject classification, which provides some basic filtering of bibliographic data by subject. ArXiv uses a coarser classification, which is however easily mapped to sets of top-level MSC 2010 categories.

There are also examples of more recently developed resources that provide collections of some mathematical objects, including the following:

- *Propositions: Wikipedia*'s List of Theorems,[38] Mizar[39];
- *Proofs:* Proofs from the Book (Aigner and Ziegler, 2010), Mizar, Coq,[40] and others[41];
- *Numbers:* A Dictionary of Real Numbers (Borwein and Borwein, 1990);
- *Sequences:* The On-Line Encyclopedia of Integer Sequences (OEIS)[42];
- *Functions:* Digital Library of Mathematical Functions,[43] Wolfram *MathWorld*,[44] Wolfram Functions Site[45];
- *Groups, rings, and fields: Wikipedia*'s List of Simple Lie Groups,[46] *Wikipedia*'s List of Finite Simple Groups,[47] Centre for Inter-

---

disciplinary Research in Computational Algebra: Finite Fields,[48] Sage's Finite Fields[49];

- *Identities:* Piezas[50]; Petkovsek et al. (1996);
- *Inequalities: Wikipedia*'s List of Inequalities,[51] DasGupta (2008); and
- *Formulas:* Springer LaTeX Search,[52] Hijikata et al. (2009), Kohlhase et al. (2012).

From a review of these lists, as well as the resources discussed in Appendix C, it is clear that authors and editors continue to be motivated to create and publish lists of various kinds of mathematical objects. Some of these lists, especially ones like tables of integrals and lists of sequences, provide very useful tools for mathematicians and other users of mathematics, especially when combined with computational resources. *Wikipedia* currently plays a key role in supporting distributed creation and maintenance of numerous lists of serious interest to mathematicians.

Lists and tables have been an essential part of mathematical research throughout history, and the vast majority of working mathematicians have made use of appropriate tables (or, more recently, the equivalent numerical or symbolic software) in the course of their research. The most basic are numerical tables (e.g., values of logarithms, trigonometric functions, various special functions, zeros of the zeta function, integer sequences). More sophisticated are lists of mathematical objects (e.g., indefinite and definite integrals, finite simple groups, Fourier transforms, partial differential equations and their solutions). Or, at even a higher level, lists of theorems, concepts, etc.

At their most basic, tables provide a simple mechanism for speeding up research. Once one identifies that an object under investigation appears in a table, one can make use of prior knowledge about said object, thereby facilitating either applications or new advances in theory. Compiling a table is an important research contribution in its own right, helping codify the knowledge in a field, point out gaps therein, and inspire new research to fill in and extend what is known. Scanning a table often enables one to spot

---

[48] CIRCA, "GAP Instructional Material," January 2003, http://www-circa.mcs.st-and.ac.uk/gapfinite.php.

[49] Sage Development Team, "Finite Fields," http://www.sagemath.org/doc/reference/rings_standard/sage/rings/finite_rings/constructor.html, accessed January 16, 2014.

[50] T. Piezas III, A Collection of Algebraic Identities, https://sites.google.com/site/tpiezas/Home/, accessed January 16, 2014.

[51] "List of Inequalities," *Wikipedia*, last modified November 28, 2013, http://en.wikipedia.org/wiki/List_of_inequalities.

[52] Springer, LaTeX Search, http://www.latexsearch.com/, accessed January 16, 2014.

otherwise obscure patterns, leading to new theorems and new directions of research.

Sara Billey and Bridget Tenner wrote that a database for mathematical theorems would "enhance experimental mathematics, help researchers make unexpected connections between areas of mathematics, and even improve the refereeing process" (Billey and Tenner, 2013, p. 1093). Extensive lists could also enhance search and retrieval of mathematical information and allow for connections to be made between mathematical topics and objects.

Currently, there are no satisfactory indexes of many mathematical objects, including symbols and their uses, formulas, equations, theorems, and proofs, and systematically labeling them is challenging and, as of yet, unsolved. In many fields where there are more specialized objects (such as groups, rings, fields), there are community efforts to index these, but they are typically not machine-readable, reusable, or easily integrated with other tools and are often lacking editorial efforts. So, the issue is how to identify existing lists that are useful and valuable and provide some central guidance for further development and maintenance of such lists.

Chapter 2 of this report discusses some of the user features that could advance mathematics research by increasing connections, and Chapter 5 discusses what collections of entity lists could start making these features and this connectivity a reality.

# 2

# Potential Value of a Digital Mathematics Library

## WHAT IS MISSING FROM THE MATHEMATICAL INFORMATION LANDSCAPE?

The current mathematical information landscape is complex and diverse, as described in Chapter 1 and Appendix C. Current digital mathematical resources provide services such as electronic access to papers (often with advanced features capable of searching and sorting based on key words, subject areas, text searches, and authors), platforms for discussion, and improved navigation across multiple data sources. What they do not do is allow a user to systematically explore the information captured within the literature and forums and readily explore connections that may not be obvious from looking at the material alone.

This inability to easily explore the mathematical ideas that exist within a mathematical paper, which cannot easily be searched for, is a detriment to the mathematical community. There is a largely unexplored network of information embedded in the connections of mathematical objects, and formalizing this network—making it easy to see, manipulate, and explore—holds the potential to vastly accelerate and expand currently mathematical research. This network would consist of information from traditional resources, such as research papers published in journals, and content dispersed in other Internet-based resources and databases. Initial development of the DML could begin immediately with the aim of providing a foundational platform on which most of the capabilities discussed in this report might imaginably be achieved in a 10- or 20-year time frame. This report discusses how the Digital Mathematics Library (DML) can make this network of information a reality.

## WHAT GAPS WOULD THE
## DIGITAL MATHEMATICS LIBRARY FILL?

The real opportunity is in offering mathematicians new and more direct ways, through the Web, to discover and explore relationships between mathematical concepts (such as axioms, definitions, theorems, proofs, formulas, equations, numbers, sets, functions) and objects (such as groups, rings) and broader knowledge (such as the evolution of a field of study; and relationships between mathematical fields, concepts, and objects). Improved discovery and interaction in the proposed DML would make it possible to find and examine material on a much finer scale than what is currently possible, making connections easier to find, shortening the needed start-up time for new research areas, and formalizing some of the logic that mathematicians are already using in their research.

In *Probability Theory: The Logic of Science*, E.T. Jaynes discusses the reasoning that many mathematicians go through when approaching their work. He describes the strong form of reasoning as variations on the following: "If A is true, then B is true. A is true; therefore, B is true." Weaker forms are assertions, such as "If A is true, then B is true. B is true; therefore, A becomes more plausible." Jaynes states that

> [George] Pólya showed that even a pure mathematician actually uses these weaker forms of reasoning most of the time. Of course, when he publishes a new theorem, he will try very hard to invent an argument which uses only the first kind; but the reasoning process which led him to the theorem in the first place almost always involves one of the weaker forms (based, for example, on following up conjectures suggested by analogies). The same idea is expressed in a remark of S. Banach (quoted by S. Ulam, 1957): "Good mathematicians see analogies between theorems; great mathematicians see analogies between analogies." (Jaynes, 2003, p. 3)

The DML could help make these analogies easier to find and use.

Box 2.1 provides an example of how a mathematics researcher would start looking into a new topic, using Gröbner bases as a specific illustration. It shows some of the initial resources that are typically used and how their information varies from, complements, and supplements the other resources. It also shows how useful it would be to be able to pull much of this information into a unified source and make additional connections to other, lesser known resources and aspects of the literature.

The DML could aggregate and make available collections of ontologies, links, and other information created and maintained by human contributors and by curators and specialized machine agents with significant editorial input from the mathematical community. The DML could afford functionalities and services over the aggregated mathematical literature.

**BOX 2.1**
**How a Mathematics Researcher May Currently**
**Approach Information Gathering**

Gröbner bases were first introduced by Bruno Buchberger for solving a range of problems in computational algebra and became an essential component of computer algebra software (Buchberger, 2006). Suppose a mathematician wanted to find out about this topic, perhaps because it was needed for a particular problem. First, when one types "Grobner basis" into MathSciNet, a list of around 2,400 chronologically ordered items appears, most of which are specialized papers. This is a potentially good resource for a specialist but is probably not ideal for the novice. If a similar search is done via Google Scholar, a list of research articles and books on the subject appear and are ordered by "popularity," which usually reflects some version of page ranking. While some of the references provided by Google Scholar can be viewed, including some books on Google books, others are behind paywalls or are books that must be purchased before reading. In Google itself, the top five links are to *Wikipedia*,[1] *MathWorld*,[2] *Scholarpedia*,[3] Mathematica code,[4] and a survey article by Bernd Sturmfels.[5]

The *Wikipedia* article is limited and only contains four references but includes the book of Cox, Little, and O'Shea (1997), which is widely recognized and a premier introductory text on the topic. *Wikipedia* also offers suggested further reading and external links. Sturmfels's article, from the "What is . . ." section of the Notices, is terse and contains only three references, but one of them is the aforementioned book. *MathWorld*'s article is short and lacks any specifics, but it contains a significantly longer list of references, survey articles, and several links to Amazon for buying books (and at least one dead link). The *Scholarpedia* article, written by Bruno Buchberger and Manuel Kauers, is more comprehensive and includes many illustrations, a wide range of applications, and a long list of references, including a Gröbner bases bibliography compiled by Buchberger and his coworkers at the Research Institute for Symbolic Computation.[6] Unfortunately, no links are supplied in the *Scholarpedia* article to the other references. In many ways, *Scholarpedia*, which bills itself as a "peer-reviewed open-access encyclopedia," could serve as one possible model for some aspects of the proposed DML.

All of these resources combined, along with the tenacity to pursue the variety of resources, can result in a good start in understanding Gröbner bases. However, suppose the researcher was working in an area that led to questions that Gröbner bases could be profitably used in, but, not being an algebraist, he/she did not know that they existed or even how to start to query any of the standard tools. Vice versa, suppose the researcher works in Gröbner basis theory and find results that could lead to advances in an area that he/she is not familiar with; how would the researcher know?

Here's a real example: Although not well known, in fact, the theory of Gröbner bases was essentially discovered in 1910-1913 by an obscure Georgian mathematician, N.M. Gjunter, in his study of the integrability of overdetermined systems of partial differential equations (Renschuch et al., 1987). It is not immediately obvious through reference searching or the standard literature that Gröbner bases are of importance in partial differential equations (although the *Scholarpedia* article does mention some applications to ordinary differential equations). Moreover, the latter area has resulted in the refined and potentially very useful concept of an involutive basis. This particular gap could be filled by editing the above mentioned articles, but this is simply one of innumerable similar cases. Making such unexpected links is not currently easy but could become so with a fully functioning DML, therefore increasing the serendipitous-like discovery of connections, which plays a role across research.

**BOX 2.1 Continued**

While it would have to store modest amounts of new knowledge structures and indices, it would not have to generally replicate mathematical literature stored elsewhere.

The committee identified a number of basic desired library capabilities, including aggregation and documentation of information, annotation, search and discovery, navigation, and visualization and analytics. Properly implemented across the domain of mathematics research literature, these capabilities and resulting enhanced functionalities would not only facilitate better and more efficient search and discovery, but also allow mathematicians to interact with the research literature in new ways and at new levels of granularity. The proposed DML is much more than an indexing service and aims to create meaningful connections between topics by utilizing lists of entities and providing coherent access to a range of tools that can speed up mathematical discovery: for example, comprehensive encyclopedia articles and review articles, lists of mathematical objects, implication diagrams, and annotated bibliographies, informal annotations, and comments on articles. These tools and others are discussed in Chapter 5.

The DML would not only result in new efficiencies, thereby freeing up researcher time, but also enable experimentation with new approaches to using and getting the maximum benefit out of the mathematics research literature. The remainder of this section describes each of these desired capabilities and illustrates how resultant improved functionality could advance mathematics research.

### Aggregation and Documentation

Mathematicians want to be able to make searchable and sharable collections or lists of various kinds of mathematical objects easily, including bibliographies of the mathematical literature, perhaps with annotations. This is an area where it should be very easy to make rapid progress. The issues of mathematical object representations are mostly about who is allowed to create, view, and update various lists and about resource management. Many of these types of lists (such as those mentioned in Chapter 1) currently exist, some with connections to the literature, but their existence is often tied to the survival of the curator's personal website. Providing a stable platform for housing and connecting these lists would also allow for this information to be incorporated in the collective knowledge of the DML.

The availability and interconnection between these lists would allow a larger network of mathematical information to be developed. This would be on a finer scale than what is currently available and facilitate higher-level features of advanced search and navigation. The world of mathematical knowledge goes much deeper than the level of research papers; it goes down into the content that is discussed within the papers, the knowledge that is assumed already to be understood by the reader, and the connections that exist between this information. If the DML could draw on this information, it would have a much more meaningful view of mathematics.

### Lists of Mathematical Objects in New Areas

While many books contain fairly comprehensive descriptions of theorems relating to a specific subject and substantial stand-alone lists of theorems have been prepared, the committee is not aware of any truly comprehensive list of theorems in any branch of mathematics. Moreover, "lists" as embodied in books are not necessarily designed to enable all the functionality envisioned for the DML. There have also been several efforts to establish a formal computer-aided proof capability (Wiedijk, 2007), but it has not had much impact on the larger mathematics community. Mizar has published the largest such collection of about 50,000 formally checked theorems.[1] New mathematical theorems and lemmas are proven and pub-lished on a routine basis.[2] There have also been efforts to identify and list in order the most important mathematical theorems (such as the list presented by Paul and Jack Abad in 1999[3]) based on assessments of their place in the literature, each theorem's proof, and the unexpectedness of the result. Even if all existing theorems and lemmas were indexed and organized in some way, there needs to be a way to continually update this list with new work.

---

[1] Mizar Home Page, last modified January 8, 2014, http://mizar.org/.

Although even a list of theorems would be valuable, or a collection of text articles about each theorem, modern knowledge representation techniques offer more ambitious possibilities. For example, collections of represented facts such as DBPedia[4] or Freebase[5] permit retrieval of data about the real world, such as populations or areas of nations and towns. Library and museum catalogs are being converted to formal Resource Description Framework (RDF)[6] statements. Having a rigorous description of a theorem enables logical deduction and comparison of that theorem with others. The generality of mathematics is one of its beauties, and when the same form appears under two different names, it implies an unsuspected applicability of each theorem.

It appears within the grasp of modern information management tools to develop a machine-readable repository of mathematical theorems and definitions in which theorems are expressed as statements about terms, terms are linked to definitions, and definitions are constructed from logical statements about other terms. This is certainly very challenging, but the first steps in this direction have been made by Wolfram|Alpha for continued fractions, with a formalism for canonical representations of theorems that appears simple and flexible enough to be more widely adopted and used for purposes of search, retrieval, and linking. The Mizar Project also has a large database of formal theorem statements and formal proofs, although this is much less easily accessible to a working mathematician. How to do this on a large scale is still an open problem, but there are indications that efforts of this kind should be rewarding (Billey and Tenner, 2013).

Only the definitions and the theorem statements need to be machine-readable—the proofs can be LaTeX or a citation. Technologies like RDF

---

[2] In his 1998 biography of Paul Erdös, Paul Hoffman reports that mathematician Ronald Graham estimated that upwards of 250,000 theorems were being published each year at that time (Hoffman, 1987).

[3] P. Abad and J. Abad, The Hundred Greatest Theorems, 1999, http://pirate.shu.edu/~kahlnath/Top100.html.

[4] DBpedia, About, last modification September 17, 2013, http://dbpedia.org/About.

[5] Freebase, http://www.freebase.com/, accessed January 16, 2014.

[6] RDF is a standard model for data interchange on the Web and facilitates data merging even in the case of differing underlying schemas. See WC3 Semantic Web, "Resource Description Framework (RDF)," last modified March 22, 2013, http://www.w3.org/RDF/.

and OWL[7] may be useful for encoding the theorems' statements and the definitions. These technologies are flexible enough to allow users to extend the ontology, while encouraging reuse of existing terms. The markup languages used by automatic theorem provers could also be useful because they are sufficiently flexible to encode many important theorems, but they might not do enough to encourage reuse of terms.

The theorem and lemma repository would benefit from being accessible to programs via an application programming interface, which is a protocol used to allow software components to easily communicate with each other and may include specifications for routines, data structures, object classes, and/or variables.

Researchers will likely submit their theorems through a Web-based interface if it helps them to get citations and to stake a claim to having proved it first. There are a lot of famous cases where theorems were proven independently by multiple individuals using different terminology. A machine-readable repository could detect duplicate terms and theorems so that researchers can focus on new results rather than proving what is already known. The main benefit, however, may come from granting programs access to the latest mathematical results through user submissions.

Another data type worthy of consideration in a DML is problems. Good problems spur research advances. Problem lists have been created and maintained at various times, most famously Hilbert's list of problems around the beginning of the 20th century. Some recent efforts at curation of problem lists are the Open Problem Garden[8] and the the American Institute of Mathematics' Problem Lists.[9] A community feature encouraging creation and maintenance of problem lists with adequate links to the literature and indications of status could be an important component of the DML.

### Annotation

Mathematicians want to be able to annotate mathematical documents in various ways and share these annotations with collaborators or students and, in some cases, publish these annotations for the benefit of a wider but closed group (a set of collaborators, or a seminar, or a cohort of doctoral students) or the general public. The ability to easily share notes could improve the learning curve for researchers in new areas, provide opportunities to learn from other researchers interested in similar topics, elucidate logic that is not explicitly stated in papers, allow authors to post corrections, and overall enrich the research discussion. Some mathematicians prefer to keep comments limited to a smaller group, while others are more comfortable posting openly. Either way, this enhancement to the traditional research paper could quicken the path toward understanding and at the same time enhance the DML's capability to traverse the literature.

---

[7] W3C, "OWL Web Ontology Language Overview," February 10, 2004, http://www.w3.org/TR/owl-features/.

[8] Open Problem Garden, http://www.openproblemgarden.org/, accessed January 16, 2014.

[9] American Institute of Mathematics, AIM Problem Lists, http://aimpl.org/, accessed January 16, 2014.

The ability to see others' annotations as well as create new annotations would make reading a paper not only easier, but potentially more interesting. Some links could point to other items residing in the digital library, while others point to popular sites such as MathOverflow and *Wikipedia* or other sites outside the DML. For researchers setting out in a new direction or for researchers in an isolated location, it is often difficult to get involved in a lively conversation with fellow researchers. Links to discussions and comments on research papers and theorems could be a way to expand research discussions to a new level. Senior mathematicians could provide some general background information to research papers, such as a basic prerequisite for understanding the paper and some suggested readings; this would assist students and people starting out in a new direction. It should be possible for individual users to tailor the writing and reading of comments. It could also be useful to be able to select or prioritize, in several possible ways to be set by each user individually, the comments that appear on one's screen while searching (e.g., so as to see most prominently the comments from other members in an existing collaboration group or from a commenter one has experienced earlier as particularly insightful on a particular topic).

An important component of successfully providing an annotation feature within the DML is separating unhelpful comments and deciding which annotations will be kept in the system. Nearly every system that allows public comments also has a way to flag unconstructive comments and responses as inappropriate for that platform. A system such as this may need to be developed for the DML and refined based on the kinds of comments and feedback that the DML receives. One example of this is how MathOverflow deals with user input that its established users deem to be spam, offensive, or in need of attention for any other reason.[10] Elected community monitors are established within MathOverflow, and experienced users are able to flag comments and posts for a moderator's attention. The moderator can then decide what action is needed (deleting spam, closing off-topic posts, removing poorly rated posts, and so on). A system like this may work well for the DML.

---

[10] MathOverflow, Help Center, Reputation and Moderation, "Who are the site moderators, and what is their role here?," http://mathoverflow.net/help/site-moderators, accessed January 16, 2014.

General support for the creation of basic text annotations has been available for some time, including for mathematics literature made available as PDF or in HTML format. Support for more sophisticated forms of semantic annotation and for the sharing of annotations across disparate content repositories is rapidly maturing through technology from other domains,[11] but these technologies have yet to be customized for use in mathematics. Adapting these technologies to the mathematical community requires adequate support for mathematical markup. Some Web services are expanding into mathematical markup. For example, Authorea[12] uses a robust source control system in the backend (git) and an engine to understand LaTeX, Markdown, and most Web formats. Authorea lets users write articles collaboratively online, and it renders them in HTML5 inside a Web browser. Authorea is a spin-off initiative of Harvard University and the Harvard-Smithsonian Center for Astrophysics.

There are numerous other tools available that provide for "wiki-like" structured discussions with attribution dates and hierarchical organization, such as PBWorks.[13] There are also tools for highlighting, summarizing, providing video and audio annotation, mapping documents, and collaborative reading; some are specialized to particular document formats, and some are not. The Mellon project on Digital Research Tools[14] has a list of more than 500 tools, of which nearly 80 are tagged as annotation systems. Some are automated (e.g., part of speech tagging), but most are tools for use by readers or writers, either individually or in groups.

Adding this capability to the readily available digital literature should not be overly complicated. There would need to be conventions established for where the annotations are stored and who is responsible for storing them, and the best default setting for privacy and sharing would also need to be established. These annotations can provide a bridge to community-sourced markup of objects or a way to pass information to editors (human- or software-based) that curate the collection, thereby further enriching the DML. This is just one way in which user and community input would play a role in the DML; many others are listed elsewhere in this report. Community support for the new digital library will be essential for its success and

---

also an essential way in which it could be much more than just a collection of mathematical information and links to other repositories and services.

> **Recommendation: A primary role of the Digital Mathematics Library should be to provide a platform that engages the mathematical community in enriching the library's knowledge base and identifies connections in the data.**

### Search and Discovery

Mathematicians want to be able to understand mathematical objects— such as an equation, theorem, or hypothesis—more effectively and with greater ease. This quest can be aided by having the ability to specify a mathematical object either in natural language or more formal notation and get information on where other uses of the object appeared in the literature, definitions of the object, or related objects of interest. For example, consider questions of the form: "Given a hypothesis, what theorems involve this hypothesis?" or "Given a partial list of hypotheses and some conclusion, what additional hypotheses are known to imply the conclusion?"

The ability to ask and receive meaningful information about questions such as these is largely out of reach of current technology. It will require considerable research and investment to get even partway there. But the committee sees first steps toward realizing such capabilities in the innovative work of Wolfram|Alpha in the restricted domain of continued fractions.[15] Wolfram|Alpha prototyped and built a technological infrastructure for collecting, tagging, storing, and searching mathematical knowledge of continued fractions and presents it through a Wolfram|Alpha-like natural language interface. The main types of knowledge provided in this work are theorems, mathematical identities, definitions and concepts, algorithms, visualizations and interactive demonstrations, and references.

The committee believes there are many other subdomains within mathematics where significant advances on such very difficult problems may be possible with some mixture of modern methods of natural language processing and machine learning, expert human analysis of the literature of the subdomain (aided by computer), and knowledge representation approaches. Beyond hints of broad feasibility, the Wolfram|Alpha experience suggests the following:

---

[15] M. Trott and E.W. Weisstein, "Computational Knowledge of Continued Fractions," *WolframAlpha Blog*, May 16, 2013, http://blog.wolframalpha.com/2013/05/16/computational-knowledge-of-continued-fractions/.

- Key characteristics may be identified to make specific subdomains more feasible;
- It is possible to understand which of those subdomains are likely to be valuable to mathematicians, if they are appropriately captured and represented; and
- It is possible to understand how to encode knowledge so it is not specific to a single computing platform.

From here, one could imagine funded investments to encode specific mathematical subdomains in parallel to investment in work on the more general problem. Such subdomain-specific campaigns could be carried out as part of larger literature analysis efforts in the subdomain, which would build up or enrich the ontology and the link databases of the DML.

Intelligent information extraction and transfer are needed. For instance, it would be helpful if a user could just highlight a formula and then click on a button that submits the formula to a DML service that responds to some obvious questions, such as the following:

- Is this a well-known formula?
- Is it close to one in some curated list of formulas?
- Does it have a name? A homepage?
- Can it be parsed directly into a rigorous format for computation? If not, can the user be provided with some indications of the ambiguities encountered in parsing, and make choices as to which meaning is intended?

Moreover, it would be useful to be able to do this for more complex objects such as theorems and hypotheses. The committee does not wish to be too prescriptive about exactly how such capabilities and services might develop. In some specific domains, such as special functions and integer sequences, the necessary database of mathematical information is largely already constructed. The remaining issues are as follows:

- *Social*—Getting data to where they can be machine processed for development of services, and
- *Technical*—Building an adequate human-computer interface to enable users to interact with such databases in their everyday mathematical work.

The committee sees enormous potential for developments in this area by some concerted research effort involving a team of people with complementary expertise in machine learning, natural language processing, human-computer interaction, and mathematical knowledge representation. The key is to attract some sustained interest of people with relevant expertise in working together to produce DML capabilities of this kind.

## Searching for/in Mathematical Equations or Formulas

Mathematical papers are rife with complex equations and formulas; a few examples (from different mathematical subfields) are given in Figures 2-1, 2-2, and 2-3. It is clear that these expressions are a daunting

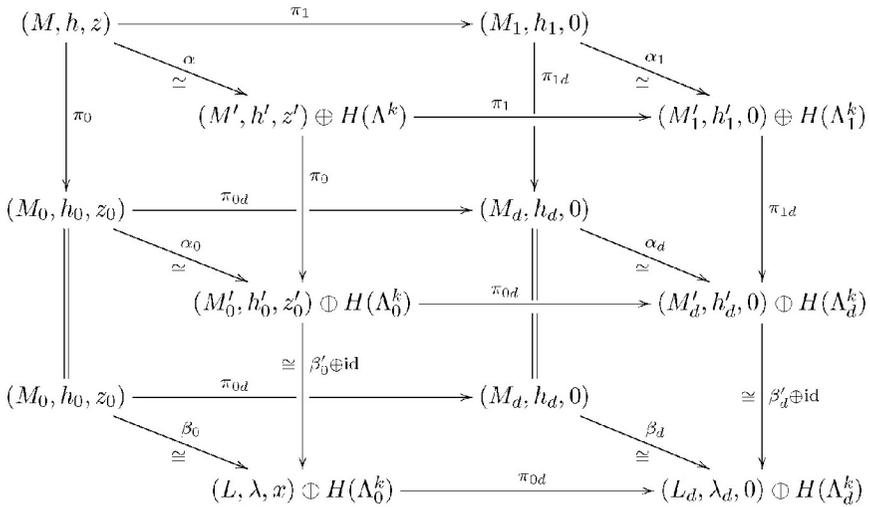

FIGURE 2-1 An example of complex mathematical typography. SOURCE: Lee and Wilczyński (1997).

$$\left| \sum_{\substack{\nu \in a + \frac{1}{2} + \mathbf{Z} \\ \nu \geq \nu_0}} \nu^{-1} e^{\pi i \nu^2 z + 2\pi i \nu (b + \frac{1}{2})} \right|$$

$$\leq \frac{e^{\pi i \nu_0^2 z}}{\nu_0} \frac{1}{\left| e^{2\pi i (b + \frac{1}{2})} - 1 \right|} + \sum_{\substack{\nu \in a + \frac{1}{2} + \mathbf{Z} \\ \nu \geq \nu_0}} \frac{1}{\left| e^{2\pi i (b + \frac{1}{2})} - 1 \right|} \left( \frac{e^{\pi i \nu^2 z}}{\nu} - \frac{e^{\pi i (\nu+1)^2 z}}{\nu + 1} \right)$$

$$= \frac{e^{\pi i \nu_0^2 z}}{\nu_0} \frac{2}{\left| e^{2\pi i (b + \frac{1}{2})} - 1 \right|} \leq \frac{2}{\nu_0 \left| e^{2\pi i (b + \frac{1}{2})} - 1 \right|}.$$

FIGURE 2-2 An example of complex mathematical typography. SOURCE: Zwegers (2008).

$$NG(q) - \overbrace{G \times \cdots \times G}^{q-times} \ni (h_1, \cdots, h_q) :$$
$$\text{face operators} \quad \varepsilon_i : NG(q) \to NG(q-1)$$

$$\varepsilon_i(h_1, \cdots, h_q) = \begin{cases} (h_2, \cdots, h_q) & i = 0 \\ (h_1, \cdots, h_i h_{i+1}, \cdots, h_q) & i = 1, \cdots, q-1 \\ (h_1, \cdots, h_{q-1}) & i = q \end{cases}$$

$$N\bar{G}(q) - \overbrace{G \times \cdots \times G}^{q+1-times} \ni (g_1, \cdots, g_{q+1}) :$$
$$\text{face operators} \quad \bar{\varepsilon}_i : N\bar{G}(q) \to N\bar{G}(q-1)$$

$$\varepsilon_i(g_1, \cdots, g_{q+1}) = (g_1, \cdots, g_i, g_{i+2}, \cdots, g_{q+1}) \qquad i = 0, 1, \cdots, q$$

We define $\gamma : N\bar{G} \to NG$ as $\gamma(g_1, \cdots, g_{q+1}) - (g_1 g_2^{-1}, \cdots, g_q g_{q+1}^{-1})$.

FIGURE 2-3 An example of complex mathematical typography. SOURCE: Suzuki (2013).

challenge to search or classify because of their complex notation, reliance on formatting not standard in other disciplines, and extreme precision. Even searching for notationally simpler concepts such as the historical evolution of an equation, whose name and notation may have been evolving for decades, is challenging. However, there exist several approaches that, separately or in combination, could be useful for this purpose. Many recognition tools are not based on a conceptual understanding of what one seeks to recognize, but rather on statistical analysis of usage patterns.

For instance, Google has a very good spellchecker ("Did you mean X?") that is based not on a deep conceptual understanding of the transformations possible on English words (i.e., morphology), but on statistical modeling of what humans typically type in as queries. For example, if one asks Google to search for "wavelette," it suggests "wavelet" as a more appropriate target, because lots of people have typed "wavelet" into the system in the past and only a few have typed "wavelette." To do this, the system needs to figure out that these two things should be connected, which is feasible given co-occurrences with other words in previous queries.

Another recognition example is Duolingo,[16] which uses a crowd-sourced model to translate content and vote on translations. Users of Duolingo sign up to learn and perfect their use of a foreign language; in the process, their feedback assists a project that produces high-quality translations of webpages. Originally, the Duolingo team expected to use insights from experts on how people learn languages. However, it is far easier and as effective to collect statistics on usage and let those guide the make-up of the successively more sophisticated tasks given to the users learning a language (von Ahn, 2011).

---

[16] DuoLingo, http://www.duolingo.com/, accessed January 16, 2014.

Because mathematical expressions are more complex than the word spellings handled by the Google spellchecker, the first step in characterizing them for clarification would be to bring them into, or describe them in, some form amenable to such statistical analysis. There exist several possible directions that could be followed, of which some are discussed below.

One approach would be to use feature vectors such as the Scale-Invariant Feature Transform (SIFT) to characterize mathematical expressions. Computer vision experts introduced SIFT for images, which can be used to "recognize" images or objects represented in images, independent of some warping, scale changes, or variations in illumination, by comparing feature vectors that keep track mainly of the statistics of many different types of local features, without trying to build a higher-level vision model of the objects depicted in the images (Lowe, 2004; Bosch et al., 2007). One may wonder whether it might be possible to likewise characterize some complicated mathematical expressions by a list of features that would then be useful in recognizing (variations of) it in other papers. Although SIFT features themselves are designed for image analysis, the general idea of attaching features to objects and then classifying them is extremely powerful. Objects can currently be retrieved or clustered based on words, tags, or descriptions in general (author, date, and so on). If properties of equations and theorems could be identified that were characteristic of their meaning, discriminating between different equations, and easy to detect, then techniques from information organization could be applied to help with browsing and searching.

A different approach would be to describe the way the formula is constructed geometrically, similar to the International Chemical Identifier (InChI)[17] for chemical substances. InChIs, which are nonproprietary identifiers for chemical substances that can be used in printed and electronic data sources, are specifically designed to enable easier linking and searching of diverse data compilations. InChIs can be computed from structural information and are human readable (with practice). Graphical representations of chemical substances are automatically converted into unique InChI labels, which can be created independently of any organization, built into any chemical structure drawing program, and created from any existing collection of chemical structures (Heller et al., 2013). It is conceivable that a similar automatic characterization could be done for mathematical formulas, which are two-dimensional graphical representations of the mathematical concepts they describe. If it were possible to determine international mathematical formula identifier labels, these could be extremely useful for searching and linking mathematical papers.

---

[17] International Union of Pure and Applied Chemistry, "The IUPAC International Chemical Identifier (InChI)," http://www.iupac.org/home/publications/e-resources/inchi.html, accessed January 16, 2014.

Yet another approach would be to create an ontology to characterize mathematical expression, similar to what Linnaeus created for biology. Linnaeus's *Systema Naturae* (Linnaeus, 1758) advocated the systematic and consistent use of a uniform botanical and zoological nomenclature in which each species is characterized by a two-word description: the first is the genus (which could, in principle, be decided for each specimen by examining just a few of its characteristics, even by people with limited expertise), and the second name specifies the exact species within the genus. In particular, Linnaeus's classification of plants was based on only their reproductive organs. By counting pistils and stamens of a plant, anyone, even without much botanical knowledge, could get a listing of genera that the plant in question should belong to. Prior to this system, botanical nomenclature was completely disorganized, and a botanist would describe a specimen by composing a long multiple-word descriptive name providing what seemed to be the list of most relevant and distinguishing characteristics, leading to different descriptive names for specimens of the same species. The two-tier system promulgated by Linnaeus provided a "tag" of which the first term (the name of the genus) could be relatively easily determined, followed by a search for the exact, lengthy description within that genus for a specific species, which made it possible to determine that the specimen was a new species. This system was first heavily opposed, but was gradually adopted by botanists.

While this system is currently under revision due to the widespread availability of genomic information, it is still a useful example of how ontologies can play an important role in establishing structures that promote the ease of comparing and searching. As compared to some scholarly domains, mathematics is fortunate to have a significant de facto standard ontology around which its research literature is organized. The Mathematical Subject Classification (MSC) standard ontology is used by both MathSciNet and zbMATH, meaning that more than 3 million research articles, chapters, and proceedings papers are already indexed using this scheme.

However, MSC is designed for indexing resources at the granularity of articles or conference proceeding papers. It exposes article topics in a broad way, but does not expose, for example, the theorems, functions, or sequences proven, used, or discussed in a paper. To support a broad range of discovery, browsing, and use cases, mathematicians need access at finer granularities, and they need descriptions created from viewpoints that go beyond the traditional library bibliographic cataloging perspective. In short, new ontologies will be required.

There is a foundation on which to base some of these supplemental ontologies. MathSciNet and zbMATH provide strong personal and journal name authority control back to the middle of the 20th century. The MathSciNet authors' database indexes more than 660,000 distinct authors. This database can be searched by name and/or affiliation with excellent recall and precision. In the context of linked data and the Semantic Web, this database also can be used (and has been used by MathSciNet) to generate such things as collaboration graphs, which are useful (for example) for calculating a "collaboration distance" between any two authors. This in turn provides another way to fine tune browsing and searching (e.g., searching for works by an individual or any other author within a certain collaboration distance of the original individual) under the assumption that such authors might have ideas that are related in some fresh way to the smaller collection of articles directly cited.

Other ontologies could be developed from resources that organize objects and concepts found within the mathematical literature. A particularly rich and promising example is the On-Line Encyclopedia of Integer Sequences (OEIS) (see Appendix C), which names and indexes more than 200,000 sequences. OEIS is particularly attractive for use as the foundation of an ontology for organizing the literature at a finer granularity, because it uses persistent identifiers for sequence entries and because most entries already include at least some bibliographic references.

Most retrieval from text is based on the words in the text; by contrast, many of the actual symbols in a mathematical equation are not useful search terms (e.g., a search for the indices $i$, $j$, and $n$ would not be useful). And strings of symbols are only partially useful as search terms; for example, the following three equations

$$f(x) = (x + 1)^2 \text{ and } g(t) = (t + 1)^2 \text{ and } h(s) = s^2 + 2s + 1$$

are all mathematically equivalent definitions of a quadratic function and should be considered so, even though the symbols have changed and elementary mathematical manipulations intervene. Similarly, the context of the variable $x$, $t$, or $s$ may itself change, indicating a real or complex number, or a quaternion, or perhaps a matrix or linear operator, although within a given work, the notation is usually gradually fixed so that, say, $x$ always represents a real variable. There are, however, certain conventions such as $\pi$ being often reserved for the irrational number and matrices usually being denoted by capital letters. Tensor analysis, in particular, utilizes strict conventions on the meanings and arrangement of subscripts and superscripts, as well as often employing the Einstein convention of implied summation on repeated indices.

As noted above, identifying the subfield of a mathematical paper may help disambiguate the notation. For instance, the (not uncommon) nota-tion π for a generic permutation of a finite number of items is less likely
to be confused with 3.141592… when the context is already identified as
pertaining to combinatorics rather than analysis or geometry. In practice,
it may well be that once sufficient statistics on usage have been collected,
such disambiguation could be done based on only statistical data.

Additional problems are posed by the historical literature: as a field
evolves, notations and terminology change, making connections to older
literature treating the same mathematical objects even more daunting.
Historically, notation conventions may vary where change is a reflection of
increased complexity and deeper understanding. In mathematics literature,
the gradual evolution in terminology and notation includes disputes that
usually (but not always) get resolved on what to call and how to represent
concepts, theorems, objects, etc. The evolution also reflects the integration
of work spanning many languages and cultures, each with their own idio-
syncrasies. Mapping back to earlier representations and concepts may not
be straightforward or direct.

In order to provide the careful typesetting modern mathematicians
require, precise typesetting and document preparation systems have been
developed, of which the most widely used is TeX,[18] together with its
descendants LaTeX, LaTeXe, etc.[19] TeX and its derived systems lead to
nicely typeset formulas (all the examples in Figures 3-1, 3-2, and 3-3 were
realized this way), and they have become an indispensable tool for math-
ematicians (most of whom do their own typesetting for papers they submit
for publication). At first sight, the LaTeX source code for a formula could
be thought a good candidate for an international mathematical formula
identifier. However, LaTeX is a presentation format, and equations in
LaTeX cannot be easily converted to a semantic representation that can be
used in other contexts. As a simple example of this problem, finding a string
in italics might mean, depending on the context and style, that it is a journal
title or a foreign word; to present the document in a different format or cre-
ate metadata, one needs to know the semantic significance underlying the
typographic display. Often, there is no one-to-one correspondence between
a mathematical formula as it appears on the printed page and the LaTeX
instructions leading to it; this nonuniqueness is even more pronounced if
one takes into account small variations in spacing (or changes of names
of variables, as illustrated above) that would not affect the reading of the
mathematical meaning of the formula by a mathematician. In this sense,
the LaTeX code for a formula would seem to fall short as a direct template
for a putative international mathematical formula identifier (as discussed earlier in this section). However, the National Institute for Standards and Technology Digital Library of Mathematical Functions[20] uses metadata embedded in the LaTeX code used to typeset the formulas to enable formula and notation search. This LaTeX metadata search, while not quite a LaTeX formula search, is fairly successful in dealing with dynamic notation and terminology change in the literature of special functions.

---

[18] "TeX," *Wikipedia*, last modified January 7, 2014, http://en.wikipedia.org/wiki/TeX.
[19] LaTeX—A document preparation system, last revised January 10, 2010, http://www.
latex-project.org/.

An option for semantic representation of mathematical formulas can be provided by MathML,[21] which allows for mathematics to be described for machine-to-machine communication and is formatted so that it can easily be displayed in webpages. There have already been some research efforts along the lines suggested above, and there are a limited number of both experimental and production systems available that involve some kind of formula search. In particular,

- There is some level of formula search in EuDML, using MIaS/WebMIaS (Math Indexer and Searcher),[22] a math-aware, full-text-based search engine developed by Petr Sojka and his group (Sojka and Líška, 2011).[23] An approach based on Presentation MathML using similarity of math subformulae is suggested and verified by implementing it as a math-aware search engine based on the state-of-the-art system Apache Lucene.[24]
- Some type of characterization of formulas is inherent to the searches underlying the Wolfram|Alpha engine. As part of a project in seeing whether mathematicians would find it useful to be able to search the literature for formulas, Michael Trott and Eric Weinstein of Wolfram|Alpha implemented some characterization of formulas for the research literature on continued fractions (essentially programming it "manually"). This small, fairly contained body of literature was chosen because most of the relevant papers are now in the public domain. However, the field of continued fractions is not very active at this point, and it may be hard to get a good sample basis of users to assess whether this search capability would lead mathematicians to new ways of using or searching the

---

[20] National Institute for Standards and Technology (NIST), Digital Library of Mathematical Functions, Version 1.0.6, release date May 6, 2013, http://dlmf.nist.gov/.

[21] W3C, "Math Home," updated November 26, 2013, http://www.w3.org/Math/, accessed January 16, 2014.

[22] EuDML@MU, "MIaS/WebMIaS," last change October 28, 2013, https://mir.fi.muni.cz/mias/.

[23] See also Petr Sojka's webpage at Masaryk University, Brno, last updated December 3, 2013, http://www.fi.muni.cz/usr/sojka/.

[24] Apache Software Foundation, "Apache Lucene Core," http://lucene.apache.org/core/, accessed January 16, 2014.

literature. It should be noted that the Mathematica-based formula characterization/search underlying Wolfram|Alpha is proprietary, in contrast to the completely nonproprietary nature of the InChI, which would also be desirable for an international mathematical formula identifier.

- Springer LaTeX Search[25] allows researchers to search for LaTeX-formatted equations in all of Springer's journals. In an issue of "Author Zone,"[26] Springer's eNewsletter for authors, Springer reveals that this free tool, which searches over a corpus of 120,000 Springer articles in mathematics and related fields, was created by 8 months of engineering a process that normalizes LaTeX equations. An open tool such as this would be valuable to the DML and to other mathematical indexing services.

**Finding: While fully automated recognition of mathematical concepts and ideas (e.g., theorems, proofs, sequences, groups) is not yet possible, significant benefit can be realized by utilizing existing scalable methods and algorithms to assist human agents in identifying important mathematical concepts contained in the research literature—even while fully automated recognition remains something to aspire to.**

## Navigation

Mathematicians want the ability to navigate and explore the corpus of mathematical documents available to them, be it through institutional library services or through free services. This goes well beyond accessing electronic versions of papers by following citations. The ability to click on an object in a document and be able to quickly find additional information about that object might help a mathematician decide whether to examine it further. Such additional information on an object might include the following:

- Other articles discussing the same object, or perhaps slightly more general or specific objects (and not necessarily with the same names);
- A description of when and where that object was first defined in the literature;
- A list of reference resources (textbooks, encyclopedia entries, survey articles) with information about the object; and

---

[25] Springer, LaTeX Search, http://www.latexsearch.com/, accessed January 16, 2014.

[26] Springer, "LaTeXSearch.com: Introducing the latest Springer eProduct in the field of Mathematics," http://www.springer.com/authors/author+zone?SGWID=0-168002-12-693906-0, accessed January 16, 2014.

- Different representations of the object (such as a LaTeX fragment or as Mathematica® code).

This is an area where it should be possible to make rapid progress, given a foundational DML investment in ontologies and links.

Improved navigation of the mathematics literature would enhance research capabilities in several ways. It would allow a researcher to find different resources and publications more easily and to find seemingly unrelated but relevant topics within the literature. It would also help a researcher to address the simply stated but inherently complex question, "Has this been done before?" Being able to answer this question would save valuable research time and simplify the problem-solving track, all while making the existing literature's structure more transparent and easy to use.

### The Citation Graph

Research articles can be viewed as the vertices in a large directed graph in which article A "points to" article B if A cites B. This citation graph is mostly tree-like: references are typically to older articles, although there are certainly cases of more or less contemporaneous articles that cite each other; some larger loops probably exist as well. Researchers interested in learning about a new direction or subject typically explore this graph; they start reading a particular research paper of interest and then climb back along the branches, reading some of its references and then some of the references of those papers, and so on. The creation of a citation index, as provided by MathSciNet within mathematics and by Google Scholar, Scopus,[27] and Web of Science across many more fields, allows the user to traverse the graph in the reverse direction, that is, to find for each paper all the articles that cite it. This very useful search tool makes it possible to easily find recent developments based on a paper of interest. Users would then be able to easily integrate or compare such information with whatever could be provided by other indexing services.

Making such comparisons or aggregations is at present very difficult. An expert user can do it in a few clicks by cutting and pasting from one browser window to another, but it is a few clicks for each resource, perhaps 12 clicks to compare returns from all three of these services. But with modern browser extension capabilities, such as those provided by Scholarometer,[28] which harvests data from Google Scholar, it is straightforward to write a dedicated browser extension for mathematical search and retrieval that would take a reference string from almost any source. The committee sees this kind of on-the-fly querying and aggregation of data from multiple services as the solution to the vexing compartmentalization problem for indexing services.

---

[27] Elsevier B.V., Scopus, http://www.scopus.com/home.url, accessed January 16, 2014.

[28] Indiana University, Scholarometer, http://scholarometer.indiana.edu/, accessed January 16, 2014.

A DML navigating tool could incorporate some mechanism for sorting and prioritizing the references it produces. A desirable feature of an open service is that such algorithms for ranking could be adjusted if so desired by the user, based on some special search criterion tailored by the searcher right then, or possibly influenced by the searcher's past preference behavior that is recorded by the system. Other basic questions that can be addressed by integration of DML data with data from various more-or-less-cooperative search service providers include the following:

- Which articles are cited in this paper? (This information is typically provided in the paper's list of references.)
- Which articles cite this paper? (This is a search that looks forwards in time, looking for papers that list this paper as a reference.)
- Which articles cite both papers A and B?
- Which articles are cited in both A and B?

Techniques for data analysis using methods such as bibliographic coupling and citation analysis are well established, and available software could be deployed for the benefit of DML users. A significant amount of citation data in mathematics and related fields is already more or less openly available from various open-access sources.

It should be possible to assemble accessible enhanced visualizations and graphical displays that capture features of a bibliographic data set that are not easy to find in a textual representation, and to make these features useful for search. Interactions between objects in a data set can be revealed by graphical displays within a browser (MacGillivray, 2013). Search results can be visualized in open formats, such as Scalable Vector Graphics (SVG),[29] and can be obtained from open search systems such as Lucene[30] or ElasticSearch.[31] Because today's widespread availability of all kinds of data is increasing attention on the need for better visualization tools, the committee anticipates that greatly improved open-source tools for graphical displays will become widely available and easily deployable to demonstrate interesting and novel features of the graphical relations in bibliographic data sets, not just those derived from citation graphs, but also those from collaboration graphs[32] and other graphs associated with relations between mathematical entities, such as implications or similarities.

---

[29] "Scalable Vector Graphics," *Wikipedia*, http://en.wikipedia.org/wiki/Scalable_Vector_Graphics, accessed January 16, 2014.

[30] Apache Software Foundation, "Welcome to Apache Lucene," http://lucene.apache.org/, accessed January 16, 2014.

[31] Elasticsearch, http://www.elasticsearch.org/, accessed January 16, 2014.

As more data about the citation and collaboration graphs in various disciplines have become available, they have also been used as a tool for ranking the impact of specific scholarly journals over time and have begun to be factored into the evaluation of individual researchers within the tenure and promotion process, where enthusiasm about their quantitative and "objective" nature has increasingly overcome very real concerns about their limitations and inaccuracies as a measure of the impact of a given scholar. A good deal of work has been done proposing various so-called alternative metrics ("alt-metrics")[33] for scholarly impact both at the article level and aggregated to characterize the contributions of a scholar (e.g., the h-index[34]). Analytics of these sorts are more likely to be useful to track topics than to measure the worth of theorems, journals, or individuals because often they are easy to manipulate and do not accurately reflect the community's view of importance (Arnold and Fowler, 2011; López-Cózar et al., 2013).

There is also real interest among working scholars in the possibility of tracking the evolution of these graphs (probably in conjunction with other data, such as popularity of articles) in order to help allocate precious reading time by identifying emergent, potentially high-impact or high-interest articles within or across specific subdisciplines, and a hope that article-based metrics can be developed to assist with this. The availability of citation and collaboration graph data, in combination with other information provided by the DML, would be an important step in advancing these research programs.

**Tracking Article-to-Article Reading**

Beyond simply exploring the citation graph, it may be desirable to obtain and exploit information about what other users of the DML have found useful as they explored the graph. For instance, what is the answer to the question, "Which articles did readers like, who are (like me) interested in A1, A2, and A3?" This way, one could find papers that do not specifically reference each other but concern the same topic. (This type of linking is a routine task, practiced by many online stores: "others who liked this also liked. . . .") It does, however, rely on a large user base to traverse the various graphs involved. Such a user base could be developed only with strong incentives for users to participate, such as superior navigation and search tools, so it is to be expected that such methods will be useful only late in DML development.

---

[32] Collaboration graphs are already attractively viewable on Microsoft Academic Search with the proprietary Microsoft Silverlight software.

[33] Altmetrics, "Altmetrics: A Manifesto," v 1.01, September 28, 2011, http://altmetrics.org/manifesto/.

[34] The *h-index* is an index that attempts to measure both the productivity and impact of the published work of a researcher based on the set of his/her most cited papers and the number of citations that they have received in other publications (Hirsch, 2005).

Recommender systems, like the one described in the previous paragraph, based on user tracking or ones based on "liking" a paper or topic within a system, are not new and are currently employed by Google Scholar and Elsevier, among others. They could also be developed within other information resources such as arXiv and MathSciNet.

These methods also raise privacy issues as users navigate a network of DML information. Concerns about privacy issues can often be addressed with customizable privacy settings (e.g., private navigation without login, public navigation with some anonymization of users, and possibly public navigation with public identity). It is important that the different models for maintaining user privacy are examined and assessed, and that a meaningful approach toward privacy be established for the DML.

Widely available machine learning algorithms can be used to predict the preference rating of as-yet-unseen articles by a customer for whom only a very partial profile is available, based on (often equally partial) profiles of other customers. A highly publicized recent success was achieved through the Netflix Prize competition in which Netflix "sought to substantially improve the accuracy of predictions about how much someone is going to enjoy a movie based on their movie preferences."[35] The final winning algorithm in that contest was an intelligent combination of strategies that alone produced insufficient improvement. This demonstrated that substantial progress can be achieved by combining different approaches that may be less spectacular when evaluated independently of one another. Such incremental improvements may not be very interesting from the perspective of machine learning research, but they are potentially useful in production applications of machine learning algorithms that the DML could provide.

### The Mathematical Concept Graph

Mathematical research can also be aided by considering mathematical objects other than papers, through exploration of their connections in a directed graph. For instance, in the answer to the question, Which theorems or papers use theorem T?, the different links would likely be references to classical results and to later improvements that were made since theorem T first appeared. The committee imagines both supervised and unsuper-vised learning approaches to these problems. In supervised learning, the machine starts from a list of known concepts, say functions or theorems, and then attempts to identify various instantiations of that concept. This is similar to automated library cataloging with a fixed structure of categories. Unsupervised learning is instead a process of clustering of instances—for example, deciding which theorems are essentially the same. At the level of LaTeX encoded formulas, some version of this capability, and a consequent search-and-discovery mechanism, is already achieved by Springer's LaTeX Search capability.

---

[35] Netflix, "Netflix Prize," http://www.netflixprize.com/, accessed January 16, 2014.

As further motivation for such efforts, which may be very challenging, the committee notes that Don Swanson identified useful public, yet undiscovered, knowledge in the biomedical domain by examining under-explored connections between clinical observations (Swanson, 1986, 1987). Despite efforts over the past few decades to automate the discovery of new scientific hypotheses based on literature analysis, insight from a human researcher is still needed. Ganiz et al. (2005) suggested that domains other than medicine should be explored. The committee believes that similar "literature discovery" methods could lead to interesting (and underexploited) connections between different mathematical fields or results.

## Visualization and Analytics

One way to help mathematicians learn from the large, complex, and rapidly growing and evolving literature base is to employ tools that are being developed to analyze data in a wide variety of settings, including both visualization tools and other analytical and statistical approaches. These tools could exploit the natural graphical structure of co-authorship and citation graphs and the relations among various kinds of mathematical objects and the parts of the literature that discuss these objects (as described in the previous section). The availability of an ontology for mathematical objects is important, and new tools are being developed that perform visualization guided by both an ontology and a set of data tagged according to the ontology (such as a collection of papers, or theorems, in a mathematical scenario). Note that in most cases, the committee expects that general-purpose graph analysis and visualization tools will be used, not tools developed by the DML.

The DML's role would be to help mathematicians find the right tools and ensure that data from the mathematical literature and knowledge base are available in forms and formats, and through interfaces, that make it easy to use these general purpose tools. Presumably, progress in this area would be quick, given the availability of the DML's underlying ontology and link collections, because it can build on other large investments that are under way already.

The committee does not expect the DML to be a contributor, but rather a testbed, for deploying methods for visualizing data. There are many widely deployed methods that can be applied to bibliographic data on the scale envisioned for the DML, which is modest compared to many big data projects. Microsoft Academic Search[36] already provides attractive displays of the collaboration graphs across its corpus using its proprietary Silverlight™ software. While open-source alternatives would be more attractive, either the DML or other agents could easily offer such displays over DML data as they are collected. This would provide an advantage over the quality of text data displays offered by the mathematical reviewing services. Similar displays could easily be provided for navigation and indication of relations between subjects at the level of MSC2010, which would greatly improve on past efforts.

### Computational Capabilities

The committee wishes to promote cooperation between the DML and computational service providers to allow users functionality, such as being able to cut a formula out of a mathematical document and paste it into a computing environment. This can already be done to some extent for simple formulas by cutting, massaging, and pasting a formula into Wolfram|Alpha, which uses natural language processing methods to match natural language queries with more formal knowledge representations.

The mathematics community uses a variety of simulation software—both numerical (such as Matlab,[37] Octave,[38] Python,[39] R,[40] Origin[41]) and symbolic (such as Maple,[42] Mathematica,[43] Sage[44]). Most software tools have different formatting requirements, and these would have to be taken into account when transporting formulas to and from them.

---

[36] Microsoft Academic Search, http://academic.research.microsoft.com/, accessed January 16, 2014.

[37] MathWorks, MATLAB, "Overview," http://www.mathworks.com/products/matlab/, accessed January 16, 2014.

[38] GNU Octave, http://www.gnu.org/software/octave/, accessed January 16, 2014.

[39] Python Software Foundation, "Python Programming Language—Official Website," http://www.python.org/, accessed January 16, 2014.

[40] R Project for Statistical Computing, http://www.r-project.org/, accessed January 16, 2014.

[41] OriginLab Corporation, "Origin," http://www.originlab.com/index.aspx?go=Products/Origin, accessed January 16, 2014.

[42] Maplesoft, "Maple 17," http://www.maplesoft.com/products/maple/, accessed January 16, 2014.

[43] Wolfram, Mathematica, http://www.wolfram.com/mathematica/, accessed January 16, 2014.

[44] Sagemath, homepage, http://www.sagemath.org/, accessed January 16, 2014.

**Recommendation: The Digital Mathematics Library should rely on citation indexing, community sourcing, and a combination of other computationally based methods for linking among articles, concepts, authors, and so on.**

## Other Useful Features

Application programming interfaces, which allow for add-on applications to be built by independent users and groups, are useful for experimentation with the processing of and understanding of mathematics. There are likely other tools that the DML could support that would be useful to the mathematics community.

For instance, there is still a need for a good pdf reader for mathematics. Most mathematicians still print out papers they really want to read, even if they own and mostly use an e-book reader for their other reading needs. When asked why they prefer reading mathematics from a print-out, researchers told the committee that they want to be able to flip back and forth, have difficulty concentrating on an electronic version, and miss the ability to annotate the paper with a pen or pencil. The DML could provide an environment to try out experimental readers.

Even prior to the existence of the DML, one could gain experience and better understanding of the feasibility and value of these technologies with the help of testbed platforms. These could serve as a framework for research programs to explore promising technologies and services, including extraction and identification of mathematical objects and applications of tagging or classification (including, perhaps, community-sourced approaches).

Experiments with structuring math knowledge into Wolfram|Alpha have been very promising and provocative. These are worth extending into other areas to gain additional understanding of effectiveness and limits. It would be of interest to select areas that are of active research interest. A key issue here, however, is understanding how to extend or share this beyond just Wolfram|Alpha and to make the investment reuseable in other settings.

# 3

# Issues to Be Addressed

## DEVELOPING PARTNERSHIPS

The ability of the proposed Digital Mathematics Library (DML) to foster and nurture a wide range of partnerships will be key to engendering and supporting the kinds of functionality and services envisioned.

Historically, many of the most important ontologies, taxonomies, and other knowledge organization systems and services in mathematics started as the research project of an individual mathematician or a small handful of mathematicians working in close collaboration. There are both social and technical challenges to establishing such partnerships.

Over time, the resources required for the individual researcher to maintain and grow these projects online can be unsustainable. Or, researchers reaching the end of their career find that they need to transfer responsibility for the knowledge organization systems they have built to someone else. Whereas in past generations it was often considered sufficient in such situations to simply instantiate a current snapshot of a mathematics ontology or taxonomy in static print form, today there is a need to sustain such knowledge organization systems online so that they can interact with the literature as it continues to be produced and so that that the knowledge organization systems can themselves continue to grow and be enhanced over time. These practical realities provide ready-made incentives for researchers to partner with a community-based entity like the DML.

Because the DML will be new, compared to existing society and commercial publishers, it would have to demonstrate that it is a worthy and stable long-term home for such research output. Firmly situating the DML within the established research mathematics community, demonstrating its commitment to openness, and adhering to technological best practices and established standards could help make the DML a new, but natural, home for this material.

The DML must be open to a range of partnerships of various degrees with mathematics researchers. In some cases, the coupling might be quite loose, with a researcher continuing to maintain and develop their knowledge organization service while it is simply used and leveraged by the DML in pursuit of the DML's broader mission. In other cases, the collaboration might be quite close, with the DML taking over (after a transition period) from the individual researcher the ongoing responsibility for an ontology or taxonomy. In the former case, adherence by both parties to a suitable application programming interface (API) will be essential.[1] In the latter case, the ability to map and automatically transform the original serialization of an ontology or taxonomy into a community-standard knowledge management serialization or encoding will be crucial.

For example, today the DML might reasonably adopt as one of its standards the Simple Knowledge Organization System (SKOS)[2] specification, an established, well-thought-of standard maintained by the World Wide Web Consortium for representing thesauri, classification schemes, subject heading lists, and taxonomies within the framework of the Semantic Web. However, few mathematicians starting out on a project to create a taxonomy of mathematical information objects pertinent to their research interests will begin by serializing their taxonomies in SKOS. Additionally, there are viable alternatives to SKOS with other vocabularies[3]—for example, ISO 25964, the international standard for thesauri and interoperability—or the metadata authority description schema (MADS[4]) promulgated by the Library of Congress.

Cognizant of the variety of ways to serialize ontologies and taxonomies and of the realities of how idiosyncratic serialization schemes adopted by individual researchers can be, it will be incumbent on the DML to work with partners to develop mappings and automated tools for transforming ontologies from one standard to another and from an idiosyncratic serialization into a standard like SKOS that is more suitable for long-term sustainability and growth. The DML will also have the challenge of doing such transformations in ways that do not foreclose on the opportunity for the partners involved (and others) to continue to help develop and refine their ontology or taxonomy.

---

[1] Currently, the Representational State Transfer (or RESTful) model of Web services enjoys broad consensus for this kind of scenario (see Pautasso et al., 2008). However, the correct approach will vary with time as technology and standards evolve.

[2] W3C Semantic Web, "SKOS Simple Knowledge Organization System—Home Page," last updated December 13, 2012, http://www.w3.org/2004/02/skos/.

[3] National Information Standards Organization, "ISO 25964—The international standard for thesauri and interoperability with other vocabularies," http://www.niso.org/schemas/iso25964/, accessed January 16, 2014.

[4] Library of Congress, "MADS Schema and Documentation," June 11, 2013, http://www.loc.gov/standards/mads/.

It is important that the DML engage members of the mathematics community from around the world. Many countries have made considerable investments in mathematical resources, and these investments should be captured, wherever possible, within the DML's outreach. There are challenges in engaging researchers across languages, but these should be addressed to the best of the DML's ability.

Partnerships with institutional entities (such as publishers and existing digital resources) are also crucial to the success of the DML. Here, the primary challenge is to be seen as complementary and enhancing, not competitive, while navigating constructive and effective partnerships with publishers, societies, Web services, and others, both specific to mathematics and those serving the much broader scholarly community. These entities control access to much of the mathematical literature under copyright. Only by establishing fruitful partnerships with such content providers and gate keepers can the DML encompass and link into and out of copyrighted scholarly literature. It is vital that users perceive that the DML is well-integrated with commercial services and commercially managed content. As described above, the committee envisions the resources, services, and tools offered by the DML as coexisting with, and often enhancing, the offerings from existing players in the mathematical information landscape.

The key to establishing such partnerships is perceived mutual benefit. Often such mutually beneficial agreements can be built around community-adopted standards and best practices. Patience may also be required, and it may be necessary to start with small agreements and collaborative undertakings. It is necessary to establish trust. Even small agreements can bear significant fruit. For example, the decision of the American Mathematical Society in 2002 to integrate OpenURL[5] into the version 8 release of MathSciNet has proven beneficial to MathSciNet, publishers, and MathSciNet users on campuses supporting OpenURL-based link resolvers. The longer-term goal, of course, is richer partnerships between publishers and services like *Wikipedia*, MathSciNet, MathOverflow, and zbMATH that would facilitate large-scale analytics, linking, and annotation.

Partnerships with academic institutions involved in the education of future research mathematicians will also be important. These should in-clude both departments of mathematics and academic libraries that serve members of these departments. The long-term sustainability of the DML is dependent on how its value is perceived by future mathematicians. As a distributed entity, many elements of the DML will almost certainly reside on campuses in mathematics departments and libraries. Additionally, alliances and partnerships with mathematics departments and academic libraries can facilitate partnerships with individual researchers and with publishers and others whose business models depend on subscriptions and memberships from academic libraries and departments.

---

[5] OpenURL is a standardized Web address format intended to enable Internet users to more easily find resources. OpenURL can be used with any kind of Internet resource but is most commonly used by libraries to connect users to subscription content.

Finally, there are rich opportunities for collaborations and partnerships with other departments and faculty within higher education, and even commercial partners that share common interests in underlying technologies and processing challenges. This would include computer science departments, schools of information, search engine developers, and others.

## ENGAGING THE MATHEMATICS COMMUNITY

As discussed throughout this report, it is essential that the DML engage the mathematics community as it works to cultivate and make sense of available mathematics knowledge. This report does not attempt to recommend how to do this but simply states that this is an important consideration for a future DML planning.

**Recommendation: Community engagement and the success of community-sourced efforts need to be continuously evaluated throughout DML development and operation to ensure that DML missions continue to align with community needs and that community engagement efforts are effective.**

Involvement of the mathematics community is being done well in a number of mathematics resources. MathOverflow does a particularly good job of providing a platform for individuals to post and respond to mathematics research questions. It rewards active users by granting them status and giving them access to additional features.

There is considerable skepticism among the mathematical community that it would be possible to encode the whole mathematical literature "by hand." It also remains an open question whether such functions can be automated, even in part. One approach to this may be developing a suitable community-sourcing algorithm, similar to the working of Duolingo[6] (see Chapter 3) or reCAPTCHA[7] (von Ahn et al., 2008). In both cases, users are asked to perform short identification or classification tasks in order to be given access to further material that they want to use or consult. Once a certain statistical consistency is achieved among several user responses for the same task, the identification or classification task is considered complete, and the result can be used (for translation of webpages in the case of Duolingo, or digital encoding of scanned documents in the case of reCAPTCHA). It might even be possible to design a "game with a purpose"[8] in which mathematicians worldwide would pair up to play entertaining games, the intermediate results of which would help recognize or characterize formulas in scanned texts, or that graduate students could use to help enhance their understanding of subject matter through chance collaborations and competitions with fellow students around the world.

---

[6] DuoLingo, http://www.duolingo.com/, accessed January 16, 2014.
[7] Google, "What Is reCAPTCHA," http://www.google.com/recaptcha/learnmore, accessed January 16, 2014.

If the DML is to be successful as a platform that enables mathematical users to access information and each other more easily in their pursuit of mathematical learning, then these users will be a huge resource to the DML. Like in *Wikipedia*, individual items such as papers, theorems, formula, comments, or open problems will be followed and maintained by volunteers. A large number of these volunteers will be students and researchers in mathematics or related fields. They could also play an important role in initiatives that mix community input and machine learning in order to provide useful tagging and links.

These and other models of community engagement should be assessed for the DML.

## MANAGING LARGE DATA SETS

From the perspective of modern data science, with data sets of petabyte scale, it is not a huge leap to move from dealing with millions of records of publications to hundreds of millions of locations of mathematical equations in the aggregated text of all mathematical documents with thousands of millions of occurrences of mathematical terms in that corpus. Based on MathSciNet by the Numbers,[9] at the time of publication there have been approximately 3 million articles and more than 696,000 authors from 1941. Traditionally, data sets of this size have been handled with relational database technology, with searches offered to users through a Web interface. More recently, it has become possible for such data sets to be usefully and easily manipulated using common technology. This means there is a dramatic increase in the potential for distributed users to contribute to substantial analysis and curation of bibliographic data sets on the order of magnitude of all mathematical books and articles.

---

[8] See von Ahn (2006) and Association of Computational Linguistics, "Games with a Purpose," last modified May 22, 2013, http://aclweb.org/aclwiki/index.php?title=Games_with_a_Purpose.

[9] American Mathematical Society, MathSciNet, "MathSciNet by the Numbers," http://www.ams.org/mathscinet/help/byTheNumbers.html, accessed January 17, 2014.

Emboldened by such technological progress, a consideration is to break mathematical papers down into their component parts, such as concepts, definitions, equations, theorems, proofs, etc., resulting in a much larger universe of mathematical artifacts, perhaps hundreds or thousands of millions of instances of these component parts. Appropriately deduplicated, these might amount to perhaps 10 million recognizable entities of mathematics, something for which one could imagine creating a webpage with pointers to at least some of its occurrences in the literature and capabilities for advanced searches and information retrieval at the level of mathematical entities rather than mathematical books and articles. A key point is that once the process of data mining some literature is done to identify mathematical entities, for example, by a process of unsupervised machine learning, these entities can be largely machine-generated but likely also manually curated for at least the most interesting and important of these entities. Managing large data sets often requires special considerations (NRC, 2013), some of which are discussed in this section.

### Dealing with Highly Distributed Data Sources

The base layer of mathematical publications is stored in a large number of widely distributed repositories owned and controlled by a variety of agents—commercial and academic publishers and various digital libraries (JSTOR, Project Euclid, arXiv, etc.)—as is the secondary indexing layer (zbMATH, MathSciNet, Google Scholar, Scirus, Microsoft Academic Search, CrossRef, etc.). Each of these sources has a distinct internal format. The European Digital Mathematics Library (EuDML) already has considerable experience in aggregation of both full text and metadata from diverse sources, and this experience should inform DML efforts.

### Tracking Data Provenance— From Data Generation Through Data Preparation

There are several distinct issues to consider as one moves into a complex digital ecosystem such as that characterizing the DML operating environment. One problem is technical and has to do with sourcing information that is aggregated, extracted, computed upon, and the like by the DML (or perhaps other services layered upon the DML services). In this case, one needs, most vitally, to be able to track where information came from; secondarily, there is a need to manage synchronization (but not always automatically preform such synchronization). If information is changed in some source repository, the DML may want to note that the information it is providing depends on an out-of-date version of the source data, unless the DML updates (recomputes) the information to reflect changes. At times it may be necessary to understand dependence and sensitivity—for example, if a given result turns out to be incorrect, what are the implications?

A second issue deals with permissions and legality and with scholarly norms of attribution. In the primary mathematical literature, citations to other papers, quotations from them, and sometimes reproduced figures, are legally covered in the United States under the doctrine of fair use, although occasionally also by explicit permission of the rights-holder of the cited material. As other types of digital uses become common, both scholarly norms of attribution and legal requirements must follow. Providing attribution is mechanical and largely covered under the source-tracing kind of provenance discussed earlier, although there are details about different levels of abstraction in cited objects that need to be sorted out, for example. From a more legalistic perspective, case-by-case analysis is needed; the first step is trying to make sure that there is enough information available to carry out the analysis, algorithmically whenever possible, due to scale and cost factors. For example, if it can be determined that sources are in the public domain, not subject to copyright, or are covered by certain kinds of well-known license (such as the Creative Commons series licenses), then much of the work is already done. In some cases, particularly involving articles published before digitization was anticipated (meaning that participation between author and publisher regarding rights to digitization are uncertain), various entities may have to explicitly give the DML permission to perform its content analysis and reuse computations. The DML will need to research and develop novel approaches to support these cases at scale.

For the primary publications work of mathematics, this problem is largely solved by widely adopted conventions of academic publication (providing authority through publication in peer-reviewed journals), the acknowledgement of primary sources through citation, and, more recently in the digital environment, the use of digital object identifiers and http links to point to sources.

Born-digital enhancements, such as the creation of derivative works from the existing base layer of book and journal data, will necessarily require indications of provenance, but the committee believes that this can be accomplished through open licensing.[10] For bibliographic data in the public domain, there is no legal requirement to acknowledge the source of a bibliographic item, although it is best academic practice to acknowledge its source, if only by a hyperlink.

---

[10] Two options include the Creative Commons Attribution—Share Alike License, which has been adopted widely and successfully by *Wikipedia* for user-contributed content such as annotations and reviews (http://en.wikipedia.org/wiki/Wikipedia:Text_of_Creative_Commons_Attribution-ShareAlike_3.0_Unported_License, last modified on May 13, 2013) and the Creative Commons Public Domain Dedication License, which is widely accepted as the appropriate license for large aggregations of bibliographic data (http://creativecommons.org/publicdomain/zero/1.0/, accessed January 16, 2014).

### Validating Data

Data validation is an important concern in any information management system but becomes especially important when aggregating multiple data sources together into a coherent knowledge base. The DML would face the challenge of addressing issues such as conflicting and incorrect data, incorrect tagging, and varying formatting syntaxes that can lead to confusion.

ChemSpider,[11] a free chemical structure database providing fast text and structure search access to more than 29 million structures from hundreds of data sources, faced data validation problems with its large aggregation of existing databases, data from peer-reviewed journals, and data provided through crowdsourced efforts. For example, each structure of chemical molecules can be described using a unique simplified molecular-input line-entry system (SMILES). However, Williams (2013) found that when the data from various sources were combined, there were instances where a unique chemical SMILES was being mapped incorrectly to multiple chemical structures. Although these inconsistencies had to be addressed through human intervention, the end result was a much more reliable database with fewer errors.

Two separate issues with validating bibliographic data need to be considered:

- The provision and maintenance of adequate schemas for the representation of mathematical bibliographic data records and the capability to check that the structure of a particular record is compliant with the schema; and
- The correctness or accuracy of particular data elements as they appear in a particular record.

Regarding the first issue, the committee expects multiple schemas for the representation of mathematical bibliographic records to coexist for a long time to come, due to a lack of heterogeneity of potential data sources and because normalizing records from different sources to confirm to a single schema would be an unnecessary cost. The challenges for the DML in utilizing metadata describing mathematical items will vary according to the type of resource being described. For instance, bibliographic metadata about library books are published in a relatively small number of schemas and are relatively consistent because of the large volume of standards and best practices published over the years by the Online Computer Library Center and most national libraries. Article-level bibliographic metadata for formally published mathematics are found in a greater variety of schemas. Augmentations to and annotations of book-level metadata, particularly in regard to digitized resources, may come from the Open Library, the EuDML, or other library or mathematics-specific community sources in some schema supported by that community.[12] These metadata inputs are even more diverse and less interoperable. The main issue for bulk processing of metadata is ensuring that every record is compliant and minimally complete to some schema and associated application profile and that every record clearly indicates the schema and profile to which it adheres. For nonbook resources especially, some resources may be required to normalize certain key properties (e.g., names). Once that is done, it is up to the processing service to achieve an acceptable level of interoperability across a modest number of schemas. The level of interoperability required varies according to service requirements. The committee anticipates that as the DML moves increasingly beyond formally published mathematics literature to also deal with nonbibliographic metadata, the resources needed for metadata remediation and higher levels of interoperability will grow. Again, community involvement in these processes will be critical.

---

[11] ChemSpider, http://www.chemspider.com/, accessed January 16, 2014.

Regarding the second issue, it is important that correctness and accuracy of data elements be monitored closely as bibliographic data acquisition and processing are undertaken by the DML. The DML collection will likely contain errors, and there should be procedures in place for users to flag these errors to draw the attention of qualified editors.

### Working with Different Data Formats and Structures

The different data formats and structures that the mathematical community finds useful for data representation will evolve over time.[13] The committee does not expect the DML to be an innovator in the field of data formats and structures, but rather to be an accommodator of the formats and structures that are widely accepted by the mathematical community and a facilitator of services for translating, when necessary, between formats.

---

[12] Examples include Marc records (see Library of Congress, "Understanding MARC," September 9, 2013, http://www.loc.gov/marc/umb/) or DublinCore (Dublin Core Metadata Initiative, http://dublincore.org/, accessed January 16, 2014).

[13] This evolution may start with legacy formats such as DublinCore, TeX and BibTeX, and progress through more advanced forms of XML including MathML, also JSON for lightweight Web services, and also incorporate formats from Mathematica and other mathematical programming languages to the extent possible.

## Ensuring Data Security

The committee sees some potential value in providing some user services that require login and storage of private data, such as for private annotations and/or the collection and mining of usage data, which might provide enhanced search and navigation features over the corpus. The committee is open to the possibility of including copyrighted data and extended metadata in the DML, with the aim of providing better services and linking to restricted-access content. These services would require enhanced data security. This would, however, impose a considerable administrative and legal burden on the organization managing the DML. Solution of this problem may depend on how monolithic or distributed the eventual DML architecture turns out to be. Having a safe, secure node in the system operated for the DML as a whole by one of the parties involved might be more feasible than having the parent DML entity responsible for it all.

## Developing Scalable and Incremental Algorithms

Literature-based data sets within mathematics are already large enough to provide some algorithmic challenges for tasks like clustering and deduplication. The problem of incremental processing is particularly important for a literature and knowledge base that continues to grow. Typically, some algorithm is applied to generate, say, a clustering or deduplication of a large data set. When new data come in, which might be recent publications or a newly digitized historical source, an update of the data processing is needed to incorporate the new data without reprocessing all of the data. This is particularly important if, subsequent to the original machine processing, there was some annotation or correction of data by human agents. Unless care is taken in managing workflows, there is the danger that these human contributions may be lost or overwritten in the reprocessing.

## Usage Tracking for Improvement and Diagnostics

Usage tracking refers to the process of capturing data on how a system is used and by whom. Such information is generally useful in identifying classes of users and their special needs, patterns of usage, beneficial workflows, underutilized areas of the system, and software bugs. Including technology to track such information would help to make the DML increasingly useful and would support diagnostics when users report errors. Types of usage tracking could include the number of times various sub-tools are used by a user during a session, the order of usage, and whether the system failed when a sub-tool was called. This usage data could then be aggregated to get system-wide usage by sub-tool and by pattern of activity. Usage data generally do not contain personally identifiable information; however, they may contain user class information—such as number of times accessed by novice, intermediate, or advanced users or related to classes of data providers, data searchers, and so on.

### User Security and Privacy Control

Systems that require users to register to use that system collect some personally identifiable information. This may just be name and email, but it can include contact information, location, and even financial information. This information is valuable from a system administration perspective because it can be used in a number of ways, from determining billing to identifying special needs by locations. In systems where the users can submit data, personal identifiers are also useful to limit the access of those who abuse the system and to provide recognition for those who provide high levels of valuable content. Such user tracking is thus of particular value when any part of the system employs contributions or community input. Keeping these data segmented from other data and not selling them or giving out user lists can preserve user privacy.

Another reason to have users register is to provide automated links to various social media systems and other online search systems, making it possible for the user to maintain a consistent user profile across tools. User desire for privacy can, in part, be maintained by having an opt-in system.

### Interoperability and Linkage to Social Networking Sites

Increasingly, scientists use social networking capabilities as a way to gather and vet data and ideas and as a way to identify and communicate with colleagues. Currently, a plethora of social networking sites are evolving independently. It may be prudent for the DML to let this functionality continue to evolve while supporting interoperability and linkage between various social networking sites to attain full functionality and to support broad usage styles even beyond what has been envisioned.

### OPEN ACCESS

The mathematical community has a limited capability to create and maintain a new information resource in an environment where a number of organizations, both commercial and noncommercial, have strong interests in owning, controlling, and profiting from the information and knowledge that potentially can be mined from mathematical publications. Scientific publishing as a whole seems to be at a crossroads regarding copyright. The committee foresees that the broad movement toward openness, mostly focused on open-access publications, open-source software, and open data, will likely encourage changes to the current copyright models used by many major publishers, as well as to scholarly practice and scholarly communication more broadly. While this report does not take a position with respect to publishing copyrights, the committee believes that all content created by and for the DML should be open to encourage the most buy-in from mathematicians and from potentially collaborating organizations.

The proposed DML organization could, for instance, oversee the creation and maintenance of a set of open resources—an ontology and collections of links—many of which rely on identification and extraction of objects or structures within the mathematical literature,[14] community input related to these objects, software used in mathematical research, and links to published literature. These object and link collections could be built up in large part by repeatedly computing over available collections of mathematical content. The initial DML creation and development will be challenging in terms of establishing the technological capabilities, engaging partners and the community, and planning for future growth. The insights about connections across the literature will be strengthened and become more useful. The process could begin with relatively open materials and willing partners; assuming that these services prove to add sufficient value, more holders of restricted-access materials may make arrangements to participate, and the net coverage of the mathematical literature would grow.

It is essential that the DML have access to and work well with all of the available mathematics literature, regardless of copyright status. While it might be tempting to build a system based on openly available material, such as mathematics heritage literature, the committee is convinced that the DML can be productive only if it has systematic input from and enthusiastic support by the mathematical community, which is unlikely to happen if the scope is restricted to open literature. In addition, it is envisioned that the DML computational services will be hospitable to new forms of mathematical scholarly communication (preprints, review papers, books, video material, etc.).

The committee is also cognizant of the current state of mathematics information resources and the systemic problems of compartmentalization, navigation, access, and maintenance. Briefly, compartmentalization—the partitioning of information and its maintenance by publisher or service provider—results in various agents having ownership and control of information and its maintenance, which can be sold to users as subscription services. Compartmentalization makes it difficult for users to navigate across boundaries, determine what information is accessible to them, and quickly access information. Unfortunately, the ability of services such as Google Search and *Wikipedia* to counter the compartmentalization problem are of only limited value in the highly structured discipline of mathematics, which requires structured information resources to provide better means of browsing and navigating the mathematical universe. Finding practical solutions to the challenges of compartmentalization, navigation, access, and maintenance—or at least compromises that allow progress—is the main challenge facing DML development.

---

[14] These mathematical objects and structures include a reference, keyword or phrase, theorem, proof, definition, equation, special function, conjecture, formula, transform, sequence, or symbol.

> **Recommendation: The Digital Mathematics Library should be open and built to cooperate with both researchers and existing services. In particular, the content (knowledge structures) of the library, at least for vocabularies, tags, and links, should also be open, although the library will link to both open and copyright-restricted literature.**

## MAINTENANCE

Many of the lists of mathematical objects described in Chapter 1 require expert and ongoing maintenance, and the DML needs to consider how to design its lists in such a way as to lessen their maintenance burden. With existing lists, it is often not clear how a user of the list can contribute new entries or edit existing ones. The problem is most obvious for lists published in copyrighted books, but it also exists for lists housed on other public sites. Some questions that need to be addressed are these: Who is responsible for maintaining this list? Is it a robot or a human? There is no established format or data schema for online publication of lists of mathematical objects, which complicates a machine's ability to read and reuse them. Rather, online representations of traditional print copies are prevalent. Often, and especially for lists contained in books, there are copyright restrictions that inhibit the process of maintenance, enlargement, enhancement, and reuse of these lists. Many of these lists do not provide links to primary or even secondary online sources.[15] Very few of these lists provide computable representations of the objects listed, such as code that can be passed to computing software.[16] These capabilities are important to mathematical research because merely knowing the formula is often insufficient; researchers also want to how it was proven, the history of the equation, and how it has evolved over time.

---

[15] *Wikipedia* supplies links where they have been provided and Online Encyclopedia of Integer Sequences (OEIS) does provide plain text references, but they are not generally hyperlinked.

[16] OEIS does provide both Mathematica and Maple code to generate most of its sequences and the Wolfram functions site offers Mathematica code for its basic functions, but not for any functional identities.

Experience to date with digital libraries and digital resources provides some insight and guidelines for how to approach the maintenance problems, specifically how to set up copyrights and licensing agreements, how to provide APIs, how to ensure that multiple copies of the information are always available, how to establish clear indications of provenance, and how to standardize and manage user contributions. These are fairly universal problems, and they should be amenable to fairly universal solutions with best practices provided by a central DML organization that is sensitive to the needs of the math community.

The maintenance strategy of the Online Encyclopedia of Integer Sequences (OEIS) seems particularly well suited for the DML. OEIS has developed a community of researchers in combinatorics who use it routinely in their research and who contribute to its maintenance. Essential here is the grass-roots nature of the effort. It was developed by one leading initiator, Neil Sloane, who had a vision of what could be done with a database of sequences and who gradually got people around him to contribute to it while enhancing the underlying software and functionality. The resource was developed in direct response to the interests and needs of a research community (and also with considerable interest from a larger community focused on recreational mathematics and pedagogy), and it was kept free and open, which engaged the community.

Another resource with similar communities of contributors/maintainers is Research Papers in Economics (RePEc). This is more of a traditional bibliographic resource than a database of entities, but the principles are very similar: find a way to make it easy for experts to contribute their domain knowledge and build up a knowledge base.

Community information projects often require both an inspired creator, often unrewarded at the start, and eventually a transition to a paid staff after the work grows beyond the capacity of an individual, even an individual assisted by a crowd-sourced effort. For example, arXiv was started by Paul Ginsparg alone at Los Alamos National Security Laboratory but is now run by the Cornell University Library. Ginsparg is still very active and involved in policy, but he cannot personally make every decision of the form, "Does this paper belong in cs.DL or cs.CY?" The Internet Archive similarly has a visionary, Brewster Kahle, founder and still in charge, but it also has a paid staff to keep operations going.

 Once the resource gets large enough to be of substantial value to the community, it has to be legally constituted to avoid issues of ownership and control. The use of the Creative Commons license[17] is an approach that the committee believes would work well for the DML. OEIS uses the

---

[17] Creative Commons, "About The Licenses," http://creativecommons.org/licenses/, accessed January 16, 2014.

Creative Commons Attribution-NonCommercial (CC BY NC) license,[18] but the DML should also consider other options such as the Creative Commons Attribution-ShareAlike (CC BY SA).[19]

In order for the DML to successfully maintain a database resource, it has to deal with the technical and human components. On the technical side, the DML has to provide adequate version-control and editorial software (similar to Wikimedia) to manage the deposit, editing, and cross-linking of documents. It is essential that this software work well and be kept up to date and well adapted to the current information environment. Some centralization of this activity seems beneficial. On the human side, the DML has to motivate people to contribute to the parts of the effort that are not easily or fully automated. One way to do that is to provide nice software that does the boring parts for them easily and allows them to focus on the parts where their expertise is really needed. Many database maintainers try to build and customize this sort of software for themselves, but then they get overwhelmed by the issue of software maintenance and spend more time on trying to deal with that than they do with contributing their domain expertise to the database. The DML could provide out-of-the-box software (or a Web service solution) for each math sub-community to curate its own material for benefit of a larger audience. The DML software would include mathematical knowledge, so that it could display properly formatted theorems and recognize structural similarities, often not possible in the numerous existing collaborative software offerings.

If the DML can provide a good software solution for managing mathematical entities, and deal with the management of that software in a central way, it can provide something that a large number of different mathematical communities could adapt for their own purposes, hopefully maintaining some centrally supported capabilities (version control, linking, math display, search, etc. ) without each sub-community having to solve these problems separately. At the very least, having some common standards for data exchange and interoperability, and some common reliable components for which there was some central support, would lessen the maintenance problem.

Some of the maintenance of the DML lists may be automated as well. The key is to find a balance between automated data mining of the literature and human annotation and curation. More work and experimentation is needed to develop editorial systems to assist this process. The main goal is to provide good tools to do largely successful cleaning and reduction of data before bringing portions of them to the attention of domain experts, whose time is limited, or possibly crowdsourcing less demanding tasks. Tools like Google Refine[20] and flexible, faceted displays of bibliographic data like BibServer[21] are very useful for this.

---

[18] Creative Commons, "Attribution-NonCommercial 2.0 Generic," http://creativecommons.org/licenses/by-nc/2.0/, accessed January 16, 2014.

[19] Creative Commons, "Attribution-ShareAlike 2.0 Generic," http://creativecommons.org/licenses/by-sa/2.0/, accessed January 16, 2014.

Both Google Scholar[22] and Microsoft Academic Search[23] do a huge amount of fully automated data processing of general academic bibliographic data. The methods behind these services could undoubtedly be brought to bear on more specialized data mining and data structuring tasks of the kind relevant to text mining the mathematical literature for formulas and the like. LaTeX Search[24] (Springer's free formula search) provides a step in this direction by allowing users to locate and view equations containing specific LaTeX code, equations similar to another LaTeX string, equations belonging to a specific digital object identifier, and equations belonging to an article or articles with a particular word or phrase in the title.

The DML will also have to develop in such a way as to learn from and complement the broader data conservation and data preservation movement, helping to organize and preserve the mathematical information it contains. It may be beneficial to cooperate with groups such as LOCKSS,[25] Portico,[26] or *HathiTrust*,[27] which do digital preservation today, and coordinate with projects such as the Data Conservancy,[28] DSpace,[29] and the linked open data movement, which are laying the groundwork for more powerful preservation techniques in the future.

---

[20] Google-refine, https://code.google.com/p/google-refine/, accessed January 16, 2014.

[21] GitHub, https://github.com/okfn/bibserver, accessed January 16, 2014.

[22] Google Scholar, http://scholar.google.com/, accessed January 16, 2014.

[23] Microsoft Academic Search, http://academic.research.microsoft.com/, accessed January 16, 2014.

[24] Springer, LaTeX Search, http://www.latexsearch.com/, accessed January 16, 2014.

[25] LOCKSS: Lots of Copies Keep Stuff Safe, http://www.lockss.org/, accessed January 16, 2014.

[26] Portico, "Services," http://www.portico.org/digital-preservation/, accessed January 16, 2014.

[27] *HathiTrust Digital Library*, http://www.hathitrust.org/, accessed January 16, 2014.

[28] DataConservancy, http://dataconservancy.org/, accessed January 16, 2014.

[29] Duraspace, DSpace, http://www.dspace.org/, accessed January 16, 2014.

# 4

# Strategic Plan

This chapter proposes a strategic plan for incremental and modular development of the Digital Mathematics Library (DML), with the aim of providing the mathematical community with at least some of the specific capabilities described in Chapter 3, as well as some further capabilities that should follow as corollaries of the basic development. The committee's strategic plan contains the following elements:

- Fundamental principles of the DML vision;
- Constitution of a nonprofit organization committed to development of the DML collection and services, called the DML organization;
- Initial development;
- Priorities for collections and service development;
- Technical considerations; and
- Resources needed.

Each of these elements is discussed in detail in the following sections.

## FUNDAMENTAL PRINCIPLES

The committee envisions the next step in advancing mathematics to go beyond traditional mathematical publications and take advantage of the mathematical information and knowledge stored in those publications to create a network of information that can be easily explored and manipulated. There is a compelling argument that through a combination of machine learning methods and editorial effort by both paid and vol-unteer editors, a significant portion of the information and knowledge in the global mathematical corpus could be made available to researchers as linked open data through the DML. The DML would help index and make discoverable collections of information created and maintained by distributed editors and specialized machine agents—much as Google now indexes and makes available information drawn from across the Web—but without the centralized processing and caching. But the DML would also need to engage substantial editorial input from the mathematical community. The DML would afford functionalities and services over the aggregated information, including capabilities for searching, browsing, navigating, linking, computing, and visualizing and analyzing, over both copyrighted and openly licensed content.

Some, but by no means all, of the proposed additional services and knowledge management utilities will rely on analysis of full content, done in a coordinated fashion. Other services will rely on analysis of metadata, which are often accessible with fewer or no restrictions. The committee feels that today—through reliance on a broad, distributed community, adherence to emerging standards and best practices, the use of new distributed collaboration and editing workflow models, and reliance on the affordances of emerging technologies such as linked open data and machine learning methods—these content and metadata analyses can be accomplished successfully in a distributed fashion—that is, without having to acquire, process, or store the entire universe of all mathematics publications centrally. While the approach outlined would require the central (or at least centrally coordinated) maintenance of key concept vocabularies and ontologies, large-scale, centralized processing and storage of mathematical publications would not be necessary.

The committee has identified a compelling opportunity for the following:

- The DML as a large, open collection of mathematical bibliographic information and mathematical concepts (e.g., axioms, definitions, theorems, proofs, formulas, equations, numbers, sets, functions) and objects (e.g., groups, rings) aggregated from diverse sources;
- Integrating and organizing the DML with existing repositories of publications and with indexing and computing services (as discussed in the Chapter 3 section on "Developing Partnerships");
- Encouraging, facilitating, and supporting the development and promulgation of novel Web and desktop services, including annotation, collection, and collaboration tools, and tools for search and literature-based discovery, that can be utilized within the DML;
- Supporting experimental and production applications of machine learning methods for the extraction of various mathematical entities, including topics, formulas, equations, and theorems, by data

    mining and large-scale data analysis of suitable portions of the mathematical corpus; and

- Supporting a combination of community input and traditional editorial workflows for validation of outputs of such machine processing and contribution of such outputs to the DML.

    The committee believes that it is necessary for the people and organizations involved in the DML to adopt some basic principles to guide the DML to reach its full potential.

### Adherence to Best Practices and Standards

    The proposed DML would benefit from adhering to broad technical standards and built-in interoperability, both for encouraging partnerships and taking advantage of non-mathematics-specific Web technologies that become available (Aalbersberg and Kähler, 2011; Gill and Miller, 2002). The DML would benefit from being developed with a modular architecture, allowing various technical development efforts to proceed in parallel with minimal coordination. There should be some initial agreements in principle about the nature of inputs and outputs of various components and Web services. One illustration of the importance of technical standards in mathematics is the value of Tex (and LaTeX), which standardized mathematical typesetting and revolutionized research mathematics publications.

    The DML architecture should adhere as much as possible to contemporary and evolving Web architecture standards for all its services, especially the standards of linked open data for publishing structured data on the Web so that it can be interlinked and become more useful. Linked open data allow a webpage to dynamically pull relevant information from related websites. For example, a website that displays local weather could pull information from an unrelated local traffic monitoring site to alert users to delays or road conditions, and it could pull from the local school district's website to alert users about potential closures. Linked open data are particularly valuable within the proposed DML because much of the value of the information comes from its connections with outside existing data. If these connections can be strengthened, the network of mathematical information will solidify, providing a clearer picture of the realm of mathematical research.

    The DML, as proposed, would not be collecting large amounts of copyrighted material; however, it would be amassing its own data collection of connections and understanding of mathematical information. These data (i.e., vocabularies, ontologies, annotations) and the DML-developed/supported software would benefit from being open source so that other researchers and developers could build upon it. The DML would need to

respect and recognize copyright limitations and work with publishers to make sure these can stay in place even while having minimal impact on the ability of users to discover and learn about resources.

The DML would also benefit from adhering to accepted norms for citations and evaluations. This may take the form of systematic application and support of the San Francisco Declaration on Research Assessment (American Society for Cell Biology, 2012) about emerging practices related to the evaluation of research articles.

> **Recommendation: The Digital Mathematics Library should serve as a nexus for the coordination of research and research outcomes, including community endorsements, and encourage best practices to facilitate knowledge management in research mathematics.**

### Competition and Cooperation with Other Organizations

To the greatest extent consistent with its goals and principles, the DML should seek to cooperate with and not to compete with existing information services and communication and desktop tools that are widely used by the mathematical community. Cooperation would include the following:

- Agreements on the structure of suitable data schemas for representation of bibliographic and mathematical information, including standards for representation of mathematics on the Web (MathML, MathJaX, etc.);
- Agreements on systematic use of identifiers and openly accessible Web services supported by other organizations (e.g., DOIs, Handles,[1] ORCIDs, MR and ZMATH identifiers, OCLC identifiers) instead of replication of these identifiers and associated services by the DML;
- Provision of agreements and conversion services, as needed, to ensure metadata interoperability and aggregation of data from various services; and
- Support for interfaces between the DML and existing information resources listed in Appendix C—for example, bibliographical, encyclopedic, content, social environments.

This cooperation also applies to arXiv, *Wikipedia*, MathSciNet, zbMATH, Google, Microsoft, and the general abstracting and indexing services, as well as to various companies with proprietary interests in mathematical communication and computation whose products the DML should seek to

---

[1] Handle System, http://www.handle.net/, accessed January 16, 2014.

enhance and make more openly accessible and reusable. This list of companies includes the following:

- Springer (with large amounts of mathematical information in SpringerLink[2] and its proprietary LaTeX search);
- Wolfram (with large amounts of mathematical information embedded in Mathematica and Wolfram|Alpha);
- Elsevier; and
- Maplesoft, a subsidiary of Cybernet Systems Co. Ltd. in Japan and a provider of software tools for engineering, science, and mathematics, especially Maple,[3] a powerful mathematical computation engine.

In areas where data standards are well established, such as for basic bibliographic data elements, such cooperation may be achieved by the DML organization with different data sources and services individually. For more complex data objects, especially those representing mathematical concepts, a community process, such as those commonly conducted by the World Wide Web Consortium,[4] should be involved in the selection and adoption of data standards by the DML. It is recognized that such data standards may typically start as ad hoc standards that eventually become codified and formalized through widespread use (e.g., Microformats Wiki[5]). The committee recognizes that some existing agents may be reluctant to cooperate with the DML in either development of data schemas, sharing of data, or both. In those cases, the DML should not allocate administrative effort on negotiating cooperation but rather find alternative agents who are willing to cooperate in providing the needed data or services in a manner consistent with DML principles.

### Collection from Diverse Sources

The DML should commit to support curation and management of mathematical information from diverse sources and facilitate access to mathematical information even though the sources are stored in different organizations. Similarly, CrossRef[6] currently tells users how to find items from different vendors, LOCKSS[7] manages shared storage across libraries, and ORCID[8] helps identify authors across publications. In particular, the DML should aim to acquire and process the following:

- Previously unindexed or partially indexed information about mathematical publications—including traditional journal papers, books, and other electronic resources—and their contents, such as their reference lists, names of their sections or chapters (table-of-contents data), their formulas, equations, theorems, and conjectures;
- Information relating to the relations of such data elements within various publications and the relations of these elements to various standardized lists of such elements; and
- Information from mathematicians' homepages, blogs, and discussion forums.

---

[2] Springer Link, http://link.springer.com/, accessed January 16, 2014.

[3] Maplesoft, "Maple 17," http://www.maplesoft.com/products/Maple/, accessed January 16, 2014.

[4] W3C, http://www.w3.org/, accessed January 16, 2014.

[5] Microformats, "The Microformats Process," last modified April 28, 2013, http://microformats.org/wiki/process.

[6] Crossref, http://www.crossref.org/, accessed January 16, 2014.

The DML should accept inputs of such data from all sources, commercial and noncommercial, subject only to copyright and licensing requirements indicated earlier, the judgment of DML-appointed editors that the material is suitable for inclusion in the DML, and the resources to process the data for ingestion into the DML. In particular, the DML should invite contributions of such content from both copyrighted and open-access sources. In all cases, the DML should commit to appropriate acknowledgement of the source and to inclusion of agreed indications of provenance in its data records.

### Support for Multiple Formats, Conversion Tools, and Best Practices

In many instances the cost of negotiating cooperation in schema standards may greatly exceed the potential reward of doing so. In such cases, it will be best for the DML to move ahead with lowest-common-denominator standards that are good enough for most applications and to which it is possible to map data from multiple alternative formats. Current examples of such standards are BibTeX, or slight enhancements thereof like BibJSON and BibXML, to which it is possible to map almost any reference text string that can be recognized as such by a human. A somewhat higher standard is provided by the European Digital Mathematics Library (EuDML) metadata schema specification[9] for typical mathematical article metadata supplied by a cooperative publisher. The DML should research and support multiple tools and services for the acquisition of data in diverse native formats and its conversion to higher-quality bibliographic formats such as those mentioned above. It should also provide guidance for best practices in managing various data and metadata formats and support basic communication spaces, such as an email list or help desk for data managers encountering issues in cleaning and converting diverse data sets of interest to the DML community. Examples of conversion tools for bibliographic data that are already very useful, although relatively unknown, are pdftotext,[10] MREF,[11] EJP-ECP Reference List Formatter,[12] inSPIRE-HELP,[13] BibSonomy Scapers,[14] Google Refine,[15] and Beautiful Soup.[16]

---

[7] LOCKSS: Lots of Copies Keep Stuff Safe, http://www.lockss.org/, accessed January 16, 2014.

[8] ORCID, http://orcid.org/, accessed January 16, 2014.

[9] European Digital Mathematics Library, EuDML Metadata Schema Specification (v2.0-final), https://project.eudml.org/eudml-metadata-schema-specification-v20-final, accessed January 16, 2014.

The creation of such data-conversion tools is typically a fairly straightforward programming task in which the difficulty depends on the complexity of the tool. However, such tools and their derivatives do impose a progressive maintenance burden to keep them compliant with changing data formats and expectations for both inputs and outputs, and with new versions of underlying software libraries and implementations. But the maintenance of such low-cost, high-reward data conversion and cleaning services, or links to the best maintained of these services and documentation of how to use them for DML purposes, is among the things the DML should commit to supporting.

### Flexibility and Extensibility of Schemas and Services

Recognizing the systemic compartmentalization problems caused by traditional database schemas and implementations, all DML schemas should adhere to current and emerging best principles of flexibility and extensibility. In particular, DML architecture should allow and encourage the following:

- Inclusion of data in a virtual collection from an essentially unlimited number of disparate and distributed resources of greatly varying

---

> sizes. Examples would include data stored on individual webpages and marked up with information, as is done with CoINS,[17] the emerging standards of schema.org or similar math-specific standards that might be developed by the DML community, or data available from various data providers via application programming interfaces (APIs) or periodic data dumps; and

- Creation of new features, tools, and services over DML data by individual and organizational participants, such as those outlined in Chapter 3, or by yet unimagined services that will develop in the future.

### Relation of the DML to Computer Algebra Systems and Formalization of Mathematics

There is a community of mathematical knowledge management, built largely around the development of formal theorem provers and reasoners (Carette and Farmer, 2009).[18,19] This community proposed an ambitious program of formalization of mathematics, following earlier efforts by Whitehead and Russel (1910, 1912, 1913), Hilbert's program,[20] and others. Some notable successes of this school are computer automated proofs of a number of important mathematical theorems, such as the famous four-color theorem. The committee anticipates further advances in this field, and perhaps some eventual synthesis of computer algebra systems (Mathematica, Maple, Sage, etc.) with the theorem provers. However, progress in this area has been slow, and there are deep cultural impediments, principally the fact that the dominant computer algebra systems are proprietary and likely to remain so for the foreseeable future.

### Summary of Principles

Consistent application of the principles in this section to the representation of mathematical information and conceptual knowledge in the World Wide Web will enable the mathematical community to achieve the most effective instantiation of the DML as an openly navigable representation of the universe of mathematical concepts, formulas, and relations. To achieve this, the DML would be just as accessible to human users as *Wikipedia* is today, with the same open license for text contributions and a public domain license for bibliographic and mathematical facts; it would be properly structured for machine access and reuse in discovery services; and it would be connected directly, through desktop software and Web services, to the mathematical research literature, current and future abstracting and indexing services, computational services such as Wolfram|Alpha, and desktop programs such as Mathematica, Maple, and Sage.

---

[17] OpenURL COinS: A Convention to Embed Bibliographic Metadata in HTML, Stable Version 1.0, http://ocoins.info/, accessed January 16, 2014.

[18] Mizar Home Page, last modified January 8, 2014, http://mizar.org/.

[19] Coq Proof Assistant, http://coq.inria.fr/, accessed January 16, 2014.

[20] "Hilbert's Program," *Wikipedia*, last modified January 3, 2014, http://en.wikipedia.org/wiki/Hilbert%27s_program.

## CONSTITUTION OF THE DIGITAL MATHEMATICS LIBRARY ORGANIZATION

The first step in this process is creating an organization that can manage and encourage the creation of a knowledge-based library of mathematical concepts and advocate for the needs of the mathematical community. The committee believes the DML effort would benefit from being spearheaded by a small centralized agent to avoid the project failing because of competing time commitments of its founders, which has happened in several cases mentioned in Appendix C. It is hoped that the DML can reach beyond this initial startup hurdle and ultimately succeed because of its core of dedicated staff, collaborators, and funders, and to ultimately create a strong, stable, and meaningful resource that is worthy of continued investment from the mathematics community.

> **Recommendation: A Digital Mathematics Library organization should be created to manage and encourage the creation of a knowledge-based library of mathematical concepts such as theorems and proofs.**

> **Recommendation: The Digital Mathematics Library organization should be an advocate for the mathematics community and help develop plans for development and funding of open information systems of use to mathematicians.**

The DML organization would benefit from being a small organization with minimal central agency and control. It is also important that the DML be able to operate in an environment of much larger organizations with big budgets and capability for sustained legal actions to achieve their ends. To survive as a small operation in a big information universe, it is important that the DML be organizationally nimble, quick to initiate pilot projects, and generally quick to learn from the experiences of both successful and unsuccessful efforts, both its own and those of others aiming to develop domain-specific knowledge bases. The DML could be largely reliant on other organizations to provide hosting for such organizational essentials as

- Basic computing and networking infrastructure, support, and services;
- Archiving (to be achieved in collaboration with existing scientific data and library archiving organizations); and
- Office space and administrative and support services of all kinds.

Management overhead can be minimized, for example, by making the executive director of the DML an employee of a supporting institution, most likely a major university library, whose time is funded either completely or in large part by a grant from the initial DML funder to that university. A modest number of initial staff positions could be funded similarly. This could be a good approach for the DML because many of the technical skills it will need are specialized and may be needed only on a part-time or fluctuating basis as various projects are taken on by the DML. The DML could at least initially avoid the management responsibility of having a large number of employees, but rather work on a contractual basis with staff employed by a variety of partner organizations with a commitment to various aspects of the DML effort.

The DML organization may also wish to consider other names before finalizing its constitution, both for itself as an organization, and for the collection and services it plans to create. One of the early administrative efforts of the DML organization would be to evaluate a number of legal and economic considerations involving branding and trademarks related to the choice of name. The committee envisions the DML organization as a coalition of member partners with commitment to the DML concept—the creation of a substantial digital representation of an open collection of mathematical information and knowledge—and to the DML development principles. The DML organization could be governed by the mathematical sciences community through an organization such as the International Mathematics Union (IMU) and, thence, through the member organizations of that union.

The DML constitution can support the general principles outlined above by including the following elements:

- Acquire and maintain a collection of digital representations of mathematical objects (e.g., theorems, functions, sequences) in machine processable formats;
- Advance mathematics by provision of useful information services over the collection;
- Maintain the DML data collection with stable URLs, an underlying Web-based open architecture, and APIs so new tools can be contributed, linked, and shared;
- Support development of a large community of users who will also help curate and contribute to the collection and its services;

- Support a community of developers of tools and services over the collection; and
- Collaborate with publishers and information providers to provide superior mathematical and information services built over the collection.

Governance of the DML could be overseen by an organization such as the IMU, with invitations to representatives from partner organizations. Initial funding of the DML for a 10-year period would be beneficial, during which long-term models for sustainable operations could be examined. The DML may benefit from including as many of the relevant organizations as are willing to participate. Some examples include MathSciNet (American Mathematical Society), the Society of Industrial and Applied Mathematics, the International Council for Industrial and Applied Mathematics, the European Mathematical Society, the Cornell University Library, Fiz Karlsruhe/Springer, Wolfram, MicroSoft, Google, *Wikipedia*, OEIS, EuDML, Elsevier, and Thomson Reuters. Publishers and volunteers will see the DML as more accurate and more tailored than other services and should recognize the gains possible from a coordinated approach to merging mathematical knowledge. As the DML grows, the community will accord respect to the volunteers who help build it. To protect itself from legal obligations regarding copyright infringement, the DML could consider a variety of approaches, including not claiming copyright on any DML material and requiring of contributions to be licensed by the contributor, or using a creative commons license.

The first step to confirm feasibility of this DML concept is to announce a proposal to the community, confirm that enough parties are willing to participate in the DML by contribution of data, expertise, or services to make the project viable, and, if so, support a meeting to resolve a basic constitution for the organization to establish its legal status in a suitable location.

## INITIAL DEVELOPMENT

Initial development of the DML would benefit from focusing on recruiting partners with potential data sources and resources, beginning a collection of mathematical entities to achieve some of the desired capabilities described in Chapter 3, and providing a foundational platform on which most of these capabilities might imaginably be achieved in a 10- or 20-year time frame.

The committee sees value in separate groups working on the technological infrastructure and on the administration of these projects, because they require different kinds of technical expertise, community input, and project management for their success.

## Recruiting Partners

The DML cultivation of partnerships would benefit from being strategic more than opportunistic. As a first step, the DML will need to assess potential partnerships in terms of the potential of the partnership to help the DML meet its goals, the likely incentives on both sides for the partnership, the maturity and stability of any technical standards required to make the partnership work, and the likely obstacles to consummating the partnership. A diversity of partnerships will be important. The advice and help of existing elements of community infrastructure could be valuable in this; for example, the IMU (in particular the CEIC) and its member societies, FIZ Karlsruhe, European Mathematical Society, the Association of Research Libraries (ARL) and similar organizations outside North America, existing mathematics digital libraries (such as *HathiTrust* and the EuDML), and prominent and influential mathematicians who have expressed an interest in the mission of the DML. Finally, while keeping in mind long-range goals and objectives, it is important to identify and pursue high-likelihood, high-potential-benefit, low-risk, near-term partnerships and agreements, even if somewhat limited in scope, as long as such partnerships can help illustrate the longer-term potential of partnering with DML. For example, a productive, beneficial partnership with arXiv might be achievable in relatively short order and at the same time be useful to illustrate some of the potential benefits of DML partnerships between the DML and content providers.

## Entity Collection

Even in advance of construction of a central repository, work could proceed immediately on development of adequate object classes for description and discovery of mathematical content in ways that complement existing capabilities—for example, at finer granularity—and on the aggregation of the lists of object instances for inclusion in the DML. The committee believes that the following mathematical objects and bibliographic entities are good targets for early DML development (each of which is discussed in more detail in Chapter 5):

- *Mathematical objects:* subject topics, sequences, functions, transforms, identities, symbols, formulas, and assorted mathematical media; and
- *Bibliographic entities:* people, homepages, journals, books, and bibliographies.

The committee recognizes that progress on aggregation, cleaning, and deduplication of these various lists will move at very different rates. Some of these lists may be completed quickly, while others that require input from many sources will mature slowly over time, and some might never be regarded as truly complete. Those in data rich areas may be ripe for initial developments. Still, the committee believes that the difficulty of completing some of these lists should not deter contributors from starting them or from converting what is already available into machine-readable formats, which can then feed various linking, navigation, search, and discovery services. Chapter 5 outlines which entity types should be targeted, at least initially, and gives some indication of the efforts required for each.

### Planning for More Complex Entities

Planning should start for the development of more complex lists where possible. These lists are outlined in Chapter 5, and some may be difficult to create and maintain. Wolfram|Alpha has a significant start on this with its continued fractions project. The potential rewards in terms of discovery and cross-linking are greatest if these mathematical objects can be adequately formalized and managed, even on a modest scale. These lists may benefit from starting small and growing slowly, to reduce the maintenance challenges before they become too burdensome, and by development of machine learning techniques for extraction of these entities from the literature.

The committee anticipates a fairly loose structure in cooperation with *Wikipedia*, with input from the Wolfram experience with continued fractions and others in managing problem lists.

### Data Structures

Initial effort is best invested in choosing an adequately flexible and extensible data structure, which needs to be easily expressible and exportable to handle diverse types of objects. The experience of Wolfram|Alpha, EuDML, and others working with metadata standardization will be essential input for this process. It is important to quickly codify the workflow for initiation of new lists of this kind and to gain a realistic assessment of the incremental cost of developing and maintaining new lists of various sizes and complexities. The intention is to lower the barriers to creation and maintenance of such lists to a point where there is substantial community enthusiasm for the activity. Simple user interfaces for the input of new entries and editing of existing ones consistent with schema restrictions are an essential requirement. The interface should be generic, much the same for all object classes, with customization as necessary for particular classes.[21]

---

[21] Prototype interfaces are provided by BibSonomy (http://www.bibsonomy.org/, accessed January 16, 2014), Zotero (http://www.zotero.org/, accessed January 16, 2014), and various library catalog tools.

## Growth and Cross-linking

Some initial effort will need to be expended on planning for eventual cross-linking of a substantial number of entries in different lists through semantic relations, such as connections between lists of authors and journals or mathematical symbols and equations. This initial step is not intended to build a complete ontology of mathematics, but obvious semantic links will need to be supported to the greatest extent possible. This would aid in the creation of a Web of mathematical information that supports further processing by modern methods of graphical data analysis and may yield unexpected visualizations and insights into the structure of the mathematical universe. The proposed development would likely benefit from starting small, demonstrating the successful ingestion of data and exposure of various facets incrementally, leveraging available ontologies and services, and building new ones as needed.

## Workflow Support

A workflow is a sequence of connected steps where each step concludes immediately before the next step begins. Workflow management systems in computer systems manage and define a series of tasks to produce a final outcome or outcomes. Once the task is complete, the workflow software ensures that the individuals responsible for the next task are notified and receive the data they need to execute their stage of the process. These systems can also automate redundant tasks and ensure that uncompleted tasks are followed up, as well as reflect the dependencies required for the completion of each task.

The DML could provide support for schema development and production software for editorial workflows involved in creation and maintenance of structured lists. It is to be expected that these workflows will evolve over time as different data sources and editorial agents become involved, and that somewhat different workflows may be required for different lists.

## RESOURCES NEEDED

As discussed throughout the report, the DML will require a small paid staff, technical infrastructure, a funded research portfolio to support relevant projects, and a governing board to ensure that the DML's components continue to function and develop properly. This section describes what is needed in each of these areas.

### Financial Resources

While it is difficult to accurately assess the necessary financial resources for the DML at this early stage, this section gives a general sense of the scale of the necessary human and technical resources. Some of these resources might be shared, too, depending on the particular arrangement developed for the DML. However, the amount of financial resources necessary is obviously an important component of evaluating the future development of the DML, and the committee provides the following recommendation for evaluating these resources before DML development.

**Recommendation: The initial DML planning group should set up a task force of suitable experts to produce a realistic plan, timeline, and prioritization of components, using this report as a high-level blueprint, to present to potential funding agencies (both public and private).**

The cost of development and upkeep for the DML will not be trivial but is currently too uncertain to be specified in this report. For some perspective on operating costs, arXiv may provide a reasonable example. In calendar year 2012, arXiv spent nearly $800,000 in expenses relating to the following:[22]

- Personnel costs (including benefits)—totaling $492,061
  - User support (2.70 full-time equivalent and 0.36 student)
  - Programming and system maintenance (2.13 full-time equivalent)
  - Management (0.50 full-time equivalent)
- Nonpersonnel costs—totaling $71,807
  - Servers (physical and virtual), hardware maintenance, storage and backup—$24,240
  - Network bandwidth and telephony—$10,867
  - Staff computers, software, and supplies—$2,700
  - Staff and arXiv Board travel—$34,000
- Indirect and in-kind costs—$208,631
  - College and department administration, staff support (26 percent of direct costs)—$146,606
  - Facilities (11 percent of direct costs)—$62,025
  - arXiv moderation (130+ moderators, varying time commitments)—volunteer efforts

---

[22] Cornell University, "Arxiv Projected Budget—Calendar Year 2012," August 29, 2012, https://confluence.cornell.edu/download/attachments/127116484/arXiv2012budget.pdf.

To give some perspective of potential costs of developing capabilities, the committee would like to draw attention to some of the resources recently devoted to developing and deploying the Wolfram|Alpha continuous fractions work[23] discussed in the Chapter 2 section "What Gaps Would the Digital Mathematics Library Fill?" Wolfram received a 1-year grant from the Alfred P. Sloan Foundation to prototype and build a technological infrastructure for collecting, tagging, storing, and searching a representative subset of mathematical knowledge (including definitions and theorems) and presenting it through a Wolfram|Alpha-like natural language interface. This work required some 3,000 hours of work from a team consisting of four professionals and one intern.

The subject of continued fractions was selected for this project because much of the relevant literature is older (therefore more representative of the type of content that can be utilized in a future system such as the DML) and is distinct from Wolfram's main computational expertise (as to lessen the bias in the results). The individuals who worked on it had no detailed prior knowledge about continued fractions, which made the work go slower than it would if it were performed by an expert in the field, but this example is likely representative of how the DML would be approached. However, three of the four team members have written multi-volume books about mathematics, as well as websites each having more than 10,000 pages, so they had some experience in covering a wide range of mathematics.

There was not enough time in a 1-year project to cover the 100,000 pages of printed continued fraction literature, so the team tried to explore and cover various content and presentation aspects to see what might be possible in future efforts. In most ways, this project succeeded in meeting its objectives but in some parts, especially fully computational representations of the content, the system still needs improvement.

In addition to having qualified people, two software infrastructure components were important in carrying out this project: Mathematica and Wolfram|Alpha. Mathematica allowed the team to check the mathematics and to generalize it, and Wolfram|Alpha allowed them to collect the information in such a way that one can access it through free-form language inputs and deliver the information in various formats, from Web to TeX.

This project is a meaningful example of how various DML features can be developed within a larger infrastructure. The following sections draw some specifics of needed human and technical resources to make the rest of the DML possible.

---

[23] M. Trott and E.W. Weisstein, "Computational Knowledge of Continued Fractions," *WolframAlpha Blog*, May 16, 2013, http://blog.wolframalpha.com/2013/05/16/computational-knowledge-of-continued-fractions/.

Lastly, the committee would like to note the importance of a sustained investment and commitment from its potential funders. The committee believes that a ramped investment pattern, starting as a prototype and scaling up, may be more beneficial than a large initial investment. The DML will require a long and sustained effort to be successful.

## Human Resources

A small paid staff will work to develop the DML vision, address issues that arise, pursue fruitful partnerships, and manage the day-to-day operations of the DML. The following is a list of staff functions that the committee sees as essential during the initial phases of the DML. These staffing needs will change as the DML grows and matures. The committee believes it is essential to include a distinguished mathematician in the senior management of the DML to provide credibility to the academic mathematics community and to gain startup funding and respectability in the nonprofit world.

- *Academic director.* A well-respected leader in both the technical and social aspects of the DML who is able to make editorial decisions and can engage and appoint editors and curators for their domain knowledge and reputation. This could be part-time position (e.g., half-time of a senior mathematician).
- *Executive director.* A manager with knowledge of large-scale data methods and digital libraries. This person would be responsible for directing the project manager, budget allocations, promotion of the project, and negotiations with partners, and also consulting with the academic director about priorities.
- *Project manager.* This person would be in charge of the creation of the DML. He/she would interface with programmers, contracting organizations, and technical partners.
- *System manager.* This person would be responsible for setting up adequate server infrastructure for day-to-day DML operations and for expanding operations as needed.
- *Data wrangler.* The person would work on an ongoing stream of specific data conversion projects and provide documentation of best practices. He/she would engage and oversee other volunteer, or possibly paid, staff and also set up and experiment with crowdsourcing tasks and implementations.
- *Rights management and legal.* This person would provide guidance on critical licensing and copyright choices for both data and software, and for possible negotiations of agreements with data and service providers. This may be a consultant position.

- *Research analyst.* This person would be responsible for keeping abreast of emerging technologies, researching solutions for identified problems, assisting the executive director with technology choices, and preparing white papers to explain proposals and processes.
- *Community liaison.* This person would be responsible for community building, advocacy of the project, intelligent responses to incoming emails, blog development, negotiations to engage and persuade partners to contribute data, and other such activities. This would likely be a full-time staff person or contractor.

### Technical Resources

A mathematics digital library requires a technical infrastructure. This infrastructure needs to support storage, backup, search, retrieval, and at least some support for analysis and visualization. Storage is needed for some documents, the software component, and the management data on the system. In general, a different storage solution will be needed for each type of data due to differences in usage associated with size, security issues, speed required, and level of backup needed. Security, in particular, will need to be carefully planned and assessed throughout the DML development to ensure that the data it stores will be well protected. Storage can be handled in-house by purchasing a number of servers or outsourced to server farms or cloud storage. The key is to plan for growth and to consider, for the operations of interest, whether it is more cost-effective to store data in multiple formats to facilitate search or to minimize storage and do data conversion on the fly. At this point it is not clear which option will be more economical.

Other resources required include machines for developing and testing software, backup facilities, and machines for monitoring and managing the system. For development, the key resource needs include high-end desktop computers, access to storage devices, Internet connectivity, backup facility, and a test bed environment for trying new features before they are launched. For managing and maintaining the ongoing system, handling the business tasks and associated financial issues, basic desktop machines with Internet connectivity, printers and associated fax, and backup to machines off the Internet for security are needed.

### Necessary Research Areas

There are many technological aspects of the proposed DML that are not currently possible. To help accelerate needed technological developments, the committee believes that several research areas can be targeted by the DML organization.

**Recommendation: The Digital Mathematics Library needs to build an ongoing relationship with the research communities spanning mathematics, computer science, information science, and related areas concerned with knowledge extraction and structuring in the context of mathematics and to help translation of developments in these areas from research to large-scale application.**

Some of these players include the following:

- National Institute of Standards and Technology,
- Cornell University Library (both Project Euclid and arXiv),
- American Mathematical Society (MathSciNet),
- Wolfram, and
- European technical partners in EuDML, including FIZ Karsruhe (zbMATH).[24]

These organizational partners would be the employers of some people engaged in DML work, funded by contracts approved by DML central administration and funded through some arrangement with DML funding sources. It may be best for the DML to collaborate with its partners to complete such work, rather than directly employing large numbers of its own people. All of the above partners have existing capabilities and services of this kind, which should not be threatened, but rather enhanced, by DML developments.

---

[24] FIZ Karsruhe is a nonprofit corporation and the largest nonuniversity institution for information infrastructure in Germany (http://www.fiz-karlsruhe.de). FIZ (together with EMS and Springer) is one of the joint owners/controllers of zbMATH (http://zbmath.org/).

# 5

# Technical Details

This chapter discusses some details of entity collections and technical considerations for the Digital Mathematics Library (DML). The lists discussed in this chapter are reasonable and obvious places for the DML to start developing its entity databases, but these may just be the starting point in an entity collection that is likely to evolve over time with the needs and capabilities of the DML. The ultimate goal of these lists is to provide interesting and nontrivial connections between topics, in particular the user features described in Chapter 2. The committee believes this is best accomplished by the DML organization overseeing the simpler entity collections first, which may have the most impact. These early lists can be managed in a straightforward, flexible, and sustainable way. Once this is achieved, the DML may benefit from moving on to more complex structures, such as ordered lists based on importance, relevance, pedagogical value, historical importance, etc., or to lists that can be (partially) ordered using different criteria and hyperlinked.

## ENTITY COLLECTION

This section discusses potential object types that the committee believes should be set up early in DML development, with details about location of relevant data sources and technical and political issues in data acquisition. These objects are divided into two categories: mathematical objects and bibliographic entities. Some of these entities are already data rich and can be developed by collaborating with existing databases and services.

Whenever possible, these areas of least resistance should be targeted first to establish a breadth of content within the DML.

## Mathematical Objects

The collection and organization of data on mathematical objects should be a high priority of any DML development. The following entities can be pursued and developed individually or jointly, but cross connections should be noted and exploited whenever possible.

### MathTopics

MathTopics would be a collection of mathematical subjects, topics, and terms that includes supporting definitions at various levels of formality and that indicates relations between topics derived from graphical analysis of book and journal data. This collection is practical to begin immediately, and some initial sources of information include MSC2010, *Wikipedia*, and tables of contents of mathematics books. As an application, MathTopics could be used to provide visualizations of the global structure of mathematical fields and their interactions.

Including information from open encyclopedic resources[1] and metadata records of entries in other encyclopedias behind paywalls[2] would be an extremely helpful service of the DML.

Encyclopedic information aggregation has been achieved before in limited cases, as in the case of the National Science Digital Library,[3] which indexed *MathWorld* and *PlanetMath* together. With the expansion of the linked open data approach, these cross-connections are happening in other domains as well. To use linked open data, interfaces need to allow cross-connections, and once an encyclopedia is available as linked open data, the data provider no longer has to be involved in the process of creating the combined resource. There are also some commercial entities providing metadata as linked open data,[4] so there is a sense that these connections may be possible in the near future.

The DML could provide dedicated search over this collection, with automated disambiguation of author names and superior subject navigation derived from graphical analysis of various forms of proximity between subjects. Google currently indexes this material but does not provide a means of browsing or navigating the material besides a simple search. Other methods of navigation, such as browsing and faceted search and browse, are very appealing and useful if available, but such systems typically do not have ontologies, and the data are not structured.

---

[1] Some encyclopedic resources include *Wikipedia*, *Springer Encyclopedia of Mathematics*, *StatProb*, *MathWorld*, *The Princeton Companion to Mathematics*, etc.

[2] *Encyclopedia of Statistical Sciences* would be of particular interest.

[3] National Science Digital Library, http://nsdl.org/, accessed January 16, 2014.

[4] See a list of linked open data encyclopedia data sets (e.g., http://datahub.io/dataset?q=encyclopedia) or search for encyclopedias in particular domains (e.g., http://datahub.io/dataset?q=biology+encyclopedia).

For topics that do not already have encyclopedia articles, the DML can flag that an article needs creation, and such an article can then be written in any of the available encyclopedia frameworks. Following the DML principle of not unnecessarily replicating data, and especially not unnecessarily replicating complex editorial structures, the DML would likely benefit from not providing its own encyclopedia publication infrastructure that would compete with and undermine the existing open encyclopedias.

## MathSequences

MathSequences would be a collection of mathematical sequences found throughout the literature. This list is already well developed in the Online Encyclopedia of Integer Sequences (OEIS). The DML could offer to help develop systematic hyperlinking of the text of all references, all author names to MathPeople, and all journal names to MathJournals, and systematic conversion of the data to standard machine-readable formats that can be understood by bibliographic and computational services. The OEIS data set, augmented with such enhancements, would be an example of what the DML should strive for in its data structures for other kinds of mathematical objects. Systematically reconstructing the OEIS as computable linked open data does not appear to be a very difficult task. Moreover, the solutions to difficulties encountered in this process should inform the choice of data schema for other similar collections. The main issue for DML involvement in the OEIS appears to be one of negotiating cooperation between the organizations.

## MathFunctions

MathFunctions would be a collection of mathematical functions found throughout the literature. This collection can begin immediately, and the National Institute of Standards and Technology Digital Library of Mathematical Functions (DLMF)[5] and the Wolfram Functions site[6] could provide the basis for a well-structured collection of mathematical special functions. This collection could then be added to MathTopics and used to tag compo-nents of articles and papers that discuss or apply specific special functions. This collection would likely take considerable time to populate extensively beyond the DLMF and Wolfram capabilities but could provide a wealth of information once reasonably established.

---

[5] National Institute for Standards and Technology, Digital Library of Mathematical Functions, Version 1.0.6, release date May 6, 2013, http://dlmf.nist.gov/.

[6] Wolfram Research, Inc., The Wolfram Functions Site, http://functions.wolfram.com/, accessed on January 16, 2014.

## MathTransforms

MathTransforms would be a searchable and browsable lookup table for classical transforms (e.g., Laplace, Fourier, Mellin) with links to computational resources. This could begin to be developed immediately in cooperation with DLMF and/or the Wolfram Functions site but would likely develop fully over a longer timeline. It is useful for mathematicians to be able to search or browse a table of transforms for various purposes: for inspiration, to see what is out there, and to see what might be adapted to a problem at hand. Moreover, such a table has the potential to be hyperlinked to the occurrences of its entries in the mathematical literature, which would be a step toward a more fine-grained indexing of the literature. Especially for rarely used functions and transforms, it is potentially rewarding for users to be able to find quickly where the same function or transform might have been used before. Special functions are often kept out of sight in higher mathematical constructions but have applications to other branches of mathematics. Making it easier for users to follow threads of their occurrences across the literature might easily lead to novel discoveries or unexpected relations between research in different branches of mathematics. Examples of such relations include the unexpected applications of Airy kernels and Painleve transcendents (functions) in random matrix theory, statistical physics, and elsewhere (Tracy and Widom, 2011; Forrester and Witte, 2012).

## MathIdentities

MathIdentities would be an organized table of classical combinatorial identities and methods of reduction and proof of such identities. There has been huge progress in recent years in computer methods for proving classical combinatorial identities, including closed-form summation formulas. This means that a great many simplifications of algebraic sums and proofs of algebraic identities can be done rapidly and with high reliability by machine.[7] For the same reasons as identified above for tables of functions and transforms, it may be instructive for mathematicians to browse through tables of identities and to follow links to applications of identities in the literature. This collection could begin immediately and progress similar to MathFunctions and MathTransforms.

---

[7] See Gould (1972); Wolfram|Alpha (http://www.wolframalpha.com/); the work of Christian Krattenthaler (http://www.mat.univie.ac.at/~kratt/, accessed January 16, 2014), including Mathematica packages HYP and HYPQ for the manipulation and identification of binomial and hypergeometric series and identities (C. Krattenthaler, "HYP and HYPQ," http://www.mat.univie.ac.at/~kratt/hyp_hypq/hyp.html, accessed January 16, 2014); and Gauthier (1999).

## MathSymbols

MathSymbols would be a collection of mathematical symbols with commonly accepted special meanings, to be cross-linked as well as possible with MathTopics, and if possible with place of first usage. Within restricted domains, symbols often acquire stable conventional meanings, and sometimes this is true across all of mathematics. Some work has been done on developing a consensus of mathematical notations across cultures (Libbrecht, 2010), and this Notation Census[8] is a meaningful precursor to what the committee envisions. The collection that the committee envisions for the DML is complex and may require multiple steps. As a first step before embarking on a complete index, the DML could partner with resources such as MathSciNet and zbMATH to create a collection of journal article titles that contain any mathematical symbols. This would provide a core set of symbols with authoritative links to the literature. The meaning of those symbols could be established by a small community-sourcing exercise. The symbols could be linked to MathTopics at the collection level, and then MathNavigator tool could serve these links to MathTopics entries from a reference to any article anywhere in the mathematical literature that has the same symbol in its title. This might be considered a preliminary exploration before attempting to do a similar but more ambitious undertaking for formulas or equations.

## MathFormulas

MathFormulas would be a collection of mathematical formulas and their variations, initially those of special interest and importance. This collection could assist in supervised machine learning processes for the creation and maintenance of a larger body of formulas and equations. This is a long-term collection goal and DLMF, Wolfram, and Springer would be desired partners, especially the data in Springer's LaTeX Search. This is an ambitious list to attempt to collect, because there are serious challenges to overcome because of superficial variations in the way every given expression might by written (as discussed in Chapter 3). Still, initial progress has been made by several teams of researchers, and the DML could provide a nexus for further research, a forum for tracking advances in this field, and eventually some attempt to create and maintain an authoritative list of at least those formulas considered interesting or important enough to be recognized and assigned an HTTP URL. Further open efforts at both supervised learning relative to these exemplars and unsupervised learning similar to the Springer effort, with linking to the literature, should also be attempted, motivated by applications to formula search as indicated in Chapter 3.

---

[8] "Notation Census Manifest," last edited March 9, 2013, http://wiki.math-bridge.org/display/ntns/Notation-Census-Manifest.

**MathMedia**

MathMedia would be a collection of images, photos, videos, and presentations—or links to such—relating to mathematics. Video clips from conferences and presentations, visualizations of results and simulations, and portraits of mathematicians who contributed to the research field could all be included in the DML and systematically integrated with the mathematics literature. Widespread collection of media entities could begin immediately and would likely continue to evolve. Many mathematics conferences are already filming and posting speakers' presentations, and it would be opportunistic for the DML to arrange for these data to be indexed and sorted based on known information such as the title of the presentation, author(s) and presenter(s), date, name of conference and/or section, etc. Other information on the contents of the presentation, which may be more difficult to automatically categorize, can be tagged by community sourcing. In terms of mathematician portraits, there are several images of mathematicians that may be of interest, such as Oberwolfach Photo Collection,[9] Portraits of Statisticians,[10] Microsoft Academic Search Profiles,[11] Halmos (Beery and Mead, 2012), and Pólya (Alexanderson, 1987).

### Bibliographic Entities

The following bibliographic data collection entities are a needed element of the DML, and their collection can begin quickly—and largely be completed—since much of the information is already available elsewhere through existing information resources. These entities can be viewed as part of the necessary infrastructure of the DML and are key areas for developing partnerships (as discussed in Chapter 3). However, the collection and development of these entities are not meaningful on their own and should only be pursued as part of a larger DML development.

---

[9] Oberwolfach Photo Collection, http://owpdb.mfo.de/, accessed January 16, 2014.

[10] Portraits and Articles from Biographical Dictionaries, revised July 10, 2013, http://www.york.ac.uk/depts/maths/histstat/people/.

[11] Microsoft Academic Search, "Overview," http://academic.research.microsoft.com/About/Help.htm, accessed January 16, 2014.

## MathPeople

MathPeople would be a lean authority file for mathematical people with links to and selected data from homepages, *Wikipedia*, MacTutor, Math Genealogy, zbMath Open Author Profiles, Celebratio.org, MacTutor, MathSciNet, and so on. There was an effort by the International Mathematical Union in 2005 to build a Federated World Directory of Mathematicians,[12] but it was abandoned due to copyright and privacy concerns and inadequate federated search technology. More recently, zbMATH Author Profiles and data in Microsoft Academic Search's Top Authors in Mathematics offer machine access to approximate authority records for about half a million authors in mathematics and related fields, with no apparent legal restriction on further processing of the data. It would be a straightforward application of machine processing and community input to deduplicate these lists, sync them also with the Virtual International Authority Files of all mathematicians, and thereby obtain a combined DML index of all mathematicians, both living and deceased, who have ever published a book or article in mathematics. This data set would include additional information about the mathematicians' fields, their collaborators, and their numbers of publications. This would then provide a graphical data set with about half a million nodes for authors and editors, and some fraction of that number of nodes for books they wrote and edited, and a few thousand subject nodes. This could be used very quickly as a test bed for application of modern graphical visualization methods to provide subject navigation, and otherwise as a major framework for organizing other facets of DML information.

MathSciNet has high-quality representation of the collaboration graph for mathematical articles, obtained through many years of manual curation of book and article metadata records, and MathSciNet offers a glimpse into this proprietary data set with its computation of minimal distance paths through the collaboration graph from one author to another. These collaborator connections are helpful and allow users to see if an author's collaborators are working in relevant areas, but they do not provide links to other relevant data. Having access to similar information in addition to the other data that the DML is proposing to collect (such as theorems, research areas, homepages), this information then becomes significantly more integrated into a larger picture of the mathematics research community.

With suitable graphical visualization, MathPeople could provide users with a sense of the "geography" of mathematics, how the subfields of mathematics are related to each other through the collaborations of authors, and how this structure has evolved over time.

---

[12] International Mathematical Union, "Personal Homepages and the World Dictionary of Mathematicians," http://www.mathunion.org/MPH-EWDM/ last updated December 13, 2012.

**MathHomepages**

MathHomepages would be a table linked to MathPeople, but with indications of depth of content (e.g., curriculum vitae, photo, bibliography). From a user perspective, this may appear to be a simple variation of MathPeople; however, a person can have more than one homepage, each of which may contain references and connections to subjects and collaborators. It would be beneficial for there to be separate tables for homepages and for people and for these to be cross-linked by a general, extensible data architecture, such that the cross-links are easily maintainable and correctible. This is not trivial, and it is illustrative of the maintenance problem for Web-based data. Much of this data could be mined from sources such as the Microsoft Academic Search API, some subject specific collections in the Web (e.g., for number theory, probability), and easily completed and maintained by Web-spidering, community input, and self-registration. While people stay the same, their homepages and affiliations may change. The relation between people and their homepages could be treated as a simple case of a dynamically changing data set, and methods and interfaces could be developed to respect this aspect.

This information would be useful to mathematics researchers because it can help find people with common names and can be useful to the larger DML because it helps with interlinking other data.

**MathJournals**

MathJournals would be a deduplicated and cleaned list of serials in mathematics, past and present, with indications of online availability and subject coverage. Most of this data currently exists and is maintained openly by a number of agents (such as Ulf Rehmann, MathSciNet, zbMATH, the Online Computer Library Center).[13] There are several lists of math journals in various places, many of them accessible and reusable, but none of these lists provides easy access to the features that researchers would like, including the following:

- Links to journal homepages whenever they exist;
- Information about the number of articles published and subject areas covered;
- Copyright and rights information for authors; and
- Simple search over the list.

---

[13] See also UlrichsWeb.com for a proprietary solution across all fields.

The entire math journals list is only a few thousand entries, but the number of readily available attributes of a journal is potentially large and, in principle, unlimited. Some desired capabilities for the journals list that will require some initial work and maintenance are the following:

- Graphical displays (e.g., nodes with size proportional to various journal metrics and locations reflective of their subject coverage, linking to MathTopics below) that could easily be derived;
- Display of journals by defined metrics (e.g., in cooperation with eigenfactor.org[14]), which uses recently developed methods of network analysis and information theory to evaluate the influence of scholarly journals and for mapping the structure of academic research; and
- Access to identities of all authors who ever published in a journal with links to MathPeople.

These are typical functionalities that the standard abstracting and indexing services could provide but currently do not offer. Aggregating and displaying this information would give users a quick overview of the whole field of mathematics from the point of view of its journal coverage, and graphical relations derived from such information could feed into tools for navigation of mathematical information. While no such navigation tools are currently available, they could easily be built over a MathJournals list, especially if cross-linked to MathPeople (e.g., "authors who published in this journal also published in these other journals").

**MathBooks**

MathBooks would be a list of mathematical books at all levels, from elementary to advanced, with links to and selected data from publishers. Some of these data already exist through services such as MathSciNet, zbMATH, OCLC, OpenLibrary, and Ulf Rehmann, but this bibliographic entity is less developed than the previous three discussed in this section. A plethora of openly accessible metadata about books in all fields has been released in the past few years by academic libraries and library cooperatives.[15] Considering just books in mathematics and related fields, the information in these data releases swamps what is currently available in MathSciNet and zbMATH both in quantity of titles and depth of information about each title.

---

[14] Eigenfactor, http://www.eigenfactor.org/, accessed January 16, 2014.
[15] Most notable are the British Library release of millions of catalog records in 2010 (British Library, 2010) and the OCLC recommendation to use Open Data Commons Attribution License (ODC-BY) for WorldCat data in August 2012 (OCLC, 2012).

A large number of elementary mathematics books in these releases are not indexed at all by MathSciNet and zbMATH, but they may be of value to students and teachers of mathematics. There is potential to index this collection in ways that would provide novel recommendation and discovery services over mathematics book data for students and teachers as well as researchers and those who apply mathematics in other fields. The process of indexing and cleaning these data, and providing enhanced discovery services over them, should be a fairly routine application of machine learning methods, which could be done as a standalone project and which would provide a first test of DML machine learning capabilities. The general methods involved would not be domain-specific, and they could be applied also to other non-math domain-specific collections. However, mathematics is special in that is already has a well-developed subject ontology for the field, the MSC2010. Cross-linking of the library books data with subject tags from either MathSciNet or zbMATH, and with author identities from MathPeople and the Virtual International Authority File,[16] should aid readers in navigating the universe of mathematical concepts by reference to the statistics of its book data. The DML could also use these data to suggest key textbooks and research texts for specific subjects or theorems.

## MathBibliographies

MathBibliographies would be a collection of bibliographies of various topics in mathematics, including personal and subject bibliographies. Initial sources for these data include Celebratio Mathematica,[17] IMS Scientific Legacy,[18] other subject bibliographies, and bibliographies from books contributed by participating publishers. This collection could be cross-linked to MathPeople and MathTopics. The structure of aggregated collections of such bibliographies could then inform search and navigation services, just as reference lists of articles do already. The key functionality for users is to make it easy for them to select, annotate, and tag bibliographic items. MathSciNet's MRLookup tool already provides a useful open interface for acquisition of modest-sized bibliographies from data in MathSciNet. Similar data are readily available from Microsoft Academic or Google Scholar, but there is not yet any tool comparable to MRLookup for acquiring data from those sources, and neither is there any good tool for aggregation and deduplication of data from multiple sources, as would typically be neces-

---

sary to develop the bibliography of any topic where mathematics reaches into other domains.

## MathArticles

MathArticles would be a collection of metadata of journal articles in various topics in mathematics. Some initial sources for these data include MathSciNet, arXiv, Web of Science, Google Scholar, and Microsoft Academic, among others. There would be considerable connections between the other bibliographic entities proposed in this section.

## TECHNICAL CONSIDERATIONS

This section lists a number of technical considerations that the committee believes will influence the development of the DML and its information management structures. Some key issues discussed are managing diverse data formats, incorporating math-aware tools and services, appropriately dealing with authority control, and managing client-side versus server-side software. None of these discussions are intended to be overly prescriptive, but to raise issues that the committee feels are very important.

### Data Formats

For annotation and sharing of data it is necessary to have a format that fulfills certain requirements as follows:

1. Easy to use and ideally human readable;
2. Can be implemented into any recording, analysis, or management tool;
3. Open and freely available;
4. Inherently extensible and flexible for science continually changes; and
5. More or less unrestricted—that is, it should not restrict the user or strictly require entries.

At some points, format conventions have to be introduced. This is the process of schema modeling and introduction, which is by now fairly well understood. It is essential to clearly separate format from content. Documentation about formats can be maintained along with the data model, and a place to record and maintain property definitions can be included. For any given list of objects, the expected internal structure of those objects and their expected relations with other objects define an ontology. There are many tools available for creating and maintaining ontologies (as discussed in Chapter 1).

Essentially the same metadata structure can be used for metadata of all kinds of objects, including documents, people, organizations, subjects, or mathematical concepts. The schema for the object is type dependent, with some sub-typing within major types like documents.[19] To the greatest extent possible, existing or adapted schemas can be used. But for mathematical concepts in particular, development of adequate schemas will be a slower process, informed by the success of partners such as Wolfram and OEIS with experience in handling such data and the experience of numerous others in development of math-aware tools and services.

## Math-Aware Tools and Services

There currently exist math-aware tools and services that can competently manage mathematical syntax and formatting. Such tools and services are essential for tasks such as conversion between formats that are different mathematically and semantic parsing of mathematical documents. However, many current resources do not functionally handle mathematical notation and syntax, and this limits how the mathematical community can use these resources. Significant interest in better utilizing math-aware tools and services is apparent in the series of Conferences on Intelligent Computer Mathematics.[20] The following is from the announcement of their digital mathematics library conference track, chaired by Petr Sojka:

> Track objective is to provide a forum for development of math-aware technologies, standards, algorithms and formats towards fulfillment of the dream of global digital mathematical library (DML[21]). Computer scientists (D) and librarians of digital age (L) are especially welcome to join mathematicians (M) and discuss many aspects of DML preparation. Track topics are all topics of mathematical knowledge management and digital libraries applicable in the context of DML building—processing of math knowledge expressed in scientific papers in natural languages, namely:
>
> - Math-aware text mining (math mining) and MSC classification;
> - Math-aware representations of mathematical knowledge;

---

[19] The basic framework for most document types is already provided by the BibTeX ontology, and easily implementable in JSON as BibJSON or in XML according to some variant of the NLM DTD (http://dtd.nlm.nih.gov/), which is currently used by the EuDML for document records.

[20] Conferences on Intelligent Computer Mathematics, last modified July 10, 2013, http://www.cicm-conference.org/2013/cicm.php.

[21] Please note that the DML discussed in this quotation is distinct from the DML vision laid in this report.

- Math-aware computational linguistics and corpora;
- Math-aware tools for [meta]data and full-text processing;
- Math-aware OCR and document analysis;
- Math-aware information retrieval;
- Math-aware indexing and search;
- Authoring languages and tools;
- MathML, OpenMath, TeX and other mathematical content standards;
- Web interfaces for DML content;
- Mathematics on the Web, math crawling and indexing;
- Math-aware document processing workflows;
- Archives of written mathematics;
- DML management, business models;
- DML rights handling, funding, sustainability; and
- DML content acquisition, validation and curation.

DML track is an opportunity to share experience and best practices between projects in many areas (MKM, NLP, OCR, IR, DL, pattern recognition, etc.) that could change the paradigm for searching, accessing, and interacting with the mathematical corpus.[22]

Integrating math-aware tools and services developed by diverse partners may be challenging but would benefit the DML. One math-aware standard of particular interest to DML developments proposed in this report is that of OpenMath,[23] which is an extensible standard for representing the semantics of mathematical objects. The OpenMath website describes its objective as follows:

OpenMath is an emerging standard for representing mathematical objects with their semantics, allowing them to be exchanged between computer programs, stored in databases, or published on the worldwide Web. While the original designers were mainly developers of Mathematical objects encoded in OpenMath can be

- Displayed in a browser,
- Exchanged between software systems,
- Cut and pasted for use in different contexts,
- Verified as being mathematically sound (or not!), and
- Used to make interactive documents really interactive.

OpenMath is highly relevant for persons working with mathematics on computers, for those working with large documents (e.g., databases, manuals) containing mathematical expressions, and for technical and mathemati-cal publishing. The worldwide OpenMath activities are coordinated within the OpenMath Society, based in Helsinki, Finland. It is coordinated by an executive committee, elected by its members. It organizes regular workshops and hosts a number of electronic discussion lists. The Society brings together tool builders, software suppliers, publishers and authors.

---

[22] Conferences on Intelligent Computer Mathematics, "Track B: DML," last modified March 4, 2013, http://www.cicm-conference.org/2013/cicm.php?event=dml.

[23] OpenMath, http://www.openmath.org/, accessed January 16, 2014.

This standard and the community that has developed around it should contribute to development of the DML.

### Authority Control

The committee expects continuing advances in authority control[24] over entities and the provision of adequate human-computer interfaces for the semi-automated curation of large digital collections. Some customization of these tools will be necessary to apply them to mathematical objects. However, once the tools are built and the editorial workflows established, these tools and workflows should be largely replicable across multiple distributed nodes in the network of bibliographic data stores contributing to the DML.

The problems of identification and deduplication of conventional bibliographic records are by now largely solved. Solutions and workflows developed by other organizations, such as OCLC and ORCID, should be adopted to the extent that these organizations are willing to share them. In mathematics, existing automated tools such as MRef and MRLookup[25] return similar matches to queries from traditional bibliographic reference data. This enables machine enhancement of reference lists by matching into the MathSciNet database. However, the universe of mathematics information resources of interest to the DML is not limited to traditionally published items alone. Neither ORCID nor MRef are comprehensive in providing identifiers for all mathematicians.

The problem of identification and deduplication of various of mathematical entities remains a research problem on which more effort will need to be expended before the fullest potential of DML navigation can be realized. Like searching for articles, exploring the citation graph in the DML will need to deal with the "identity problem"—that is, the problem of deciding that two citations are actually to the same article, although the names of authors can be different (e.g., initial instead of full first name) the journal names can be altered (e.g., abbreviations or misspellings), the ordering of terms can be changed, and so on. Another aspect of this problem is determining to what degree lightweight authorities (e.g., MathPeople, men-tioned above) can be relied on as supplements to more traditional authorities. It is interesting to note that Google Scholar and Microsoft Academic Search deal with this problem reasonably well by using a statistical modeling approach rather than the more in-depth approach of writing down all possible transformations and then unraveling those transformations.

---

[24] In library science, authority control is a process that organizes bibliographic information by using a single, distinct name for each topic.

[25] American Mathematical Society, MRLookup, http://www.ams.org/mrlookup, accessed January 16, 2014.

### Client-Side Software

The DML would likely benefit from using a combination of client-side software and Web services to provide its content to users. Client-side software can be thought of as a computer application, such as a Web browser, that runs on a user's local workstation and connects to a server as necessary. If part of the DML were run client-side, a user would download a DML application that would carry out much of the data processing on the users machine, thereby lessening the server load on the DML. However, it is not always clear what resources are available on the user's machine, and users may not like the DML application using their machine's potentially limited storage and processing ability. Another concern is DML security; if too much of the DML data and processing is pushed client side, it may become an easy target for unintended manipulation. To balance the security and processing load concerns, the DML may benefit from moving much of the processing layer client-side while keeping the data layer server-side (or accessible only as a Web service that cannot be easily manipulated).

There are a number of services that use a mix of client-side software and Web services to provide enhanced document navigation capabilities, some of which may serve as an example of how to set up the DML:

- BibSonomy[26] (very open data and services, great scrapers for acquiring bibliographic metadata from publisher sites),
- CiteULike[27] (could easily go the way of Mendeley, which had a partnership with Springer at one time but has since stalled),
- Connotea,[28]
- Delicious,[29]
- JabRef[30] (desktop bibliography manager, syncs with BibSonomy),
- Mendeley,[31]

---

[26] BibSonomy, http://www.bibsonomy.org/, accessed January 16, 2014.

[27] CiteULike, http://www.citeulike.org/, accessed January 16, 2014.

[28] Nature Publishing Group, Connotea, http://www.connotea.org/, discontinued service on March 12, 2013.

[29] Delicious, https://delicious.com/, accessed January 16, 2014.

[30] JabRef, last updated October 29, 2013, http://jabref.sourceforge.net/.

[31] Mendeley, http://www.mendeley.com/, accessed January 16, 2014.

- MindMaps,[32]
- Papers,[33]
- Scholarometer,[34] and
- Zotero[35] (open source, but focused on the humanities).

---

[32] For general concept, see "Mind Map," *Wikipedia*, last modified January 14, 2014, http://en.wikipedia.org/wiki/Mind_map. For a specific software implementation, see Docear—The Academic Literature Suite (http://www.docear.org/). Docear (pronounced "dog-ear") is a free and open academic literature suite that integrates tools to search, organize, and create academic literature into a single application. Docear works seamlessly with many existing tools like Mendeley, Microsoft Word, and Foxit Reader.

[33] Mekentosj B.V., Papers 3, http://www.papersapp.com, accessed January 16, 2014.

[34] Indiana University, Scholarometer, http://scholarometer.indiana.edu/, accessed January 16, 2014.

[35] Zotero, http://www.zotero.org/, accessed January 16, 2014.

# Appendixes

# A

# Meeting Agendas and Other Inputs to the Study

**MEETING 1**
**NOVEMBER 27-28, 2012**
**WASHINGTON, D.C.**

| | |
|---|---|
| Discussion of Study Goals with Sponsor | *Daniel Goroff*, Alfred P. Sloan Foundation |
| Update on European Digital Mathematics Library Project | *Thierry Bouche*, Scientific Coordinator, EuDML Project, and Cellule MathDoc and Institut Fourier, Université de Grenoble (via WebEx) |
| Discussion with Representatives of Mathematical Information Resources | *François G. Dorais*, Moderator, MathOverflow, and John Wesley Young Research Instructor, Department of Mathematics, Dartmouth College |
| | *Michael Trott*, Senior Researcher, Wolfram |
| | *Paul Ginsparg*, Founder, arXiv.org, and Professor of Physics, Cornell University (via WebEx) |
| Discussion of Study Goals with Major Professional Societies | *Donald McClure*, Executive Director, American Mathematical Society |
| | *James Crowley*, Executive Director, Society for Industrial and Applied Mathematics |

| Discussion with Representatives of Digital Libraries Outside of Mathematics | *David Lipman*, Director of the National Center for Biotechnology Information<br>*Wayne Graves*, Director of the Office of Information Systems at the Association for Computing Machinery |
|---|---|

### MEETING 2
### FEBRUARY 19-20, 2013
### WASHINGTON, D.C.

| Current State of Semantic Libraries and Active Documents | *Michael Kohlhase*, Professor of Computer Science, Jacobs University in Bremen, Germany |
|---|---|
| Searching Outside of Mathematics | *Herb Roitblat*, Chief Scientist and Chief Technology Officer, OrcaTec |
| Building Infrastructure for Digital Libraries | *Andrew McCallum*, Professor of Computer Science, University of Massachusetts, Amherst |
| Crowdsourcing in Chemistry | *Antony Williams*, Vice President of Strategic Development, Royal Society of Chemistry |

### MEETING 3
### MAY 6, 2013
### MINNEAPOLIS, MINNESOTA

| Input on Desired Capabilities from the Mathematics Community | *Participants included:*<br>Kris Fowler<br>Andrew Odlyzko<br>George Sell |
|---|---|

### MEETING 4
### MAY 30-31, 2013
### EVANSTON, ILLINOIS

| Discussion with Representatives of Mathematical Information Resources | *Michael Trott*, Senior Researcher, Wolfram |
|---|---|

Scale and Cost of
Running a Large
Digital Library

*John Wilkin*, Associate University Librarian,
    Library Information Technology, University
    of Michigan

Input on Desired
Capabilities from
the Mathematics
Community

*Participants included:*
Patrick Allen
Dean Baskin
Anna Marie Bohmann
Yanxia Deng
Honghaw Gai
Elton Hsu
Ben Knudsen
Chao Liang
Clark Robinson
Melissa Tacy

# B

# Biographical Sketches of Committee Members and Staff

INGRID DAUBECHIES, *Co-Chair,* is a professor of mathematics at Duke University. She completed her undergraduate studies in physics at the Vrije Universiteit Brussel in 1975. She obtained her Ph.D. in theoretical physics in 1980 and continued her research career at that institution until 1987, rising through the ranks to positions roughly equivalent with research assistant professor in 1981 and research associate professor in 1985. Dr. Daubechies then moved to the United States, taking a position at the AT&T Bell Laboratories' facility in Murray Hill, New Jersey. Earlier that same year, she had made her best-known discovery: the construction of compactly supported continuous wavelets. From 1993 to 2011, Dr. Daubechies was a professor at Princeton University, where she was especially active within the Program in Applied and Computational Mathematics. She was the first female full professor of mathematics at Princeton. In January 2011 she moved to Duke University to serve as a professor of mathematics. She is a member of the National Academy of Sciences.

CLIFFORD LYNCH, *Co-Chair,* is executive director of the Coalition for Networked Information (CNI), which he has led since 1997. CNI, jointly sponsored by the Association of Research Libraries and EDUCAUSE, includes about 200 member organizations concerned with the intelligent uses of information technology and networked information to enhance scholarship and intellectual life. CNI's wide-ranging agenda includes work in digital preservation, data intensive scholarship, teaching, learning and technology, and infrastructure and standards development. Prior to joining CNI, Dr. Lynch spent 18 years at the University of California Office of the President, the last 10 as director of library automation. Dr. Lynch, who holds a Ph.D. in computer science from the University of California, Berkeley, is an adjunct professor at Berkeley's School of Information. In 2011, he was appointed co-chair of the National Research Council's (NRC's) Board on Research Data and Information. He serves on numerous advisory boards and visiting committees. His work has been recognized by the American Library Association's Lippincott Award, the EDUCAUSE Leadership Award in Public Policy and Practice, and the American Society for Engineering Education's Homer Bernhardt Award.

KATHLEEN M. CARLEY is a professor in the School of Computer Science at Carnegie Mellon University. She is the director of the Center for Computational Analysis of Social and Organizational Systems (CASOS), a university-wide interdisciplinary center that brings together network analysis, computer science, and organization science and has an associated National Science Foundation (NSF)-funded training program for Ph.D. students. Dr. Carley's research combines cognitive science, social networks, and computer science to address complex social and organizational problems. Her specific research areas are dynamic network analysis, computational social and organization theory, adaptation and evolution, text mining, and the impact of telecommunication technologies and policy on communication, information diffusion, and disease contagion and response within and among groups, particularly in disaster or crisis situations. Dr. Carley and her team have developed infrastructure tools for analyzing large-scale dynamic networks and various multi-agent simulation systems. The infrastructure tools include the ORA, a statistical toolkit for analyzing and visualizing multi-dimensional networks. Another tool is AutoMap, a text-mining system for extracting semantic networks from texts and then cross-classifying them using an organizational ontology into the underlying social, knowledge, resource, and task networks. Dr. Carley is the founding co-editor with Al Wallace of *Computational Organization Theory* and has co-edited several books in the computational organizations and dynamic network area.

TIMOTHY W. COLE is a professor of library and information science and the head of the Mathematics Library at the University of Illinois, Urbana-Champaign. He received a B.S. in aeronautical and astronautical engineering and a M.S. in library and information science from the University of Illinois, Urbana-Champaign. His research interests include metadata best practices, digital library system design, digital library interoperability protocols, and the use of XML for encoding metadata and digitized scholarly resources in science, mathematics and literature.

JUDITH L. KLAVANS is the principal investigator on the Mellon-funded Computational Linguistics for Metadata Building (CLiMB) research project, now based at the College of Information Studies at the University of Maryland (UMD). In addition to leading the project, Dr. Klavans is involved in developing analysis and filtering techniques for the extraction of metadata, particularly through thesaurus-driven disambiguation. She is also involved in the Defense Advanced Research Projects Agency (DARPA)-funded TIDES multilingual multimedia summarization project in which her primary technical role is the areas of utility evaluation and in coherence for summarization. Dr. Klavans is currently a research professor at the College of Information Studies at UMD. Dr. Klavans holds a Ph.D. in linguistics from the University of London and has worked on numerous computer science, digital library, and digital government projects. In particular, she has served as principal investigator on several other large research projects, including the NSF-funded PERSIVAL medical digital library, the NSF and Bureau of Labor Statistics-supported Digital Government Research Center joint project with University of Southern California-ISI, and DARPA-funded TIDES multilingual summarization project. Her research interests include linguistics, digital library research, language, and natural language systems. Dr. Klavans initiated the CLiMB project at Columbia University in 2002.

YANN LeCUN is a professor of computer science at the Courant Institute of Mathematical Sciences at New York University (NYU) since 2003 and was named Silver Professor in 2008. Dr. LeCun received a Ph.D. in computer science from the Université Pierre et Marie Curie, Paris, in 1987. He joined the Adaptive Systems Research Department at AT&T Bell Laboratories in Holmdel, New Jersey, in 1988, where he later became head of the Image Processing Research Department, part of Larry Rabiner's Speech and Image Processing Research Laboratory at AT&T Labs-Research in Red Bank, New Jersey. In 2002, he became a fellow of the NEC Research Institute (now NEC Labs America) in Princeton, New Jersey. He then began his tenure at NYU, where he remains. Dr. LeCun's research focuses on machine learning, computer vision, pattern recognition, neural networks, handwriting recognition, image compression, document understanding, image processing, VLSI design, and information theory. His handwriting recognition technology is used by several banks around the world, and his image compression technology is used by hundreds of websites and publishers and millions of users to access scanned documents on the Web.

MICHAEL LESK is a professor of library and information science at Rutgers—The State University of New Jersey and past department chair (2005-2008). After receiving a Ph.D. in chemical physics, Dr. Lesk joined the computer science research group at Bell Laboratories and from 1984 to 1995 managed computer science research at Bellcore. He was then head of the division of information and intelligent systems at NSF (1998- 2002), and then joined Rutgers. He is best known for work in electronic libraries, and his book *Practical Digital Libraries* was published in 1997 by Morgan Kaufmann and the revision *Understanding Digital Libraries* appeared in 2004. His research has included the CORE project for chemical information, and he wrote some Unix system utilities including those for table printing (tbl), lexical analyzers (lex), and intersystem mail (uucp). His other technical interests include document production and retrieval software, computer networks, computer languages, and human-computer interfaces. He is a fellow of the Association for Computing Machinery, received the Flame award from the Usenix association, and in 2005 was elected to the National Academy of Engineering.

PETER J. OLVER is the head of, and a professor in, the School of Mathematics at the University of Minnesota. Before joining the University of Minnesota, he was a Dickson Instructor at the University of Chicago and a postdoc at the University of Oxford. He is currently the chair of two committees with the International Mathematical Union: the Committee on Electronic Information and Communication and the Moderating Group of the Blog on Mathematical Journals. Dr. Olver is also a member of the American Mathematical Society and the Society for Industrial and Applied Mathematics. His research interests revolve around the applications of symmetry and Lie groups to differential equations. He is the author of four books and 130 papers published in refereed journals that include applications in computer vision, fluid mechanics, elasticity, quantum mechanics, Hamiltonian systems, the calculus of variations, geometric numerical methods, differential geometry, algebra, and classical invariant theory. He received a bachelor's degree in applied mathematics from Brown University and a Ph.D. in mathematics from Harvard University.

JIM PITMAN is a professor in the departments of statistics and mathematics at the University of California, Berkeley. Before joining the faculty at UC Berkeley, Dr. Pitman held a position in the Department of Mathematics and Mathematical Statistics at the University of Cambridge, England. Dr. Pitman has devoted much effort to promote the development of open access resources in the fields of probability and statistics. As a member of the Executive Committee of the Institute of Mathematical Statistics (IMS) from 2005 to 2008, he guided the IMS through implementation of a policy to promote open access to all of its professional journals, through systematic deposit of peer-reviewed final versions of all articles on arXiv.org and to provide technical support to other organizations willing to do the same. He has a continuing interest in the technical management of scientific information in ways that encourage individuals and small organizations to maintain high-quality knowledge repositories that are openly accessible. Dr. Pitman holds a B.Sc. in statistics from the Australian National University, Canberra, and a Ph.D. in probability and statistics from Sheffield University.

ZHIHONG (JEFF) XIA is an Arthur and Gladys Pancoe Professor of Mathematics at Northwestern University. He joined Northwestern in 1994 after serving as an associate professor at both the Georgia Institute of Technology and Harvard University. His research interests include dynamical systems, Hamiltonian dynamics, celestial mechanics, and ergodic theory. Dr. Xia received a B.S. in astronomy from Nanjing University in China and a Ph.D. in mathematics from Northwestern University.

## *Staff*

MICHELLE SCHWALBE is a program officer with the Board on Mathematical Sciences and Their Applications (BMSA) within the NRC. She has been with the National Academies since 2010, when she participated in the Christine Mirzayan Science and Technology Policy Graduate Fellowship Program with BMSA. She then joined the Report Review Committee of the National Academies before re-joining BMSA. With BMSA, she has worked on assignments relating to verification, validation, and uncertainty quantification; the future of mathematical science libraries; the mathematical sciences in 2025; and the Committee on Applied and Theoretical Statistics. Her interests lie broadly in mathematics, statistics, and their many applications. She received a B.S. in applied mathematics specializing in computing from the University of California, Los Angeles, an M.S. in applied mathematics from Northwestern University, and a Ph.D. in mechanical engineering from Northwestern University.

SCOTT T. WEIDMAN is the director of the NRC's BMSA. He joined the NRC in 1989 with the Board on Mathematical Sciences and moved to the Board on Chemical Sciences and Technology in 1992. In 1996 he established a new board to conduct annual peer reviews of the Army Research Laboratory, which conducts a broad array of science, engineering, and human factors research and analysis, and he later directed a similar board that reviews the National Institute of Standards and Technology. Dr. Weidman has been full-time with the BMSA since mid-2004. During his NRC career, he has staffed studies on a wide variety of topics related to mathematical, chemical, and materials sciences, laboratory assessment, risk analysis, and science and technology policy. His current focus is on building up the NRC's capabilities and portfolio related to all areas of analysis and computational science. He holds bachelor's degrees in mathematics and materials science from Northwestern University and M.S. and Ph.D. degrees in applied mathematics from the University of Virginia. Prior to joining the NRC, he had positions with General Electric, General Accident Insurance Company, Exxon Research and Engineering, and MRJ, Inc.

# C

# The Landscape of
# Digital Information Resources in
# Mathematics and Selected Other Fields

The following is a brief overview of some of the many information resources and tools currently available in mathematics and selected other fields, which offer some insight into the diverse ways that mathematics literature can be used.

## GENERAL BIBLIOGRAPHIC RESOURCES

Library information services have well-established conceptual schemas and database tools for handling the first five classes of bibliographic objects listed in Chapter 2 (documents, people, organizations, events, and subjects) and the most common relations between objects in these classes. These library services are exemplified by the following cross-disciplinary databases and portals:

- WorldCat[1]—a union catalog that itemizes the collections of 72,000 libraries in 170 countries and territories that participate in the On-line Computer Library Center (OCLC) global cooperative;
- Library of Congress—index of books, both academic and nonacademic;
- SciVerse Scopus—index of abstracts and citations for journal articles[2];

---

[1] OCLC, WorldCat, http://www.worldcat.org/, accessed January 16, 2014.
[2] Elsevier, "Scopus," http://www.info.sciverse.com/scopus, accessed January 16, 2014.

- Web of Science[3]—index of abstracts and citations for journal articles;
- Google Scholar[4]—a search engine for research literature capable of examining full text of articles (not just metadata and abstracts), ranking returns by citation counts and other criteria, and providing links to related papers and accessible versions;
- Scopus[5]—a bibliographic data service covering all academic fields, offering citation analysis tools, owned by Elsevier;
- Web of Science[6]—a bibliographic data service covering all academic fields, offering citation analysis tools, owned by Thompson Reuters; and
- Microsoft Academic Search[7]—a relatively new, free search engine for academic papers and resources, with the capability to identify papers, authors, conferences, journals, and organizations as first class objects; display relations between these objects; and the displays of "citation in context" with snippets from citing documents.

Larger, more loosely defined data structures and services use methods of massive data analysis (NRC, 2013) for search and discovery on the vastly larger scale of the World Wide Web. These services have become essential tools for information retrieval in mathematics as in every other field. They include the following:

- Google Web Search,[8]
- Bing,[9]
- Google Scholar[10] (an index of an unknown and not easily estimated number of academic books and articles), and
- Microsoft Academic Search[11] (an index of 48 million publications and more than 20 million authors across a variety of domains with updates added each week).

---

[3] Thomson Reuters, "Web of Science Core Collection," http://thomsonreuters.com/web-of-science/, accessed January 16, 2014.

[4] Google Scholar, http://scholar.google.com/, accessed January 16, 2014.

[5] Elsevier, Scopus, http://www.scopus.com/home.url, accessed January 16, 2014.

[6] Thomson Reuters, "Web of Science," http://thomsonreuters.com/products_services/science/science_products/a-z/web_of_science/, accessed January 16, 2014.

[7] Microsoft Academic Search, http://academic.research.microsoft.com/, accessed January 16, 2014.

[8] Google, https://www.google.com/, accessed January 16, 2014.

[9] "Bing," *Wikipedia*, last modified January 9, 2014, http://en.wikipedia.org/wiki/Bing.

[10] Google Scholar, http://scholar.google.com/, accessed January 16, 2014.

[11] "Microsoft Academic Search," *Wikipedia*, last modified January 12, 2014, http://en.wikipedia.org/wiki/Microsoft_Academic_Search.

Other, more specialized indexes provide essential Web services to participating partners. These services provide data that are consumed to varying extents in machine processing by the above services in preparation of data for display to human users. These indexes include the following:

- CrossRef[12] index of Digital Object Identifiers (DOIs),[13] available only to participating publishers; and
- ORCID (Open Researcher and Contributor ID) index of nonproprietary alphanumeric codes that uniquely identify academic authors with annual open data dumps.

## RESOURCES FOR THE MATHEMATICAL SCIENCES

### Specialized Mathematical Databases

Specialized mathematical databases are examples of "bottom up" attempts by the mathematics community to create relatively open, accessible databases of mathematical facts. There are many specialized databases of formal information that are of interest to specific communities, such as those described below.

- *On-Line Encyclopedia of Integer Sequences (OEIS)*[14]—This searchable database of integer sequences provides a brief description for each sequence, including how that sequence is defined and how it arises in various contexts, and related formulas, generating functions, code, links, and references (Sloan, 1973; Sloan and Plouffe, 1995). This resource is extremely valuable for researchers in number theory and combinatorics, where sequences arise naturally. It is very useful for a researcher encountering an unfamiliar sequence to check quickly if this sequence has been encountered before and, if so, what is known about it. OEIS has an active user community, which it relies on heavily for user contributions. It is licensed under the Creative Commons Attribution Non-Commercial 3.0 license.[15]
- *EZFace interface for evaluation of Euler sums*[16]—This specialized computational tool provides for the evaluation of multiple Euler sums, also known as multiple zeta values. Multiple zeta values are

---

functions of a finite sequence of positive integers and are known to satisfy a myriad of tricky identities. They can sometimes be reduced to polynomial functions of evaluations of the Riemann zeta function at integer values. This tool helps researchers who may encounter such sums to evaluate them using known reduction algorithms.

— *Distributome: An Interactive Web-based Resource for Probability Distributions*[17]—This is an open-source, open content-development project for exploring, discovering, navigating, learning, and computationally utilizing diverse probability distributions. Probability distributions are highly structured mathematical objects with fairly universal features, depending on the space over which a given probability distribution is defined (discrete, continuous, univariate, multivariate, Euclidian, non-Euclidean, etc.), such as a probability mass or density function, distribution function, quantile function, probability and moment generating function, etc. The interactive Distributome graphical user Navigator and the Distributome-Editor provide the following core functions:

o Visually traverse the space of all well-defined (named) distributions;
o Explore the relations between different distributions;
o Distribution search by keyword, property, and type;
o Obtain qualitative (e.g., analytic form of density function) and quantitative (e.g., critical and probability values) information about each distribution;
o Discover references and additional distribution resources; and
o Revise, add, and edit the properties, interrelations, and meta-data for various distributions.

Complete Java source code is available with the LGPL license.

• *Modular Forms Database*[18]—This database consists of tables related to modular forms, elliptic curves, and abelian varieties, which are specialized data of interest to number theorists.

• *Multiple Zeta Value Data Mine*[19]—These pages contain tables with multiple zeta values and Euler sums to allow people to look for relations, systematics, and patterns. They are expressed in terms of a basis.

---

[17] Distributome, http://www.distributome.org/, accessed January 16, 2014.

[18] William A. Stein, The Modular Forms Database, http://modular.math.washington.edu/Tables/, accessed January 16, 2014.

[19] Multiple Zeta Value Data Mine, http://www.nikhef.nl/~form/datamine/datamine.html, accessed January 16, 2014.

- *NIST Digital Library of Mathematical Functions (DLMF)*[20]—This is the Web version of the authoritative 1,046-page *Handbook of Mathematical Functions with Formulas, Graphs, and Mathematical Tables* (Abramowitz and Stegun, 1972). The DLMF has been constructed specifically for effective Web usage and contains features unique to Web presentation. The webpages contain many active links, for example, to the definitions of symbols within the DLMF, and to external sources of reviews, full texts of articles, and items of mathematical software. Advanced capabilities have been developed at the National Institute of Standards and Technology for the DLMF and also as part of a larger research effort intended to promote the use of the Web as a tool for doing mathematics. Among these capabilities are the following: a facility to allow users to download LaTeX and MathML encodings of every formula into document processors and software packages; a search engine that allows users to locate formulas based on queries expressed in mathematical notation; and user-manipulatable three-dimensional color graphics.
- *Information on Enumerative Combinatorics*[21]—This website contains a number of supplements to the two-volume textbook *Enumerative Combinatorics*,[22] including a *Catalan Addendum*, a 94-page PDF listing 204 combinatorial interpretations of the sequence of Catalan numbers. This site structures and curates the information and makes it available in machine-readable formats to allow various means of searching, browsing, and reuse.
- *Wolfram Functions Site*[23]—This website provides a substantial collection of formulas and graphics about mathematical functions. The information is fragmented into small packages (which makes it difficult to browse) and does not include references to original sources, and it is available only in proprietary formats (Mathematica® Notebook and PDF).

Currently, there is no unified way to exchange information between these specialized databases, and it is not clear that there are any incentives to make these databases talk to each other. Libraries have approached the interoperability issues at multiple levels. The highest-level and simplest ap-proach is the Open Archives Initiative (OAI), which provides for metadata exchange and federated search. The Protocol for Metadata Harvesting (OAI-PMH) enables spiders to gather up the cataloging information from multiple websites and then build a central search engine. The best known such service is OAISTER, now run by OCLC, which provides a search of more than 25 million records contributed by more than 1,100 institutions. For example, a search for a map of Polynesia from the 19th century turns up an 1839 map from the U.S. Hydrographic Office in the Harvard Map Collection (corrected to 1872). Entries in OAISTER typically have detailed but conventional library cataloging and refer to whole documents or objects.

---

[20] NIST Digital Library of Mathematical Functions, 2013, http://dlmf.nist.gov/.

[21] Information on Enumerative Combinatorics, http://www-math.mit.edu/~rstan/ec/, accessed January 16, 2014.

[22] Stanley, R.P., *Enumerative Combinatorics*, Cambridge University Press, Cambridge, Volume 1 (2nd edition, 2011) and Volume 2 (2001).

[23] Wolfram Research, Inc., The Wolfram Functions Site, http://functions.wolfram.com/, accessed on January 16, 2014.

More detailed interoperability methods include the linked open data movement, which tries to connect individual pieces of data using RDF (resource description format). RDF entries name two items and a relation between them, and are thus called "triples." So, to take an example from "data.gov.uk": the triple "John works for Ordnance Survey" would look something like:

http://www.johngoodwin.me.uk/me →
    http://www.intelleo.eu/ontologies/organization/ns/worksFor →
        http://data.ordnancesurvey.co.uk/id/ordnancesurvey
(John Goodwin, http://data.gov.uk/blog/what-is-linked-data)

In this example, the triple contains two items which identify John Goodwin and the Ordnance Survey, and a link between them labeled with "works for" as a relational concept. In this case, URLs are used for each item, with the relation taken from an ontology of organizational relations defined by a European project on learning. Other relations are defined by groups like the Dublin Core Metadata Initiative, which has cataloging-type relations such as publication date, author, and so on. The ontology for music (musicontology.com) describes relations such as conductor or artist.

Linked data are an example of the general concept of the Semantic Web introduced by Tim Berners-Lee and are in use in some very large organizations such as the British Museum. In general, these cooperative catalogs are based on volunteer contributions and run by some kind of nonprofit group.

### Bibliographic Resources

There are currently many bibliographic resources available within the mathematical sciences as well as the larger scientific community. Some examples of these mathematical bibliographic resources include the following:

- *MathSciNet*[24] is the online interface to the database of *Mathematical Reviews* maintained by the American Mathematical Society (AMS) since 1940. It is a carefully maintained and easily searchable database of reviews, abstracts, and bibliographic information for much of the mathematical sciences literature. More than 100,000 new items are added each year, most of them classified according to the Mathematics Subject Classification (MSC). Authors are uniquely identified, enabling a search for publications by individual author rather than by name string. Expert reviewers are selected by a staff of professional mathematicians to write reviews of the current published literature; more than 80,000 reviews are added to the database each year. MathSciNet contains more than 2.8 million items and more than 1.6 million direct links to original articles. Bibliographic data from retro-digitized articles dates back to the early 1800s. Reference lists are collected and matched internally from approximately 500 journals, and citation data for journals, authors, articles, and reviews is provided. This Web of citations allows users to track the history and influence of research publications in the mathematical sciences. MathSciNet is a major revenue generator for AMS, for which reason the database contents are closely protected by copyright and licensing.

- *Zentralblatt MATH (zbMATH)*[25] is a thorough and long-running abstracting and reviewing service in pure and applied mathematics. The zbMATH database contains more than 3 million bibliographic entries with reviews or abstracts drawn from more than 3,500 journals and 1,100 serials and covers the period from 1826 to the present. Reviews are written by more than 10,000 international experts, and the entries are classified according to the MSC scheme (MSC 2010). zbMATH covers published and refereed articles, books, and conferences as well as other publication formats (CD-ROM, DVD, videotapes, Web documents). Within current electronic library activities retrospective data of journals are made available even prior to 1868. The bibliographic information and links to the full text are stored within zbMATH if available. The current number of new items added to zbMATH is about 120,000 per year. More than 50 percent of the items core areas are independent reviews by experts, the remainder are abstracts and summaries of comparable quality. zbMATH is run jointly by the European Mathematical Society, FIZ Karlsruhe, and Springer-

Verlag. zbMATH is a subscription service but allows nonsubscribers to ask queries and access the zbMATH author profile pages,[26] which are freely accessible.

- *Ulf Rehmann's DML page*[27] lists retro-digitized mathematics links to nearly 5,000 digitized books and to nearly 600 digitized journals/seminars. This is a major resource for discovering information that has already been digitized. The webpage also lists more than 2,800 journals that have been digitized whole or in part and notes whether they are free or require a paid subscription.

- *AMS Digital Mathematics Registry*[28] provides centralized access to certain collections of digitized publications in the mathematical sciences. The registry is primarily focused on older material from journals and journal-like book series that originally appeared in print but now are available in digital form.

- *AMS eBooks*[29] includes retrospective digitization of *Contemporary Mathematics* back to the beginning of the series in 1980.

- *European Digital Mathematics Library (EuDML)*[30] makes a significant portion of European mathematics literature available online: more than 200,000 publications, in the form of an enduring digital collection, developed and maintained by a network of institutions. A unified metadata schema was developed and adopted by all providers. The library offers a number of features including the following:

  — Metadata search over the entire corpus,
  — Reference and citation lists,
  — Capability for users to make lists and annotations,
  — An API for metadata search over the entire corpus, and
  — Some capability for formula search.

### Encyclopedia Resources

Some encyclopedia resources are listed below.

---

- *MacTutor History of Mathematics Archive*[31] contains biographies of several thousand historical and contemporary mathematicians as well as an index of famous curves and histories of various mathematical topics. The full text is freely available, with no formal copyright or licensing restrictions.
- *On-Line Encyclopedia of Integer Sequences* is described above in the "Specialized Mathematical Databases" discussion.
- *Mathematics Genealogy Project*[32] aims to list all individuals who have received a doctorate in mathematics, providing the following information:

  — The complete name of the degree recipient,
  — The name of the university that awarded the degree,
  — The year in which the degree was awarded,
  — The complete title of the dissertation, and
  — The complete name(s) of the advisor(s).

  The Mathematics Genealogy Project contains more than 170,000 records. Individual pages can be freely copied without explicit licensing or copyright restriction, but data are not made available in bulk, and there is no API.
- *Wolfram's MathWorld*[33] is a comprehensive and interactive encyclopedia of mathematical equations, terms, derivations, and more, for students, educators, math enthusiasts, and researchers.
- *Wikipedia*[34] is perhaps the most well known of all online encyclopedia resources. It also houses a wide array of mathematical content, generally very useful as a first place to look for the definition of a mathematical concept. *Wikipedia* uses Creative Commons Attribution-ShareAlike (CC-BY-SA) license.
- *Encyclopedia of Mathematics*[35] is an open access wiki that includes original articles from the online *Encyclopedia of Mathematics* (2002) as well as user-added articles, totaling more than 8,000 entries and nearly 50,000 notions in mathematics. Springer, in cooperation with the European Mathematical Society, has made the content of this encyclopedia freely open to the public. The original articles from the *Encyclopedia of Mathematics* remain copyrighted to Springer, but any new articles added and any changes made to existing articles within encyclopediaofmath.org will come under the CC-BY-SA license. An editorial board, under the management of the European Mathematical Society, monitors any changes to articles and has full scientific authority over alterations and deletions. This wiki is a MediaWiki that uses the MathJax extension, making it possible to insert mathematical equations in TeX and LaTeX.

---

[31] University of St Andrews, Scotland, The MacTutor History of Mathematics Archive, October 2013, http://www-history.mcs.st-and.ac.uk/.

[32] North Dakota State University, Mathematics Genealogy Project, http://genealogy.math.ndsu.nodak.edu/, accessed January 16, 2014.

[33] Wolfram MathWorld, http://mathworld.wolfram.com/, accessed January 16, 2014.

[34] *Wikipedia*, http://www.wikipedia.org/, accessed January 16, 2014.

[35] *Encyclopedia of Mathematics*, http://www.encyclopediaofmath.org/index.php/, accessed January 16, 2014.

- *The Stacks Project*[36] website is an open source textbook and reference work on algebraic stacks and the algebraic geometry needed to define them. *The Stacks Project* aims to build up enough basic algebraic geometry to serve as foundations for algebraic stacks.

### Specialized Mathematical Resources

Several specialized mathematical resources are available to the mathematics community. Some of these resources include the following:

- *MathOverflow*[37] is an online resource that allows users to ask and answer research-level mathematics questions such as arise when writing or reading articles or graduate-level books. Users gain writing authority on the site by building up reputation points. Mathematics display support is provided with MathJaX from LaTeX source. MathOverflow runs on Stack Exchange, the hosted service that provides the same software as the popular programming Q&A site Stack Overflow. The hosting cost is paid from the research funds of Ravi Vakil at Stanford University. User-contributed content is licensed under Creative Commons Attribution-Share Alike.
- *Wolfram|Alpha*[38] is a "computational knowledge engine" developed as an online service by Wolfram Research. It answers factual queries by computation of the answer from an internal database of mathematical and factual data acquired from diverse data sources. Both free and premium services are available. Underlying software combines natural language processing of queries with symbolic computation using Mathematica. Numerous mathematical concepts, such as sequences, functions, and probability distributions are recognized and displayed in ways that respect their mathematical structure.

---

[36] The Stacks Project, http://stacks.math.columbia.edu/, accessed January 16, 2014.

[37] MathOverflow, http://mathoverflow.net/, accessed January 16, 2014.

[38] Wolfram|Alpha, search engine, http://www.wolframalpha.com/, accessed January 16, 2014.

- *Selected Papers Network*[39,40] is a free, open-source project aimed at improving the way people find, read, and share academic papers. This project is not a website with a system for reviewing, evaluating, rewarding, etc. Rather, it is an environment that makes it easy to build one's own systems, which allows for more flexibility when needed.
- *Tricki*[41] is a Wiki-style site intended to develop into a large store of useful mathematical problem-solving techniques. Some of these techniques are very general, and others concern particular subareas of mathematics spanning all levels of experience. This project is largely inactive now after failing to acquire critical mass of users.

## SELECTED RELATED EFFORTS

Many disciplines have ongoing efforts that aim to bring diverse discipline-specific information together, and many of these hold valuable lessons for the mathematics community. The following are a few illustrations of such efforts.

- *Digital Library Federation Aquifer (DLF Aquifer)*[42] promotes effective use of distributed digital library content for teaching, learning, and research in the area of American culture and life. It supports scholarly discovery and access by developing schemas, protocols, and communities of practice to make digital content available to scholars and students where they do their work, and by developing the best possible systems for finding, identifying, and using digital resources in context.
- *Project Bamboo*[43] is a partnership of 10 research universities building shared infrastructure for humanities research. Led by the University of California, Berkeley, one of the goals of this project is to design research environments where scholars may discover, analyze, and curate digital texts across the 450 years of print culture in English from 1473 until 1923, along with the texts from the Classical world upon which that print culture is based.
- *Research Papers in Economics*[44] is a collaborative effort of hundreds of volunteers in 75 countries to enhance the dissemination

---

[39] SelectedPapers, https://selectedpapers.net/, accessed January 16, 2014.
[40] "The Selected-Papers Network," *Gower's Weblog*, June 16, 2013, http://gowers.wordpress.com/2013/06/16/the-selected-papers-network/.
[41] Tricki, http://www.tricki.org/, accessed January 16, 2014.
[42] DLF Aquifer, http://old.diglib.org/aquifer/ (no longer maintained as of June 2010).
[43] Project Bamboo, http://www.projectbamboo.org/, accessed January 16, 2014.
[44] RePEc, http://repec.org/, accessed January 16, 2014.

of research in economics and related sciences. The heart of the project is a decentralized bibliographic database of working papers, journal articles, books, books chapters, and software components, all maintained by volunteers. The collected data are then used in various services.

- *Digital Library of Chemistry Education*[45] provides an integrated guide to chemistry textbooks and allows both students and educators to explore chemistry. The ChemEd DL repository can be searched for resource groups within particular domains of chemistry, such as organic or physical. Resource groups relate to specific topics, such as bonding or kinetics, and are associated with specific elements. ChemEd also allows users to search by topics and look up definitions of terms. The provided glossary is extensive.

- *Digital Library of Biochemistry and Molecular Biology.* BioMoleculesAlive.org is a collection of digital resources sponsored by the American Society for Biochemistry and Molecular Biology. It is part of a larger effort called the BioSciEdNet (BEN) initiative.[46] The collection includes resources in five areas: software, visual resources, curriculum resources, reviews, and articles from the *Biochemistry and Molecular Biology Education* journal. Efforts on the Web interface, database design, and tools and guidelines for submission to BioMoleculesAlive.org began in 2003 and are still ongoing.

- *Astrophysics Data Service (ADS).*[47] Also known as the Digital Library for Physics and Astronomy, this library is maintained by the Smithsonian Astrophysical Observatory, working with NASA and the community of astronomers and astrophysicists, and links to more than 10 million papers in astronomy and related areas. An unusual aspect of this system is that it not only catalogs papers, but also tries to link papers to the astronomical objects to which they refer. A user can see papers that refer to a specific star or galaxy, via volunteer tagging of all papers with star catalog entries. NASA provides the base funding for ADS.

- *U.S. Virtual Astronomical Observatory.* Astronomers have access to a variety of sky images, including some interfaces designed for the general public, such as Google Sky or the WorldWide Telescope (Microsoft). Digital imagery exists at multiple wavelengths, includ-

---

ing the Sloan Digital Sky Survey showing visible light, the Two Micron All Sky Survey (2MASS), the Chandra X-ray Observatory, and so on. These databases are unified via the Virtual Observatory program, including the Euro-VO in Europe and others. Funding for the U.S. Virtual Observatory has been provided by NSF and NASA, but the organization is attempting to find a new support model.

- *National Center for Biotechnology Information (NCBI)*. The National Library of Medicine maintains many important bio-medical data resources. Full articles are stored in PubMed Central,[48] which receives medical articles deposited by authors working on NIH-funded research (after an embargo period). It currently contains 2.8 million articles. More detailed data is stored in several specific resources such as GenBank or OMIM (Online Mendelian Inheritance in Man). NCBI also provides software tools such as BLAST (Basic Local Alignment Search Tool). All these resources are funded by NIH in the United States. A number of other organizations support tools for molecular biology. For example, EMBL (the European Molecular Biology Laboratory) provides bio-informatic services including tools for sequencing, structural analysis, microscopy, and so on. Other groups that provide molecular biology tools include the Wellcome Trust Sanger Institute, the Craig Venter Institute, and commercial suppliers. EMBL is funded by a consortium of nations not exactly overlapping the European Union, but close. The Wellcome Trust is endowed under the will of Sir Henry Wellcome, the Venter Institute is supported by J. Craig Venter and others, and so on.

- *Digital Public Library of America*. Numerous libraries have provided cataloging information to the Digital Public Library of America, which provides links to more than 2 million items in its member libraries. There are currently more than 400 participating libraries, including the many libraries aggregated by state library systems. The organization is a cooperative of its members, aggregated into "hubs."

- *Chemical Abstracts*.[49] The American Chemical Society operates one of the largest and oldest scientific information services. Chemical Abstracts Service indexes and abstracts the chemical literature and maintains an authority file of chemical compounds, with more than 70 million entries. It also keeps track of reactions, suppliers, and other chemical information resources. *Chemical Abstracts* dates back to 1907 and is one of the most exhaustive services, with a history of seeking out all important chemical information, wherever it is published. In its early years, it was largely supported by major chemical companies, but for decades has been funded by users, typically university libraries or industrial organizations in chemistry, chemical engineering, biomedicine, or related areas.

---

[48] National Library of Medicine, PubMed, http://www.ncbi.nlm.nih.gov/pubmed, accessed January 16, 2014.

[49] American Chemical Society, Chemical Abstracts Service (CAS), http://www.cas.org/, accessed January 16, 2014.

- *Internet Public Library.* The Internet Public Library is a resource to provide answers to questions, particularly questions from students and educators. It also maintains some information collections. Originally operated at the University of Michigan with funding from the W. K. Kellogg Foundation, it is now run by Drexel University with support from a group of about 20 universities.